\documentclass{article}
\makeatletter
\newcommand{\@affiliation}{}
\makeatother

\usepackage{arxiv}
\usepackage{comment}
\usepackage[utf8]{inputenc} 
\usepackage[T1]{fontenc}    
\usepackage{hyperref}       
\usepackage{url}            
\usepackage{booktabs}       
\usepackage{amsfonts}       
\usepackage{microtype}      
\usepackage{xcolor}
\usepackage{lipsum}         
\usepackage{graphicx}
\usepackage{doi}
\usepackage[
backend=biber, style=numeric, sortcites=true, citestyle=numeric-comp, giveninits=true, sorting=nyt, doi=false, isbn=false, url=false, maxcitenames=3, maxbibnames=10, minalphanames=3, block=none, date=year
]{biblatex}
\newenvironment{proof}{\paragraph{Proof:}}{\hfill$\square$}

\usepackage{amsmath,amsfonts}
\usepackage{amssymb}
\usepackage{nicefrac}
\usepackage{mathtools}
\numberwithin{equation}{section}
\usepackage{subcaption}
\usepackage{cleveref}  
\usepackage{todonotes}
\usepackage{xfrac}
\usepackage{csvsimple}

\usepackage{gensymb}

\newcommand{\ds}{\displaystyle}

\newcommand{\cI}{\mathcal{I}}

\newcommand{\cO}{\mathcal{O}}

\newcommand{\cR}{\mathcal{R}}

\newcommand{\eps}{\varepsilon}
\renewcommand{\phi}{\varphi}

\AtBeginDocument{

}

\makeatletter
\newcommand{\oset}[3][0ex]{%
  \mathrel{\mathop{#3}\limits^{
    \vbox to#1{\kern-2\ex@
    \hbox{$\scriptstyle#2$}\vss}}}}
\makeatother

\DeclarePairedDelimiterX{\norm}[1]{\lVert}{\rVert}{
\ifblank{#1}{\:\cdot\:}{#1}
}

\makeatletter
\newcommand{\opnorm}{\@ifstar\@opnorms\@opnorm}
\newcommand{\@opnorms}[1]{%
  \left|\mkern-1.5mu\left|\mkern-1.5mu\left|
   #1
  \right|\mkern-1.5mu\right|\mkern-1.5mu\right|
}
\newcommand{\@opnorm}[2][]{%
  \mathopen{#1|\mkern-1.5mu#1|\mkern-1.5mu#1|}
  #2
  \mathclose{#1|\mkern-1.5mu#1|\mkern-1.5mu#1|}
}
\makeatother

\newcommand{\DoubWellPot}   {W}

\newcommand{\sphere}{{{\mathcal S}^{d-1}}}

\newcommand{\vecstyle}[1]{{\bf{#1}}}
\newcommand{\vecu}{\vecstyle{u}}
\newcommand{\vece}{\vecstyle{e}}

\newcommand{\vecd}{\vecstyle{d}}

\newcommand{\vecx}{\vecstyle{x}}

\newcommand{\vecn}{\vecstyle{n}}

\newcommand{\vecv}{\vecstyle{v}}
\newcommand{\setstyle}[1]{{\mathbb #1}}
\newcommand{\R}{\setstyle{R}}

\def\setR {\setstyle{R}}

\newcommand{\visc}         {\nu}
\newcommand{\Dt}           {\Delta t}
\newcommand{\domain}       {\Omega}
\newcommand{\xmin}         {x_{min}}
\newcommand{\xmax}         {x_{max}}
\newcommand{\ymin}         {y_{min}}
\newcommand{\ymax}         {y_{max}}
\newcommand{\nCells}       {N}
\newcommand{\Nx}           {N_x}
\newcommand{\Ny}           {N_y}
\newcommand{\vol}[1]       {\mathrm{vol}\left( #1 \right)}
\newcommand{\domainCell}   {\domain_{i,j}}
\newcommand{\domainCellx}  {\domain_{i-\frac{1}{2},j}}
\newcommand{\domainCelly}  {\domain_{i,j-\frac{1}{2}}}
\newcommand{\Dx}           {\triangle x}
\newcommand{\Dy}           {\triangle y}
\newcommand{\Vecx}         {\mathbf{x}}
\newcommand{\xindex}       {i}
\newcommand{\yindex}       {j}
\newcommand{\phas}         {c}
\newcommand{\pres}         {p}
\newcommand{\velv}         {\mathbf{u}}
\newcommand{\velx}         {u_1}
\newcommand{\vely}         {u_2}
\newcommand{\grad}         {\nabla}
\newcommand{\lap}          {\Delta}

\newcommand{\vecmu}        {\mathbf{v}}
\newcommand{\vecmux}       {v_1}
\newcommand{\vecmuy}       {v_2}
\newcommand{\chempot}      {\mu}
\newcommand{\ordp}         {\omega}
\newcommand{\OperatorORDP} {\cI_{\ordp}}
\newcommand{\oot}          {\frac{1}{2}}
\newcommand{\numepsilon}   {{\tilde{\varepsilon}}}

\newcommand{\transpose}    {^\top}
\newcommand{\pfrac}[2]     {\frac{\partial #1}{\partial #2}}
\newcommand{\pppfrac}[2]   {\frac{\partial^3 #1}{\partial #2^3}}
\newcommand{\pppxyyfrac}[3]{\frac{\partial^3 #1}{\partial #2 \partial #3^2}}
\newcommand{\pppxxyfrac}[3]{\frac{\partial^3 #1}{\partial #2^2 \partial #3}}

\renewcommand{\div}        {\nabla \cdot}

\DeclareMathAlphabet{\mathmybb}{U}{bbold}{m}{n}
\newcommand{\1}{\mathmybb{1}}

 \newtheorem{thm}{Theorem}[section]

 \newtheorem{Lemma}[thm]{Lemma}
 \newtheorem{Definition}[thm]{Definition}
 \newtheorem{Proposition}[thm]{Proposition}
 \newtheorem{rem}[thm]{Remark}

\title{Justification of a relaxation approximation for the Navier--Stokes--Cahn--Hilliard system}

\newif\ifuniqueAffiliation
\uniqueAffiliationfalse

\ifuniqueAffiliation 
\author{ \href{https://orcid.org/0000-0000-0000-0000}{\includegraphics[scale=0.06]{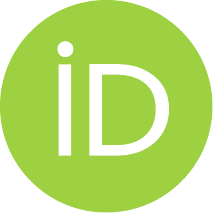}\hspace{1mm}David S.~Hippocampus}\thanks{Use footnote for providing further
		information about author (webpage, alternative
		address)---\emph{not} for acknowledging funding agencies.} \\
	Department of Computer Science\\
	Cranberry-Lemon University\\
	Pittsburgh, PA 15213 \\
	\texttt{hippo@cs.cranberry-lemon.edu} \\
	\And
	\href{https://orcid.org/0000-0000-0000-0000}{\includegraphics[scale=0.06]{orcid.pdf}\hspace{1mm}Elias D.~Striatum} \\
	Department of Electrical Engineering\\
	Mount-Sheikh University\\
	Santa Narimana, Levand \\
	\texttt{stariate@ee.mount-sheikh.edu} \\
}
\else
\usepackage{authblk}

\newbox{\orcid}\sbox{\orcid}{\includegraphics[scale=0.06]{orcid.pdf}} 
\author{%
	Jan Giesselmann\thanks{Numerical Analysis and Scientific Computing,
Department of Mathematics,
Technical University of Darmstadt, email: \texttt{Jan.Giesselmann@tu-darmstadt.de}}
~\&
Jens Keim\thanks{Institute of Applied Analysis and Numerical Simulation \& Institute of Aerodynamics and Gas Dynamics, University of Stuttgart, email: \texttt{Jens.Keim@mathematik.uni-stuttgart.de}}%
~ \& 
	Fabio Leotta\thanks{Numerical Analysis and Scientific Computing,
Department of Mathematics,
Technical University of Darmstadt, email: \texttt{Fabio.Leotta@tu-darmstadt.de}}   
~ \& 
Christian Rohde\thanks{Institute of Applied Analysis and Numerical Simulation  \& Stuttgart Center for Simulation Science, University of Stuttgart, email: \texttt{Christian.Rohde@mathematik.uni-stuttgart.de}}}%

\fi

\renewcommand{\headeright}{Relaxation approximation of the Navier-Stokes--Cahn--Hilliard system}
\renewcommand{\undertitle}{}
\renewcommand{\shorttitle}{}

\hypersetup{
pdftitle={A template for the arxiv style},
pdfsubject={q-bio.NC, q-bio.QM},
pdfauthor={David S.~Hippocampus, Elias D.~Striatum},
pdfkeywords={First keyword, Second keyword, More},
}

\begin{document}
\maketitle

\begin{abstract} The  Navier-Stokes-Cahn--Hilliard (NSCH) system governs the diffuse-interface dynamics of two incompressible and immiscible fluids. We consider a relaxation approximation  of the 
NSCH system  that is composed by  a system of first-order hyperbolic balance laws and second-order elliptic operators. We prove first that 
the solutions of an initial boundary value problem for the   approximation recover the limiting 
NSCH system for vanishing relaxation parameters. To cope with the singular limit we exploit the fact  that the  approximate solutions dissipate an almost quadratic    energy,  and employ the relative 
entropy-framework. In the  second part of the work we provide numerical evidence for the analytical results, even in flow regimes not covered by the  assumptions needed for the theoretical results. 
Using a novel marker-and-cell conservative finite-difference approach for both  the approximation and the limit system, we are able to compute physically relevant interfacial flow
problems including Ostwald ripening and high-velocity flow.
\end{abstract}




\section{Introduction}




We consider  the  viscous  motion of two incompressible, immiscible fluids at constant temperature. 
For the  modelling of  such two-phase flows there are two  major  approaches: the sharp-interface  approach and the diffuse-interface  approach. We consider  the   diffuse-interface approach such that  interfaces are represented as steep but continuous transition zones, with interfacial width
controlled by a small parameter  $\gamma >0$. Thereby we focus on the classical 
Navier--Stokes--Cahn--Hilliard (NSCH) system as it has been proposed  in \cite{hohenberg1977theory} to describe two fluids with matching density, see also \cite{Lowengrub1998}. This model  is thermodynamically consistent in the sense that it dissipates a free energy functional  
that is composed of kinetic energy and the Van-der-Waals energy. It has proven to govern  flows involving topological changes like the coalescence and breakup of bubbles and droplets.
In the last decades the NSCH  system has been extended  to 
multiple physical settings; we only mention here 
extensions to cope  with inhomogeneous   densities of the two fluids, see, e.g.,~\cite{ABG,Boyer01,Eikelder}.  For the well-posedness of the resulting  initial boundary problems for the NSCH system  in the  framework of weak and strong solutions 
we  refer to  \cite{Abels09,Giorgini22,Giorgini19}; without any attempt to cover the 
rich literature. \\
Despite significant progress  (see, e.g.,~\cite{brunk2023second,chen2016efficient, diegel2017convergence,Feng}), the numerics for the  NSCH system still poses  major challenges, most notably in the context of   high-velocity regimes.   These difficulties are related to 
 the nonlocality induced by the  incompressibility  constraint  and the fourth-order operators in the Cahn--Hilliard evolution. As a remedy in the framework of the incompressible Navier--Stokes equations, approximations based on the artificial compressibility approach   have been suggested, see \cite{Temam1969I}.  Notably, this approximation  involves a nonlinear  first-order hyperbolic transport 
operator (see \cite{HRYZ25} for a version with a linear transport operator but nonlinear zeroth-order terms) such that it is accessible for robust numerical methods from the field of hyperbolic balance laws. An artificial compressibility  approximation for the NSCH system can be found in the work \cite{GiorginiTemamVu} that relaxes the divergence constraint by an evolution equation for the pressure but 
keeps the Cahn--Hilliard evolution unchanged. Lower-order hyperbolic approximations or approximations  based  on the hyperbolic  wave operator of the  Cahn--Hilliard equations have been suggested in  e.g.~\cite{Dhaouadi24,Grasselli05}. The work of \cite{Dhaouadi24} relies on a two-stage approximation procedure. First,  the Cahn--Hilliard equation is embedded into a system composed of  a conservation law for the phase-field variable and a third-order system of evolution  equations for the  Cahn--Hilliard flux with frictional damping. The Cahn--Hilliard equation is recovered in the high-friction limit.  Next, the third-order operators are relaxed using a screened Poisson equation. This results in an  approximate system  of first-order evolution equations,  similar to previously found approximations of  the  Navier--Stokes--Korteweg models, see  \cite{GiesselmannLattanzioTzavaras2017,keim2020relaxation,rohde2010local}.\\
Based on  a combination of the  artificial compressibility approach for the Navier--Stokes equations, the friction-type approximation, and relaxation of the  Cahn--Hilliard equations,   a lower-order relaxation 
approximation of the entire Navier--Stokes--Cahn--Hilliard model has been 
suggested in \cite{KKR25}.  This relaxation approximation comes with three a-priori independent approximation parameters such that the NSCH system is formally recovered if they tend to zero.
We note that the relaxation approximation is thermodynamically consistent, i.e.,  there is a free energy functional that  is non-increasing in-time for 
smooth solutions. For vanishing approximation parameters, the free energy formally collapses to the free energy of the NSCH system.
Under reasonable assumptions it has been shown that the involved first-order operator in the lower-order approximate system is hyperbolic. In fact, for this sub-system, the free energy acts as a   convex (mathematical) entropy. Due to the hyperbolic structure one  can  apply  also in this case numerical  methods from the field of systems of hyperbolic 
balance laws. First  results in one spatial  dimension, which have been presented in \cite{KKR25}, suggest that it is also possible to address in this way high-velocity flows.

In the present work, we provide a rigorous analytical justification of a slightly modified version of the relaxation approximation introduced in \cite{KKR25} as an approximation of the Navier--Stokes--Cahn--Hilliard system. 
The modification consists of additional dissipative terms, which can be used to absorb error terms in the relative-entropy inequality, without altering the hyperbolic structure of the first-order operator.
All involved systems are shortly reviewed in Section \ref{chap:p2:NSAC}.
To the best of our knowledge, this is the first result establishing convergence of a fully lower-order relaxation model with hyperbolic first order operator to the NSCH system in the regime of smooth solutions.
In addition, our result does not rely on viscosity, i.e., it trivially extends to the Euler--Cahn--Hilliard case.

For the approximate system, we employ a low-regularity entropy solution 
concept. Our main result, Theorem \ref{thm:high order convergence result}  in  Section \ref{sec:analysis}, proves
that we recover in the limit of vanishing approximation parameters the velocity and phase field  solving the initial boundary value problem for the  
NSCH system, as long as the latter  admits a sufficiently  regular  strong solution.
In the proof, we exploit the energetic structure of the relaxation approximation and the NSCH system. We utilize 
the relative-entropy method  (see \cite{DafermosBook,DiPerna1979})   and derive quantitative expressions for the error in terms of all three approximation parameters yielding explicit convergence rates in the regime of smooth solutions.
A key difficulty, especially in the setting of weak solutions for the relaxation system, arises from the non-convexity of the energy functional of the relaxation system, which prevents a direct application of standard relative-entropy arguments. 
Another key ingredient of the proof is the construction of an approximate solution to the relaxation system from a given smooth solution of the NSCH system. We use a nontrivial construction, since a straightforward choice would lead to reduced convergence rates, see Remark \ref{rem:reduced} for details.
\\
Section \ref{sec:numerics} is devoted to  the numerical justification of the relaxation approximation  by a novel 
numerical scheme. It is proposed in Section \ref{sec:num} as a  first-order time-integration scheme  for both the NSCH system and the relaxation system such that Chorin's projection method is recovered in the limit of vanishing approximation parameters.
In space, a staggered-mesh conservative  
finite-difference method is employed. More precisely, 
we introduce  a marker-and-cell (MAC)     finite difference scheme for both the target NSCH system and the relaxation approximation. Within this second-order approach, the pressure and the phase field are approximated at the centers of a Cartesian mesh volume whereas the velocity and the Cahn--Hilliard flux are evaluated on the mid points of the boundaries of the cell volumes.
We test the full discretization method on one- and two-dimensional settings. The numerical experiments  confirm the convergence rates
predicted  by Theorem \ref{thm:high order convergence result}. Moreover, we display results of numerical experiments featuring  Ostwald ripening and bubble
coalescence scenarios for the relaxation approximation. The latter setup with velocity fields of large magnitude clearly shows that our approach is robust for   high-speed flow fields . \\
Let us note  that our approach is different from the 
 MAC finite-difference method for the NSCH
system from \cite{LiShen20} which relies on the 
classical work \cite{HW65} for the incompressible Navier--Stokes equations. The method in \cite{LiShen20} builds on the order reduction by the scalar auxiliary variable idea. We follow \cite{Shu2006} and   discretize the fourth-order  operator directly.


We conclude the introduction with some remarks on the use of the relative-entropy method in related frameworks.  Closely related to the artificial-compressibility ansatz is the low-Mach number asymptotics, i.e., the passage from the compressible regime to incompressible  Navier--Stokes or Euler flow. Rigorous analysis for vanishing Mach number   based  on relative entropies  can be found in \cite{ALN24}. The high-friction limit that connects the damped  Navier-Stokes-Korteweg equations to the Cahn--Hilliard equation has been analyzed  via relative entropies for convex free energies in, e.g., \cite{CL20,GiesselmannLattanzioTzavaras2017} and  for non-convex cases in  \cite{CG25}. With regard to the  relaxation of the originally  third-order term in the friction-type approximation, we refer to
\cite{chaudhuri2025weakstronguniquenessrelaxationlimit}.  Finally, we mention the closely related works \cite{barthwal2025hyperbolicapproximationsclassdispersive,GiesselmannRanocha2025} to approximate third-order operators in dispersive transport.


\section{Phase-field modelling for  incompressible two-phase flow}\label{chap:p2:NSAC}
We present a short review of the NSCH system and the relaxation approximation together with the 
respective energy dissipation statements. 
\subsection{The Navier--Stokes--Cahn--Hilliard system}\label{subsec:NSCH}
For $d\in \{2,3\}$, let $\Omega \subset \setR^d$ be an open, bounded set with smooth boundary $\partial \Omega$. This trivially includes the  periodic torus case $\Omega= \mathbb{T}^d$ which is used in the analytical Section \ref{sec:analysis}. 
We denote the outer normal of $\partial \Omega$ by $\vecn \in \sphere$.\\
For time $T>0$, the classical Navier--Stokes--Cahn--Hilliard (NSCH) system for matched densities from \cite{hohenberg1977theory} is given 
by
\begin{equation}\label{eq:NSCH}
\begin{array}{rcl}
\nabla \cdot \mathbf{u} 
		&=& 0, \\[1.2ex]
		\mathbf{u}_t 
		+  (\mathbf{u} \cdot \nabla)\mathbf{u} + \nabla p 
		&=&  -c \nabla \mu(c) + \nu \Delta \mathbf{u}, \\[1.2ex]
		c_t 
		+  \nabla \cdot (c \mathbf{u}) 
		&=&  \nabla \cdot \big(\nabla \mu(c)\big)
	\end{array}\qquad 
	\text{ in $ \Omega_T:=(0,T) \times  \Omega   $}.
\end{equation}
The hydromechanical unknowns in \eqref{eq:NSCH} are the pressure (perturbation) $p=p(t,\vecx)\in \setR$ and the velocity field $\vecu = \vecu(t, \vecx)\in \setR^d$. For the Cahn--Hilliard part, 
we  have the phase-field variable 
$c=c(t, \vecx)\in \setR$, which acts as phase indicator. \\
In \eqref{eq:NSCH}, we denote  by $\nu \ge 0$ the  constant kinematic  viscosity parameter. Thus, inviscid Euler flow is included by choosing $\nu =0$.  The parameter $\gamma >0$ is  the phase-field parameter which is fixed throughout the paper.  
For the free energy function $W:\setR\to \setR$, we make the simple polynomial double-well choice
\begin{equation}\label{doublewell} 
{W} (c) = \dfrac{1}{4} (c^2-1)^2 \qquad   (c \in \setR).
\end{equation}
Having $W$ defined,  the chemical potential in \eqref{eq:NSCH} is given by $       \mu(c) :=  W'(c) - \gamma \Delta c$.\\
Note that we consider here a special version of the  
Cahn--Hilliard equation where  densities of both phases  are equal and  the mobility is supposed to be constant (all set to $1$).  A generalization of our  results to general mobilities  and in particular the analysis for  the case of different densities, see \cite{ABG}, is much more intricate.\\
The system \eqref{eq:NSCH} is complemented with
the initial conditions
\begin{equation}\label{eq:iniNSCH}
	\mathbf{u}(0,\cdot) =\vecu_0,\quad c(0,\cdot)=c_0
    \text{ in } \Omega.
\end{equation}
If $\Omega \neq \emptyset $ holds, a natural choice for boundary conditions is given  by 
\begin{equation}\label{eq:boundNSCH}
	\mathbf{u}=\mathbf{0},\quad\nabla c\cdot\mathbf{n}=0,\quad\nabla\mu(c)\cdot\mathbf{n}=0
  \text{ on } (0,T) \times \partial\Omega.
\end{equation}
%
The homogeneous  boundary condition    $\nabla c \cdot \vecn =0 $ ensures that the contact angle 
of a diffuse interface  at the boundary $\partial \Omega$ is fixed to  be $90^o$. 
In Section \ref{sec:analysis}, we are interested in classical solutions for the initial value problem \eqref{eq:NSCH}, \eqref{eq:iniNSCH} on $\Omega= {\mathbb{T}}^d$. By a classical solution
 we mean a  triple    $(p,\vecu,c): \bar\Omega_T \to \setR^{d+2}$ 
 with 
\begin{align*}
&  p \in C^{0}(\Omega_T), p(\cdot,t) \in C^1(\Omega),\\ 
&  \vecu \in C^1((0,T); C(\Omega)^d)  \cap C^0([0,T]; C^2(\Omega)),   \\                                
&c \in C^1((0,T); C(\Omega))    \cap  C([0,T]; C^4(\Omega)), 
\end{align*}
 such that all equations in \eqref{eq:NSCH}, in \eqref{eq:iniNSCH}, and as the case may be, in \eqref{eq:boundNSCH} hold pointwise.\\
For the wellposedness   of solutions for the initial boundary value problem \eqref{eq:NSCH}, \eqref{eq:iniNSCH}, \eqref{eq:boundNSCH} let us  refer to \cite{Giorgini22}. The authors  establish a 
global wellposedness result for a strong  solution concept in 2D and for weak solutions in 2D/3D.\\
Our analysis 
in Section \ref{sec:analysis} relies on energy techniques. 
For the NSCH system \eqref{eq:NSCH}, we define the energy expression
\begin{equation}\label{eq:energy NSCH}
     E[\vecu,c,\vecd] =  \frac12 {|\vecu|^2}  + W(c) + \frac{\gamma}{2} {|\vecd|^{2}}
     \quad \left(\vecu \in \setR^d, \, c\in \setR, \, \vecd \in \setR^d  \right),
\end{equation}
which is composed as a sum of the kinetic energy and   the classical van der Waals energy if  $\vecd$ is identified with the gradient of the phase field variable.
Then 
we recall  the following well-known statement on energy dissipation. Note that  all expressions in \eqref{EnergydissipationNSCH} are bounded quantities for a classical solution. The result is formulated for the boundary conditions in  \eqref{eq:boundNSCH} but holds also for $\Omega={\mathbb{T}}^d$.

\begin{thm}[Energy dissipation for \eqref{eq:NSCH}]
\label{NSCH-theorem}
Let  $  (p,\vecu,c)^T: \bar \Omega_T \to \setR^{d+2}$ be a classical  solution of the initial boundary value problem for (\ref{eq:NSCH}) satisfying 
\eqref{eq:iniNSCH}, \eqref{eq:boundNSCH}.\\
Then, we have for all $t \in [0,T]$ the energy balance relation  
\end{thm}
\begin{equation}\label{EnergydissipationNSCH}
\int_{\Omega} E [\vecu,c, \nabla c](t, \vecx) \; d\vecx + \int_0^t\int_{\Omega} \left( |\nabla \mu(c(s,\vecx))|^2 + \nu |\nabla \vecu(s,\vecx)|^{2} \right) d\vecx ds 
=  \int_{\Omega} E [\vecu_0,c_0, \nabla c_0](\vecx) \; d\vecx .
\end{equation}


\subsection{A low-order relaxation approximation of the NSCH system} \label{subsec:approx}

In \cite{KKR25}, an approximation of the NSCH system \eqref{eq:NSCH} has been suggested that combines the artificial compressibility idea of \cite{Temam1969I} and  the combined friction-type/screened Poisson ansatz as in \cite{Dhaouadi24,rohde2010local}.
These three  approximations come with the  tuple $\varepsilon=(\alpha,\beta,\delta,\tau)\in(0,1)^4$
of small approximation parameters. Here we rely on a modification of the model from \cite{KKR25} that
we refer to in the sequel as relaxation approximation.\\ 
The relaxation approximation   of the NSCH system \eqref{eq:NSCH} writes as 
\begin{equation}\label{eq:NSCHr}
	\begin{array}{rcl}
		p_t^\varepsilon 
		+  \ds \frac{1}{\alpha} \nabla \cdot \mathbf{u}^\varepsilon 
		&=& ~ 0, \\[1.2ex]
	\ds 	\mathbf{u}_t^\varepsilon 
		+  (\mathbf{u}^\varepsilon \cdot \nabla)\mathbf{u}^\varepsilon 
		+ \frac{1}{2} (\nabla \cdot \mathbf{u}^\varepsilon)\mathbf{u}^\varepsilon 
		+ \nabla p^\varepsilon  
		&=& \ds ~ -c^\varepsilon \nabla \left( W'(c^\varepsilon)+ \frac{1}{\beta} (c^\varepsilon-\omega^\varepsilon) \right)
		\\[2.3ex]
      &&		\ds \mathcolor{white}{=}~ + \nu \Delta \mathbf{u}^\varepsilon+(\lambda+\nu)\nabla\left(\nabla\cdot\mathbf{u}^\varepsilon\right), \\[1.2ex]
	\ds 	c_t^\varepsilon 
		+  \nabla \cdot (c^\varepsilon \mathbf{u}^\varepsilon) 
		+ \nabla \cdot \mathbf{v}^\varepsilon 
		&=& ~ \tau\Delta c^\varepsilon, \\[1.2ex]
	\ds 	\mathbf{v}_t^\varepsilon 
		+ ~ \frac{1}{\delta} \nabla \left(W'(c^\varepsilon) + \frac{1}{\beta} c^\varepsilon \right) 
		&=& \ds ~ \ds -\frac{\mathbf{v}^\varepsilon}{\delta} 
		+ \frac{1}{\delta \beta} \nabla \omega^\varepsilon+\nabla\left(\nabla\cdot\mathbf{v}^\varepsilon\right), \\[2.3ex]
		\ds  -\gamma \Delta \omega^\varepsilon 
		+ \frac{1}{\beta} \omega^\varepsilon 
		&=& \ds  ~ \frac{1}{\beta} c^\varepsilon
	\end{array} \qquad \text{ in $\Omega_T$.}
\end{equation}
The unknowns are  the   pressure   $ p^{\eps}=p^{\eps}(t, \vecx)$,  the velocity 
$\vecu^{\eps} = \vecu^{\eps}(t, \vecx)$, the phase-field variable  
$c^{\eps} = c^{\eps}(t, \vecx)$, the Cahn-Hilliard flux  $\vecv^{\eps} = \vecv^{\eps}(t, \vecx)$, and 
the relaxation variable $\omega^{\eps}=\omega^{\eps}(t, \vecx)$. 
\\
The modifications in the relaxation approximation \eqref{eq:NSCHr} refer to three items: 1) we consider for the  viscous part of the stress tensor the kinematic shear viscosity and the (artificial) bulk viscosity
$\lambda>0 $. 
The latter has not been accounted for in \cite{KKR25}; 2) we add the term  $\tau\Delta c^\varepsilon$ in the 
conservation law for the phase field variable $c^\varepsilon$; 3) we add the regularization term $\nabla\left(\nabla\cdot\mathbf{v}^\varepsilon\right)$ in the equations for the Cahn-Hilliard flux.\\ 
%
These three extensions are needed to provide  
additional dissipation terms in 
 the proof of Theorem
\ref{thm:high order convergence result} below. We will see that the terms $(\lambda+\nu)\nabla\left(\nabla\cdot\mathbf{u}^\varepsilon\right)$, $\tau \Delta c^\varepsilon$ and 
$\nabla\left(\nabla\cdot\mathbf{v}^\varepsilon\right)$
vanish in the limit when $|\varepsilon|\to 0$.

The system \eqref{eq:NSCHr} is endowed with the initial conditions 
\begin{equation}\label{eq:iniNSCHr}
	p^\varepsilon(0,\cdot)=p^\eps_0, \quad  \mathbf{u}^\varepsilon(0,\cdot) =\vecu^\eps_0,\quad  \mathbf{v}^\varepsilon(0,\cdot) =\vecv^\eps_0,    \quad c^\varepsilon(0,\cdot)=c^\eps_0
    \text{ in } \Omega,
\end{equation}
and, if applicable,   the boundary conditions
\begin{equation}\label{eq:NSCHrBC}
	\mathbf{u}^\varepsilon=\mathbf{0},\quad\nabla\omega^\varepsilon\cdot\mathbf{n}=0, \quad\nabla c^\varepsilon\cdot\mathbf{n}=0,\quad\mathbf{v}^\varepsilon\cdot\mathbf{n}=0   \text{ on }   (0,T) \times \partial\Omega.
\end{equation}
Note that the Neumann condition for the phase field and the  Neumann condition for the chemical potential $\mu$ in \eqref{eq:boundNSCH} are transferred
to a  Neumann condition for $\omega^\varepsilon$ and a no-flux condition for the Cahn-Hilliard flux, respectively.  Analogously  as for the NSCH system,  a classical solution of  the  relaxation approximation  \eqref{eq:NSCHr} that satisfies   \eqref{eq:iniNSCHr} 
is a tuple of functions $\mathbf{U}^\eps=(p^\eps, \mathbf{u}^\eps, c, \mathbf{v}^\eps, \omega^\eps)$ such that
	\begin{align*}
		p^\eps\in&~C^1\left((0,T)\times\Omega\right)\cap C\left([0,T]\times\overline{\Omega}\right),\\
		\mathbf{u}^\eps\in&~C^1\left((0,T);C(\Omega;\R^d)\right)\cap C\left([0,T];C^2\left(\overline{\Omega};\R^d\right)\right),\\
		c^\eps\in&~C^1\left((0,T);C(\Omega)\right)\cap C\left([0,T];C^1\left(\overline{\Omega}\right)\right),\\
		\mathbf{v}^\eps\in&~C^1\left((0,T);C(\Omega;\R^d)\right)\cap C\left([0,T];C^2\left(\overline{\Omega}\right)\right),\\
		\omega^\eps\in&~C^1\left((0,T)\times\Omega\right)\cap C\left([0,T];C^2\left(\overline{\Omega}\right)\right)
	\end{align*}
    holds,  
	and such that the equations in \eqref{eq:NSCHr} and  \eqref{eq:iniNSCHr} 
    are satisfied. We are not aware of any global well-posedness results for classical solutions but note that the initial value problem on  $\setR^d$ can be analyzed following the results on hyperbolic-parabolic systems in \cite{Serrebook}.
For our analysis, we will work with a weaker solution concept that is provided in Definition \ref{def:weak solutions NSCHr} and relies on the energy structure of the relaxation approximation. 
Indeed, for the
energy term 
\begin{equation}\label{eq:energy}
   \ds   E^{\epsilon}[p,\vecu,c,\vecv, \omega, \vece] =
   \ds    \frac{\alpha}{2}{p^{2} } + \frac12 {|\vecu|^2} + \dfrac{\delta}{2}|\vecv|^{2} + W(c) + 
    \frac1{2\beta} (c-\omega)^2+
     \frac{\gamma}{2} {|\vece|^{2}} \qquad 
\left((p,\vecu,c,\vecv,\omega,\vece) \in \setR^{3d+3}\right),
\end{equation}
we obtain as in \cite{KKR25} the statement
%
%
%
%
%
\begin{thm}[Energy dissipation for \eqref{eq:NSCHr}]
\label{EDrelaxedNSE}
 For $\eps \in (0,1)^4$, let  $
 (p^\eps,\vecu^\eps,c^\eps,\vecv^\eps, \omega^\eps )^T  : \bar \Omega_{T} \to \setR^{2d+3}$ be a classical solution of  (\ref{eq:NSCHr}), \eqref{eq:iniNSCHr} and \eqref{eq:NSCHrBC}:
Then we have for all $t \in [0,T]$ the energy balance
equation
\[
\begin{array}{rcl}
 \lefteqn{\ds \int_{\Omega} E^{\varepsilon}\left[p^{\varepsilon},   \mathbf{u}^\eps, c^\varepsilon, \mathbf{v}^\eps, \omega^\eps, \nabla\omega^\eps\right]\; d\vecx}\\[1.5ex] 
  &=  & \ds \int_{\Omega}  E^{\varepsilon}\left[p^{\varepsilon}_0,   \mathbf{u}^\eps_0, c^\varepsilon_0, \mathbf{v}^\eps_0, \omega^\eps(0,\cdot),  \nabla\omega^\eps(0,\cdot)\right] \; d\vecx \\[1.7ex]
  && {} \ds -   \int_0^t \int_{\Omega} \nu |\nabla \vecu^\eps|^{2}+(\lambda+\nu)(\nabla\cdot\mathbf{u}^\varepsilon)^2 + \tau|\nabla c^\varepsilon|^2 +  |\vecv^\eps|^{2} + \delta (\nabla\cdot\vecv^\varepsilon)^2 d\vecx ds.
\end{array}
\]
\end{thm} 

One readily sees that the  energy dissipation rate statement formally  reduces to \eqref{EnergydissipationNSCH} if the 
approximation parameters vanish. Note that the energy \eqref{eq:energy} is non-convex due to the presence of the double-well function $W$.


\section{Convergence of the relaxation approximation }\label{sec:analysis}

Throughout this section, we will consider the case of periodic boundary conditions, i.e., $\Omega=\mathbb{T}^d$ the torus. Note however, that a generalization of the presented results to the case of physical boundary conditions is possible; cf.~Remark \ref{rem:convergence in physical domains}. Our main goal is to show that solutions of \eqref{eq:NSCHr} converge to \eqref{eq:NSCH} in the $|\varepsilon |\rightarrow 0$ limit, as long as \eqref{eq:NSCH} admits a classical solution with additional regularity; cf.~Theorem \ref{thm:high order convergence result}. A crucial step in this endeavor is to study the stability of solutions to \eqref{eq:NSCHr}, i.e., their closeness to solutions of
\begin{equation}\label{eq:residualNSCHr}
	\left\{
	\begin{aligned}
		p_t^\varepsilon 
		&+  \frac{1}{\alpha} \nabla \cdot \mathbf{u}^\varepsilon 
		&&= ~ -\mathcal{R}_1, \\
		\mathbf{u}_t^\varepsilon 
		&+ (\mathbf{u}^\varepsilon \cdot \nabla)\mathbf{u}^\varepsilon 
		+ \frac{1}{2} (\nabla \cdot \mathbf{u}^\varepsilon)\mathbf{u}^\varepsilon 
		+ \nabla p^\varepsilon 
		&&= ~ -c^\varepsilon \nabla \left( W'(c^\varepsilon) + \frac{1}{\beta} (c^\varepsilon-\omega^\varepsilon) \right)
		\\~&~&&\mathcolor{white}{=} ~ + \nu \Delta \mathbf{u}^\varepsilon+(\lambda+\nu)\nabla\left(\nabla\cdot\mathbf{u}^\varepsilon\right)-\mathcal{R}_2, \\
		c_t^\varepsilon 
		&+  \nabla \cdot (c^\varepsilon \mathbf{u}^\varepsilon) 
		+ \nabla \cdot \mathbf{v}^\varepsilon 
		&&= ~\tau\Delta c^\varepsilon -\mathcal{R}_3, \\
		\mathbf{v}_t^\varepsilon 
		&+  \frac{1}{\delta} \nabla \left(W'(c^\varepsilon) + \frac{1}{\beta} c^\varepsilon \right) 
		&&= ~ -\frac{\mathbf{v}^\varepsilon}{\delta} 
		+ \frac{1}{\delta \beta} \nabla \omega^\varepsilon+\nabla\left(\nabla\cdot\mathbf{v}^\varepsilon\right)-\mathcal{R}_4, \\
		&  -\gamma \Delta \omega^\varepsilon 
		+ \frac{1}{\beta} \omega^\varepsilon 
		&&= ~ \frac{1}{\beta} c^\varepsilon-\mathcal{R}_5
	\end{aligned}
	\right.,
\end{equation}
w.r.t.~given residuals $\mathcal{R}_i\in C\left([0,T]\times\Omega\right),~i\in\{1,2,3,4,5\}$, and initial data. 

In the following, solutions to the relaxation system \eqref{eq:NSCHr} will be denoted as $\mathbf{U}=[p, \mathbf{u}, c, \mathbf{v}, \omega]$ while solutions to the relaxation system with residuals \eqref{eq:residualNSCHr} will be denoted as $\bar{\mathbf{U}}=[\bar{p}, \bar{\mathbf{u}}, \bar{c}, \bar{\mathbf{v}}, \bar{\omega}]$, i.e., in both cases we drop the $^\varepsilon$ to increase readability.

We will study the stability of \eqref{eq:NSCHr} via a relative energy analysis for the energy given in \eqref{eq:energy}.
In this analysis, we will compare higher-regularity solutions to the NSCH system \eqref{eq:NSCH} -- interpreted appropriately as classical solutions to the relaxation system with residuals \eqref{eq:residualNSCHr}, details will be provided in Proposition \ref{prop:high-order} -- with weak solutions to the relaxation system \eqref{eq:NSCHr}. We thus formulate the following solution concepts.

\begin{Definition}[Classical solutions to the relaxation system with residuals]\label{eq:strong solutions NSCHr}
A classical solution to \eqref{eq:residualNSCHr} is a tuple of functions $\bar{\mathbf{U}}=[\bar{p}, \bar{\mathbf{u}}, \bar{c}, \bar{\mathbf{v}}, \bar{\omega}]$ such that
	\begin{align*}
		\bar{p}\in&~C^1\left((0,T)\times\Omega\right)\cap C\left([0,T]\times\Omega\right)\\
		\bar{\mathbf{u}}\in&~C^1\left((0,T);C(\Omega;\R^d)\right)\cap C\left([0,T];C^2\left(\Omega;\R^d\right)\right)\\
		\bar{c}\in&~C^1\left((0,T);C(\Omega)\right)\cap C\left([0,T];C^1\left(\Omega\right)\right)\\
		\bar{\mathbf{v}}\in&~C^1\left((0,T);C(\Omega;\R^d)\right)\cap C\left([0,T];C^2\left(\Omega\right)\right)\\
		\bar{\omega}\in&~C^1\left((0,T)\times\Omega\right)\cap C\left([0,T];C^2\left(\Omega\right)\right),
	\end{align*}
	and the equations in \eqref{eq:residualNSCHr} and initial conditions hold in a pointwise sense.
\end{Definition}

\begin{Definition}[Weak solutions to the relaxation system]\label{def:weak solutions NSCHr}
A weak solution to \eqref{eq:NSCHr} is a tuple of functions $\mathbf{U}=[p, \mathbf{u}, c, \mathbf{v}, \omega]$ with
	\begin{align*}
		p\in&~H^1\left(0,T;L^2\left(\Omega\right)\right)\cap L^\infty\left(0,T; L^2\left(\Omega\right)\right)\\
		\mathbf{u}\in&~L^2\left(0,T ; H^1(\Omega;\R^d)\right)\cap L^\infty\left(0,T;L^2\left(\Omega;\R^d\right)\right)\\
		c\in&~L^2\left(0,T; H^1(\Omega)\right)\cap L^\infty\left(0,T; L^4\left(\Omega\right)\right)\\
		\mathbf{v}\in&~L^2\left(0,T; H_{\operatorname{div}}\left(\Omega;\R^d\right)\right)\cap L^\infty\left(0,T; L^2\left(\Omega\right)\right)\\
		\omega\in&~H^1\left((0,T)\times\Omega\right),
	\end{align*}
	where $H_{\operatorname{div}}(\Omega;\R^d)=\left\{\mathbf{v}\in L^2(\Omega;\R^d)~|~\nabla\cdot\mathbf{v}\in L^2(\Omega;\R^d)\right\}$ and such that for any $0\leq t_1< t_2\leq T$ there holds
    \begin{equation}\label{eq:weak NSCHr}
\begin{aligned}
p_t + \frac{1}{\alpha}\nabla\cdot\mathbf{u} &= 0, \\[1.2ex]
\int_{t_1}^{t_2}\!\!\int_\Omega \mathbf{u}\cdot\boldsymbol{\phi}_t
- (\mathbf{u}\cdot\nabla)\mathbf{u}\cdot\boldsymbol{\phi}
- \tfrac12(\nabla\cdot\mathbf{u})\mathbf{u}\cdot\boldsymbol{\phi}
+ p\nabla\cdot\boldsymbol{\phi}~d\vecx dt
&= \int_{t_1}^{t_2}\!\!\int_\Omega c\nabla\!\left(W'(c)+\tfrac1\beta(c-\omega)\right)\!\cdot\boldsymbol{\phi} \\
&\quad + \nu\nabla\mathbf{u}:\nabla\boldsymbol{\phi}
+ (\lambda+\nu)\nabla\cdot\mathbf{u}\nabla\cdot\boldsymbol{\phi}\,d\mathbf{x}dt \\
&\quad + \left[\int_\Omega \mathbf{u}\cdot\boldsymbol{\phi}\,d\mathbf{x}\right]_{t_1}^{t_2},
\\[1.2ex]
\int_{t_1}^{t_2}\!\!\int_\Omega c\phi_t
- \phi\nabla\cdot(c\mathbf{u})
- \phi\nabla\cdot\mathbf{v}~d\vecx dt
&= \int_{t_1}^{t_2}\!\!\int_\Omega \tau\nabla c\cdot\nabla\phi\,d\mathbf{x}dt
+ \left[\int_\Omega c\phi\,d\mathbf{x}\right]_{t_1}^{t_2},
\\[1.2ex]
\int_{t_1}^{t_2}\!\!\int_\Omega \mathbf{v}\cdot\boldsymbol{\psi}_t
+ \frac1\delta\!\left(W'(c)+\tfrac1\beta(c-\omega)\right)\!\nabla\cdot\boldsymbol{\psi}~d\vecx dt
&= \int_{t_1}^{t_2}\!\!\int_\Omega
\frac{1}{\delta}\mathbf{v}\cdot\boldsymbol{\psi}
+ \nabla\cdot\mathbf{v}\nabla\cdot\boldsymbol{\psi}\,d\mathbf{x}dt \\
&\quad + \left[\int_\Omega \mathbf{v}\cdot\boldsymbol{\psi}\,d\mathbf{x}\right]_{t_1}^{t_2},
\\[1.2ex]
\int_\Omega \gamma\nabla\omega\cdot\nabla\xi\,d\mathbf{x}
&= \int_\Omega \tfrac1\beta(c-\omega)\xi\,d\mathbf{x}
\end{aligned}
\end{equation}
for arbitrary test functions $\boldsymbol{\phi}, \boldsymbol{\psi}\in C^\infty\left([0,T]\times\Omega;\R^d\right), \phi\in C^\infty\left([0,T]\times\Omega\right), \xi\in C^\infty(\Omega)$. 

The initial conditions are imposed strongly for the pressure and 
weakly for the other terms, i.e.,
\begin{align}
    p(0,\cdot)=p_0,\qquad\lim\limits_{t\to 0+}\int_\Omega (\mathbf{u},c,\mathbf{v})\cdot\boldsymbol{\Phi}~d\vecx=\int_\Omega (\mathbf{u}_0,c_0,\mathbf{v}_0)\cdot\boldsymbol{\Phi}~d\vecx,
\end{align}
for all $\boldsymbol{\Phi}\in C^\infty\left(\Omega;\R^d\right)\times C^\infty\left(\Omega\right)\times C^\infty\left(\Omega;\R^d\right)$.

Furthermore, for all $\psi\in C^\infty_c\left([0,T)\times\Omega\right)$, $\psi\geq 0$, we demand that the energy inequality
\begin{align}\label{eq:weak energy inequality}
    \int_0^T\int_\Omega E\psi_t+\mathbf{Q}\cdot\nabla\psi~d\vecx dt\geq \int_0^T\int_\Omega S\psi~d\vecx dt -\int_\Omega [E\psi](0,\cdot)~d\vecx,
\end{align}
holds with
\begin{align}
    E:=&~E^\varepsilon(\mathbf{U}),\\
     \mathbf{Q}=\mathbf{Q}[\mathbf{U}]:=&~p\mathbf{u}+\frac{1}{2}|\mathbf{u}|^2\mathbf{u}+\left(W(c)+\frac{1}{\beta}(c-\omega)\right)(c\mathbf{u}+\mathbf{v})\nonumber\\
    &-\left(\nu\nabla\mathbf{u}+(\lambda+\nu)(\nabla\cdot\mathbf{u})\right)\mathbf{u}-\tau c\nabla c-(\nabla\cdot\mathbf{v})\mathbf{v}-\gamma\omega_t\nabla\omega,\\
    S=S[\mathbf{U}]:=&~\lambda|\nabla\mathbf{u}|^2+(\lambda+\nu)(\nabla\cdot\mathbf{u})^2+\tau|\nabla c|^2+\frac{1}{\delta}|\mathbf{v}|^2+(\nabla\cdot\mathbf{v})^2.
\end{align}
\end{Definition}

\subsection{The relative energy analysis}
Note that in the following, we write $a\lesssim b$ if there exists a constant $k\geq0$ that is independent of $\varepsilon, \nu, \lambda, \gamma$ such that $a\leq kb$.

The (full) relative energy between two states $\mathbf{U}=(p,\vecu,c,\vecv, \omega, \vece), \bar{\mathbf{U}}=(\bar{p},\bar{\vecu},\bar{c},\bar{\vecv}, \bar{\omega}, \bar{\vece})\in \setR^{3d+3}$ is given by
\begin{align}\label{eq:relenergy}
	\widetilde{\eta}^\varepsilon(\mathbf{U},\bar{\mathbf{U}}):=&~\frac{\alpha}{2}(p-\bar{p})^2+W(c)-W(\bar{c})-W'(\bar{c})(c-\bar{c})+\frac{1}{2\beta}(c-\omega-(\bar{c}-\bar{\omega}))^2\nonumber\\
	&+\frac{1}{2}|\mathbf{u}-\bar{\mathbf{u}}|^2+\frac{\delta}{2}|\mathbf{v}-\bar{\mathbf{v}}|^2+\frac{\gamma}{2}|\mathbf{e}-\bar{\mathbf{e}}|^2.
\end{align}
Due to the double-well structure of $W$, the relative energy $\widetilde{\eta}^\varepsilon$ is not necessarily nonnegative. This is unfavorable for a subsequent Gronwall argument. However, by disregarding the non-convex part of $W$ and adding $\frac{1}{2}(c-\bar{c})^2$, we can define the \textit{reduced relative energy}
\begin{align}\label{eq:relenergy}
	\eta^\varepsilon(\mathbf{U},\bar{\mathbf{U}}):=&~\frac{\alpha}{2}(p-\bar{p})^2+\frac{1}{4}c^4-\frac{1}{4}\bar{c}^4-\bar{c}^3(c-\bar{c})+\frac{1}{2\beta}(c-\omega-(\bar{c}-\bar{\omega}))^2\nonumber\\
	&+\frac{1}{2}|\mathbf{u}-\bar{\mathbf{u}}|^2+\frac{\delta}{2}|\mathbf{v}-\bar{\mathbf{v}}|^2+\frac{\gamma}{2}|\mathbf{e}-\bar{\mathbf{e}}|^2+\frac{1}{2}(c-\bar{c})^2,
\end{align}
which is nonnegative. Also, note that the reduced relative energy allows us to bound
\begin{align*}
    (c-\bar{c})^4\lesssim \eta^\varepsilon,
\end{align*}
due to the following proposition.
\begin{Proposition}\label{prop:W}
	There holds
	\begin{align}
			0\leq(2-\sqrt{3})(c-\bar{c})^4\leq c^4-\bar{c}^4-4\bar{c}^3(c-\bar{c}),
	\end{align}
	for all $c,\bar{c}\in\R$.
\end{Proposition}
\begin{proof}
	First note that 
	\begin{align}
		c^4-\bar{c}^4=(c^3+c^2\bar{c}+\bar{c}^2c+\bar{c}^3)(c-\bar{c}),
	\end{align}
	thus
	\begin{align}
		c^4-\bar{c}^4-4\bar{c}^3(c-\bar{c})=&~(c^3+c^2\bar{c}+\bar{c}^2c-3\bar{c}^3)(c-\bar{c})\nonumber\\
		=&~(c^2+2c\bar{c}+3\bar{c}^2)(c-\bar{c})^2,
	\end{align}
	while
	\begin{align}
		c^2+2c\bar{c}+3\bar{c}^2=\begin{pmatrix}
			c&\bar{c}
		\end{pmatrix}
		\begin{pmatrix}
			1&1\\1&3
		\end{pmatrix}
		\begin{pmatrix}
			c\\\bar{c}
		\end{pmatrix}.
	\end{align}
	Finally, the matrix $\begin{pmatrix}
		1&1\\1&3
	\end{pmatrix}$ is diagonizable with eigenvalues $\lambda_\pm=2\pm\sqrt{3}$ and thus the claim follows.
\end{proof}

Furthermore, since
\begin{align}\label{eq:reduced rel en}
	\eta^\varepsilon=\widetilde{\eta}^\varepsilon+(c-\bar{c})^2,
\end{align}
a Gronwall argument for the reduced relative energy $\eta^\varepsilon$ becomes available, once we appropriately bound the time evolution of the full relative energy $\widetilde{\eta}^\varepsilon$ and of $(c-\bar{c})^2$. We start with the former.

\begin{Lemma}[Growth rate of the full relative energy]
	 Let $\mathbf{U}$ be a weak solution to \eqref{eq:NSCHr} and let $\bar{\mathbf{U}}$ be a classical solution to \eqref{eq:residualNSCHr}. Then, for every $t\leq T$ there holds the following.
	\begin{align}\label{eq:growth rate full relative energy}
		\int_{\Omega}\widetilde{\eta}^\varepsilon(\mathbf{U},\bar{\mathbf{U}})(t,\mathbf{x})~d\mathbf{x}\leq ~&\int_{\Omega}\widetilde{\eta}^\varepsilon(\mathbf{U},\bar{\mathbf{U}})(0,\mathbf{x})~d\mathbf{x}
       + \int_0^t\int_\Omega-6\tau(c-\bar{c})\nabla(c-\bar{c})\cdot\nabla\bar{c}-3\tau(c-\bar{c})^2\nabla(c-\bar{c})\cdot\nabla\bar{c}\nonumber\\
		&+\tau|\nabla(c-\bar{c})|^2-(W'(c)-W'(\bar{c})-W''(\bar{c})(c-\bar{c}))\nabla\cdot\left(\bar{c}\bar{\mathbf{u}}+\bar{\mathbf{v}}\right)-\delta(\nabla\cdot(\mathbf{v}-\bar{\mathbf{v}}))^2\nonumber\\
		&-\left(\frac{1}{2}(c-\bar{c})^2W''(\bar{c})+\frac{1}{3}(c-\bar{c})^3W'''(\bar{c})+\frac{3}{4}(c-\bar{c})^4\right)\nabla\cdot\bar{\mathbf{u}}-|\mathbf{v}-\bar{\mathbf{v}}|^2\nonumber\\
		&+\left(\frac{1}{2}(c-\bar{c})^2W'''(\bar{c})+(c-\bar{c})^3\right)\bar{\mathbf{u}}\cdot\nabla \bar{c}-(c-\bar{c})(\mathbf{u}-\bar{\mathbf{u}})\cdot\nabla W'(\bar{c})\nonumber\\
		&-(c-\bar{c})\left(\mathbf{u}-\bar{\mathbf{u}}\right)\cdot\nabla(\bar{c}-\bar{\omega})-\frac{1}{2\beta}(c-\omega-(\bar{c}-\bar{\omega}))^2\nabla\cdot\bar{\mathbf{u}}-\frac{\tau}{2\beta}|\nabla(c-\bar{c})|^2\nonumber\\
        &+\gamma\frac{|\nabla(\omega-\bar{\omega})|^2}{2}\nabla\cdot\bar{\mathbf{u}}+\gamma\left(\nabla(\omega-\bar{\omega})\otimes\nabla(\omega-\bar{\omega})\right):\nabla\bar{\mathbf{u}}+\frac{\tau}{2\beta}|\nabla(\omega-\bar{\omega})|^2\nonumber\\
        &+\gamma(\omega-\bar{\omega})\nabla(\omega-\bar{\omega})\cdot\nabla\left(\nabla\cdot\bar{\mathbf{u}}\right)-\nabla\bar{\mathbf{u}}:((\mathbf{u}-\bar{\mathbf{u}})\otimes(\mathbf{u}-\bar{\mathbf{u}}))\nonumber\\
		&-\frac{1}{2}(\nabla\cdot(\mathbf{u}-\bar{\mathbf{u}}))\bar{\mathbf{u}}\cdot(\mathbf{u}-\bar{\mathbf{u}})-\nu|\nabla(\mathbf{u}-\bar{\mathbf{u}})|^2-(\lambda+\nu)(\nabla\cdot(\mathbf{u}-\bar{\mathbf{u}}))^2\nonumber\\
		&+\alpha\mathcal{R}_1(p-\bar{p})+\mathcal{R}_2(\mathbf{u}-\bar{\mathbf{u}})+\frac{1}{\beta}\mathcal{R}_3(c-\omega-(\bar{c}-\bar{\omega}))+\delta\mathcal{R}_4(\mathbf{v}-\bar{\mathbf{v}})&\nonumber\\
		&+\mathcal{R}_5(\omega_t-\bar{\omega}_t)+\mathcal{R}_3 W''(\bar{c})(c-\bar{c})-(\omega-\bar{\omega})\bar{\mathbf{u}}\cdot\nabla\mathcal{R}_5~d\mathbf{x}ds.
	\end{align}
\end{Lemma}
\begin{proof}
        For ease of presentation, we will undertake all computations in this proof for  a classical solution $\mathbf{U}$. The generalization to weak solutions $\mathbf{U}$ is straightforward; cf.~e.g.~the weak-strong uniqueness proof in \cite{DafermosBook}. 
        
    	We first compute the time evolution of $\widetilde{\eta}^\varepsilon(\mathbf{U},\bar{\mathbf{U}})~\widehat{=}~\widetilde{\eta}^\varepsilon$ rather naively and then regroup appropriately. According to the evolution equations and after noticing some simple cancellations, the growth rate of the relative energy takes the preliminary form
	\begin{align}\label{eq:naiv growth rate}
		\partial_t\widetilde{\eta}^\varepsilon=&-(p-\bar{p})\nabla\cdot(\mathbf{u}-\bar{\mathbf{u}})+\alpha\mathcal{R}_1(p-\bar{p})\nonumber\\
        &-(W'(c)-W'(\bar{c}))\nabla\cdot(c\mathbf{u})+W''(\bar{c})(c-\bar{c})\nabla\cdot(\bar{c}\bar{\mathbf{u}})\nonumber\\
		&-(W'(c)-W'(\bar{c}))\nabla\cdot\mathbf{v}+W''(\bar{c})(c-\bar{c})\nabla\cdot\bar{\mathbf{v}}+\mathcal{R}_3 W''(\bar{c})(c-\bar{c})\nonumber\\
        &+\tau(W'(c)-W'(\bar{c}))\Delta c-\tau W''(\bar{c})(c-\bar{c})\Delta \bar{c} + \frac{\tau}{\beta}(c-\omega-(\bar{c}-\bar{\omega}))\Delta(c-\bar{c})\nonumber\\
		&-\frac{1}{\beta}(c-\omega-(\bar{c}-\bar{\omega}))\nabla\cdot(c\mathbf{u}-\bar{c}\bar{\mathbf{u}})-\frac{1}{\beta}(c-\omega-(\bar{c}-\bar{\omega}))\nabla\cdot(\mathbf{v}-\bar{\mathbf{v}})+\frac{1}{\beta}\mathcal{R}_3(c-\omega-(\bar{c}-\bar{\omega}))\nonumber\\
		&-\left((\mathbf{u}\cdot\nabla)\mathbf{u}-(\bar{\mathbf{u}}\cdot\nabla)\bar{\mathbf{u}}\right)\cdot(\mathbf{u}-\bar{\mathbf{u}})-\frac{1}{2}\left((\nabla\cdot\mathbf{u})\mathbf{u}-(\nabla\cdot\bar{\mathbf{u}})\bar{\mathbf{u}}\right)\cdot(\mathbf{u}-\bar{\mathbf{u}})\nonumber\\
		&-(\mathbf{u}-\bar{\mathbf{u}})\cdot\nabla(p-\bar{p})-\left(c\nabla W'(c)-\bar{c}\nabla W'(\bar{c})\right)\cdot(\mathbf{u}-\bar{\mathbf{u}})\nonumber\\
		&-\frac{1}{\beta}\left(c\nabla(c-\omega)-\bar{c}\nabla(\bar{c}-\bar{\omega})\right)\cdot(\mathbf{u}-\bar{\mathbf{u}})+\nu\Delta(\mathbf{u}-\bar{\mathbf{u}})\cdot(\mathbf{u}-\bar{\mathbf{u}})+\mathcal{R}_2(\mathbf{u}-\bar{\mathbf{u}})\nonumber\\
		&-\nabla( W'(c)- W'(\bar{c}))\cdot(\mathbf{v}-\bar{\mathbf{v}})-|\mathbf{v}-\bar{\mathbf{v}}|^2-\frac{1}{\beta}(\mathbf{v}-\bar{\mathbf{v}})\cdot\nabla(c-\omega-(\bar{c}-\bar{\omega}))\nonumber\\
		&+\delta\Delta(\mathbf{v}-\bar{\mathbf{v}})\cdot(\mathbf{v}-\bar{\mathbf{v}})+\delta\mathcal{R}_4(\mathbf{v}-\bar{\mathbf{v}})+\mathcal{R}_5(\omega_t-\bar{\omega}_t).
	\end{align}
	Let us group some related terms and study them separately, i.e.,
	\begin{align}
		I_1:=& -(W'(c)-W'(\bar{c}))\nabla\cdot\mathbf{v}+W''(\bar{c})(c-\bar{c})\nabla\cdot\bar{\mathbf{v}}-\nabla( W'(c)- W'(\bar{c}))\cdot(\mathbf{v}-\bar{\mathbf{v}}),\\
		I_2:=&-(W'(c)-W'(\bar{c}))\nabla\cdot(c\mathbf{u})+W''(\bar{c})(c-\bar{c})\nabla\cdot(\bar{c}\bar{\mathbf{u}})-(c\nabla W'(c)-\bar{c}\nabla W'(\bar{c}))\cdot(\mathbf{u}-\bar{\mathbf{u}}),\\
        I_3:=&~\tau(W'(c)-W'(\bar{c}))\Delta c-\tau W''(\bar{c})(c-\bar{c})\Delta \bar{c},\\
		I_4:=&-\frac{1}{\beta}(c\nabla(c-\omega)-\bar{c}\nabla(\bar{c}-\bar{\omega}))\cdot(\mathbf{u}-\bar{\mathbf{u}})-\frac{1}{\beta}(c-\omega-(\bar{c}-\bar{\omega}))\nabla\cdot(c\mathbf{u}-\bar{c}\bar{\mathbf{u}}).
	\end{align}
	
	Now,
	\begin{align}\label{eq:blue}
		I_1=-\nabla\cdot[(W'(c)-W'(\bar{c}))(\mathbf{v}-\bar{\mathbf{v}})]-(W'(c)-W'(\bar{c})-W''(\bar{c})(c-\bar{c}))\nabla\cdot\bar{\mathbf{v}},
	\end{align}
	and, similarly,
	\begin{align}
		I_2=&-\nabla\cdot[(W'(c)-W'(\bar{c}))(c\mathbf{u}-\bar{c}\bar{\mathbf{u}})]-(W'(c)-W'(\bar{c})-W''(\bar{c})(c-\bar{c}))\nabla\cdot(\bar{c}\bar{\mathbf{u}})\nonumber\\
		&+(c-\bar{c})\bar{\mathbf{u}}\cdot\nabla(W'(c)-W'(\bar{c}))-(c-\bar{c})(\mathbf{u}-\bar{\mathbf{u}})\cdot\nabla W'(\bar{c}),
	\end{align}
	where, by Taylor expansion,
	\begin{align}
		(c-\bar{c})\nabla(W'(c)-W'(\bar{c}))=&~(c-\bar{c})\nabla\left(W''(\bar{c})(c-\bar{c})+\frac{W'''(\bar{c})}{2}(c-\bar{c})^2+(c-\bar{c})^3\right)\nonumber\\
		=&~\nabla\left[\frac{1}{2}(c-\bar{c})^2W''(\bar{c})+\frac{1}{3}(c-\bar{c})^3W'''(\bar{c})+(c-\bar{c})^4\right]\nonumber\\
		&+\frac{1}{2}(c-\bar{c})^2\nabla W''(\bar{c})+\frac{1}{6}(c-\bar{c})^3\nabla W'''(\bar{c}),
	\end{align}
	since $W$ is a fourth-order polynomial. Thus, we can rearrange $I_2$ as
	\begin{align}
		I_2=&-\nabla\cdot\left[(W'(c)-W'(\bar{c}))(c\mathbf{u}-\bar{c}\bar{\mathbf{u}})\right]\nonumber\\
		&+\nabla\cdot\left[\left(\frac{1}{2}(c-\bar{c})^2W''(\bar{c})+\frac{1}{3}(c-\bar{c})^3W'''(\bar{c})+\frac{1}{4}(c-\bar{c})^4\right)\bar{\mathbf{u}}\right]\nonumber\\
		&-(W'(c)-W'(\bar{c})-W''(\bar{c})(c-\bar{c}))\nabla\cdot(\bar{c}\bar{\mathbf{u}})\nonumber\\
		&-\left(\frac{1}{2}(c-\bar{c})^2W''(\bar{c})+\frac{1}{3}(c-\bar{c})^3W'''(\bar{c})+\frac{1}{4}(c-\bar{c})^4\right)\nabla\cdot\bar{\mathbf{u}}\nonumber\\
		&+\left(\frac{1}{2}(c-\bar{c})^2\nabla W''(\bar{c})+\frac{1}{6}(c-\bar{c})^3\nabla W'''(\bar{c})\right)\cdot\bar{\mathbf{u}}\nonumber\\
		&-(c-\bar{c})(\mathbf{u}-\bar{\mathbf{u}})\cdot\nabla W'(\bar{c}).\label{eq:I2}
	\end{align}

    The third term can be rewritten as follows,
    \begin{align}
        I_3=&~\tau\left(W'(c)-W'(\bar{c})-W''(\bar{c})(c-\bar{c})-\frac{W'''(\bar{c})}{2}(c-\bar{c})^2\right)\Delta(c-\bar{c})\nonumber\\
        &+\tau W''(\bar{c})(c-\bar{c})\Delta(c-\bar{c})+\frac{W'''(\bar{c})}{2}(c-\bar{c})^2\Delta(c-\bar{c}),
    \end{align}
    i.e.,
    \begin{align}
        I_3=&~\tau(c-\bar{c})^3\Delta(c-\bar{c})+\tau(3\bar{c}^2-1)(c-\bar{c})\Delta(c-\bar{c})+\tau3\bar{c}(c-\bar{c})^2\Delta(c-\bar{c})\nonumber\\
        =&~\nabla\cdot\left[\tau\left((c-\bar{c})^3+(3\bar{c}^2-1)(c-\bar{c})+3\bar{c}(c-\bar{c})^2\right)\nabla(c-\bar{c})\right]\nonumber\\
        &-3\tau(c-\bar{c})^2|\nabla(c-\bar{c})|^2-6\tau(c-\bar{c})\nabla(c-\bar{c})\cdot\nabla\bar{c}-\tau(3\bar{c}^2-1)|\nabla(c-\bar{c})|^2\nonumber\\
        &-3\tau(c-\bar{c})^2\nabla(c-\bar{c})\cdot\nabla\bar{c}-6\tau\bar{c}(c-\bar{c})|\nabla(c-\bar{c})|^2.
    \end{align}
    Note that
    \begin{align}
        &-3\tau(c-\bar{c})^2|\nabla(c-\bar{c})|^2-\tau(3\bar{c}^2-1)|\nabla(c-\bar{c})|^2-6\tau\bar{c}(c-\bar{c})|\nabla(c-\bar{c})|^2\leq \tau|\nabla(c-\bar{c})|^2,
    \end{align}
   so that
   \begin{multline}\label{eq:I3}
       I_3\leq  ~\nabla\cdot\left[\tau\left((c-\bar{c})^3+(3\bar{c}^2-1)(c-\bar{c})+3\bar{c}(c-\bar{c})^2\right)\nabla(c-\bar{c})\right]\\
        -6\tau(c-\bar{c})\nabla(c-\bar{c})\cdot\nabla\bar{c}
        -3\tau(c-\bar{c})^2\nabla(c-\bar{c})\cdot\nabla\bar{c}+ \tau|\nabla(c-\bar{c})|^2.
   \end{multline}
	
	Moving forward,
	\begin{align}
		I_4=-\nabla\cdot\left[\frac{1}{\beta}(c-\omega-(\bar{c}-\bar{\omega}))(c\mathbf{u}-\bar{c}\bar{\mathbf{u}})\right]+\frac{1}{\beta}(c-\bar{c})\left(\bar{\mathbf{u}}\cdot\nabla(c-\omega)-\mathbf{u}\cdot\nabla(\bar{c}-\bar{\omega})\right),
	\end{align}
	where
	\begin{equation}
		\frac{1}{\beta}(c-\bar{c})\left(\bar{\mathbf{u}}\cdot\nabla(c-\omega)-\mathbf{u}\cdot\nabla(\bar{c}-\bar{\omega})\right)=~\frac{1}{\beta}(c-\bar{c})\bar{\mathbf{u}}\cdot\nabla\left(c-\omega-(\bar{c}-\bar{\omega})\right)
        -(c-\bar{c})\left(\mathbf{u}-\bar{\mathbf{u}}\right)\cdot\nabla(\bar{c}-\bar{\omega}),
	\end{equation}
	and we write
	\begin{align}
		\frac{1}{\beta}(c-\bar{c})\bar{\mathbf{u}}\cdot\nabla\left(c-\omega-(\bar{c}-\bar{\omega})\right)=\frac{1}{2\beta}\bar{\mathbf{u}}\cdot\nabla(c-\omega-(\bar{c}-\bar{\omega}))^2+\frac{1}{\beta}(\omega-\bar{\omega})\bar{\mathbf{u}}\cdot\nabla\left(c-\omega-(\bar{c}-\bar{\omega})\right).
	\end{align}
	Furthermore, by the elliptic coupling, we have
	\begin{align}\label{eq:higher derivatives of omega}
		\frac{1}{\beta}(\omega-\bar{\omega})\bar{\mathbf{u}}\cdot\nabla\left(c-\omega-(\bar{c}-\bar{\omega})\right)=&~-\gamma(\omega-\bar{\omega})\bar{\mathbf{u}}\cdot\nabla\Delta(\omega-\bar{\omega})-(\omega-\bar{\omega})\bar{\mathbf{u}}\cdot\nabla\mathcal{R}_5.
	\end{align}
	Now, since for any scalar field $\phi$ and any vector field $\boldsymbol{\Phi}$,
	\begin{align}
		\phi\boldsymbol{\Phi}\cdot\nabla\Delta\phi=&~\nabla\cdot\left[\phi\Delta\phi\boldsymbol{\Phi}-\nabla\cdot\left(\phi\boldsymbol{\Phi}\right)\nabla\phi+\frac{|\nabla\phi|^2}{2}\boldsymbol{\Phi}\right]\nonumber\\
        &+\frac{|\nabla\phi|^2}{2}\nabla\cdot\boldsymbol{\Phi}+\left(\nabla\phi\otimes\nabla\phi\right):\nabla\boldsymbol{\Phi}+\phi\nabla\phi\cdot\nabla\left(\nabla\cdot\boldsymbol{\Phi}\right),
	\end{align}
	we thereby obtain
	\begin{align}
		I_4=&-\nabla\cdot\left[\frac{1}{\beta}(c-\omega-(\bar{c}-\bar{\omega}))(c\mathbf{u}-\bar{c}\bar{\mathbf{u}})-\frac{1}{2\beta}\left(c-\omega-(\bar{c}-\bar{\omega})\right)^2\bar{\mathbf{u}}\right]\nonumber\\
        &+\nabla\cdot\left[\gamma(\omega-\bar{\omega})\Delta(\omega-\bar{\omega})\bar{\mathbf{u}}-\gamma\nabla\cdot\left((\omega-\bar{\omega})\bar{\mathbf{u}}\right)\nabla(\omega-\bar{\omega})+\gamma\frac{|\nabla(\omega-\bar{\omega})|^2}{2}\bar{\mathbf{u}}\right]\nonumber\\
        &-(c-\bar{c})\left(\mathbf{u}-\bar{\mathbf{u}}\right)\cdot\nabla(\bar{c}-\bar{\omega})-\frac{1}{2\beta}(c-\omega-(\bar{c}-\bar{\omega}))^2\nabla\cdot\bar{\mathbf{u}}-(\omega-\bar{\omega})\bar{\mathbf{u}}\cdot\nabla\mathcal{R}_5\nonumber\\
        &+\gamma\frac{|\nabla(\omega-\bar{\omega})|^2}{2}\nabla\cdot\bar{\mathbf{u}}+\gamma\left(\nabla(\omega-\bar{\omega})\otimes\nabla(\omega-\bar{\omega})\right):\nabla\bar{\mathbf{u}}+\gamma(\omega-\bar{\omega})\nabla(\omega-\bar{\omega})\cdot\nabla\left(\nabla\cdot\bar{\mathbf{u}}\right)\label{eq:I4}.
	\end{align}
	
	Finally, the convective terms in the relative energy balance \eqref{eq:naiv growth rate} amount to
	\begin{align}\label{eq:convective terms}
		&-[(\mathbf{u}\cdot\nabla)\mathbf{u}-(\bar{\mathbf{u}}\cdot\nabla)\bar{\mathbf{u}}]\cdot(\mathbf{u}-\bar{\mathbf{u}})-\frac{1}{2}[(\nabla\cdot\mathbf{u})\mathbf{u}-(\nabla\cdot\bar{\mathbf{u}})\bar{\mathbf{u}}]\cdot(\mathbf{u}-\bar{\mathbf{u}})\nonumber\\
		=&-\nabla\cdot\left[\frac{|\mathbf{u}-\bar{\mathbf{u}}|^2}{2}\mathbf{u}\right]-((\mathbf{u}-\bar{\mathbf{u}})\cdot\nabla)\bar{\mathbf{u}}\cdot(\mathbf{u}-\bar{\mathbf{u}})-\frac{1}{2}(\nabla\cdot(\mathbf{u}-\bar{\mathbf{u}}))\bar{\mathbf{u}}\cdot(\mathbf{u}-\bar{\mathbf{u}}).
	\end{align}
The proof can be completed by adding \eqref{eq:blue}, \eqref{eq:I2}, \eqref{eq:I3}, \eqref{eq:I4} and \eqref{eq:convective terms}.
    
\end{proof}

\begin{rem}\label{rem:powers of c}
    Since the terms $(c-\bar{c})^4$ and $(c-\bar{c})^2$ are controlled by definition of the reduced relative energy, we can simply use Young's inequality to handle the cubic term $(c-\bar{c})^3$.
\end{rem}

\begin{rem}
    Note that the integral of the term $\gamma(\omega-\bar{\omega})\nabla(\omega-\bar{\omega})\cdot\nabla(\nabla\cdot\bar{\mathbf{u}})$ in \eqref{eq:growth rate full relative energy} can be controlled by a multiple of the integral of the reduced relative entropy without resorting to Poincare's inequality, since
    \begin{align*}
        (\omega-\bar{\omega})^2\leq \frac{1}{2\beta}(\omega-\bar{\omega})^2+\frac{\beta}{2}(c-\omega-(\bar{c}-\bar{\omega}))^2+\frac{1}{2}(\omega-\bar{\omega})^2+\frac{1}{2}(c-\bar{c})^2,
    \end{align*}
    so that 
      \begin{align*}
       \frac{1}{2} (\omega-\bar{\omega})^2\leq \frac{1}{2\beta}(\omega-\bar{\omega})^2+\frac{\beta}{2}(c-\omega-(\bar{c}-\bar{\omega}))^2+\frac{1}{2}(c-\bar{c})^2.
    \end{align*}
\end{rem}

To control the growth rate of the reduced relative entropy $\eta^\varepsilon$, we recall \eqref{eq:reduced rel en}. Our next step is to compute the rate of $(c-\bar{c})^2$. Note that this term is nonlinear in the weak solution $\mathbf{U}$ and, in particular, not controlled in the energy inequality \eqref{eq:weak energy inequality}. This means that some care is needed to carry it out without resorting to classical solutions $\mathbf{U}$.

\begin{Proposition}\label{prop:c-bc weak evolution} 
Let $\mathbf{U}$ be a weak solution to \eqref{eq:NSCHr} and $\bar{\mathbf{U}}$ be a classical solution to \eqref{eq:residualNSCHr}. Then, for every $t\leq T$ there holds
    \begin{align}\label{eq:growth rate of c-barc}
        \int_{\Omega}\frac{(c-\bar{c})^2}{2}(t,\vecx)~d\mathbf{x} = & \int_{\Omega}\frac{(c-\bar{c})^2}{2}(0,\vecx)~d\mathbf{x}+\int_0^t\int_\Omega-\tau|\nabla(c-\bar{c})|^2-\frac{1}{2}(c-\bar{c})^2\nabla\cdot(\mathbf{u}-\bar{\mathbf{u}})
	\nonumber\\
    &-(c-\bar{c})(\mathbf{u}-\bar{\mathbf{u}})\cdot\nabla\bar{c}-\bar{c}(c-\bar{c})\nabla\cdot(\mathbf{u}-\bar{\mathbf{u}})-\frac{1}{2}(c-\bar{c})^2\nabla\cdot\bar{\mathbf{u}}\nonumber\\
    &-(c-\omega-(\bar{c}-\bar{\omega}))\nabla\cdot(\mathbf{v}-\bar{\mathbf{v}})+(\mathbf{v}-\bar{\mathbf{v}})\nabla(\omega-\bar{\omega})+\mathcal{R}_3(c-\bar{c})~d\mathbf{x}ds.
    \end{align}
\end{Proposition}

\begin{proof}
First we note that, by testing the evolution equation for $c$ with $\phi\in C^\infty_c\left((0,T)\times\Omega\right)$, it is evident that $c\in H^1\left((0,T);H^{-1}(\Omega)\right)$. Thus, it follows that for any $\psi\in C^\infty_c\left((0,T)\times\Omega\right)$, $\psi\geq0$, there holds
\begin{align}
    \int_0^T\int_\Omega \frac{(c-\bar{c})^2}{2}\psi_t~d\mathbf{x}dt=& -\int_0^T\left\langle c_t-\bar{c}_t,(c-\bar{c})\psi\right\rangle_{H^{-1}(\Omega),H^1(\Omega)} dt,
\end{align}
which allows us to use the evolution equation of $c-\bar{c}$ directly, i.e.,
\begin{align}
    -\int_0^T\left\langle(c_t-\bar{c}_t),(c-\bar{c})\psi\right\rangle_{H^{-1}(\Omega),H^1(\Omega)} dt=&~\int_0^T\int_\Omega \psi(c-\bar{c})\nabla\cdot(c\mathbf{u}-\bar{c}\bar{\mathbf{u}})+\psi(c-\bar{c})\nabla\cdot(\mathbf{v}-\bar{\mathbf{v}})\nonumber\\
    &~\qquad\quad -\tau\nabla(c-\bar{c})\cdot\nabla\left(\psi(c-\bar{c})\right)-\psi\mathcal{R}_3(c-\bar{c})~d\mathbf{x}dt.
\end{align}
Now, on the one hand
\begin{align}
    \psi(c-\bar{c})\nabla\cdot(c\mathbf{u}-\bar{c}\bar{\mathbf{u}})=&~\psi(c-\bar{c})\nabla\cdot\left((c-\bar{c})(\mathbf{u}-\bar{\mathbf{u}})+(c-\bar{c})\bar{\mathbf{u}}+\bar{c}(\mathbf{u}-\bar{\mathbf{u}})\right)\nonumber\\
    =&~\nabla\cdot\left[\psi\frac{(c-\bar{c})^2}{2}\mathbf{u}\right]+\psi\frac{(c-\bar{c})^2}{2}\nabla\cdot(\mathbf{u}-\bar{\mathbf{u}})-\frac{(c-\bar{c})^2}{2}(\mathbf{u}-\bar{\mathbf{u}})\cdot\nabla\psi\nonumber\\
    &+\psi\frac{(c-\bar{c})^2}{2}\bar{\mathbf{u}}-\frac{(c-\bar{c})^2}{2}\bar{\mathbf{u}}\cdot\nabla\psi+\psi(c-\bar{c})\left((\mathbf{u}-\bar{\mathbf{u}})\cdot\nabla\bar{c}+\bar{c}\nabla\cdot(\mathbf{u}-\bar{\mathbf{u}})\right),
\end{align}
while
\begin{align}
    \psi(c-\bar{c})\nabla\cdot(\mathbf{v}-\bar{\mathbf{v}})=&~\psi(\omega-\bar{\omega})\nabla\cdot(\mathbf{v}-\bar{\mathbf{v}})+\psi\left(c-\omega-(\bar{c}-\bar{\omega})\right)\nabla\cdot(\mathbf{v}-\bar{\mathbf{v}})\nonumber\\
    =&~\nabla\cdot\left[\psi(\omega-\bar{\omega})(\mathbf{v}-\bar{\mathbf{v}})\right]-(\omega-\bar{\omega})(\mathbf{v}-\bar{\mathbf{v}})\cdot\nabla\psi-\psi(\mathbf{v}-\bar{\mathbf{v}})\cdot\nabla(\omega-\bar{\omega}).
\end{align}
Using the above identities, we thus obtain \eqref{eq:growth rate of c-barc} by an approximation $\psi\to \1_{[0,T]\times\Omega}$.
\end{proof}

Now, we can, after carefully bounding the terms in \eqref{eq:growth rate full relative energy} and \eqref{eq:growth rate of c-barc} and by  using   damping and bulk viscosity, apply Gronwall's lemma to the relative entropy between ${\bf U}$ and $\overline{\bf U}$.

\begin{Lemma}[Bound for the growth rate of the reduced relative energy]\label{lem:growth rate reduced relen}
     Let $\mathbf{U}$ be a weak solution to \eqref{eq:NSCHr} and let $\bar{\mathbf{U}}$ be a classical solution to \eqref{eq:residualNSCHr}.	Let $\varepsilon\in(0,1)^4$, $\tau\lesssim\beta\lesssim \delta$ and assume that $\bar{c}\in L^\infty((0,T)\times\Omega)$, uniformly in $\varepsilon$.
     Then there holds for all $t\leq T$,
	\begin{align}\label{eq:growth rate bound}
		\int\limits_\Omega \eta^\varepsilon(\mathbf{U},\bar{\mathbf{U}})(t,\mathbf{x})~d\mathbf{x}\leq&~ \int\limits_\Omega \eta^\varepsilon(\mathbf{U},\bar{\mathbf{U}})(0,\mathbf{x})~d\mathbf{x}+\int_0^tK\int\limits_\Omega \eta^\varepsilon(\mathbf{U},\bar{\mathbf{U}})(s,\mathbf{x})~d\mathbf{x}+\int\limits_\Omega\mathcal{R}~d\mathbf{x}\nonumber\\
        &-\nu|\mathbf{u}-\bar{\mathbf{u}}|_{H^1(\Omega)}^2-\frac{(\lambda+\nu)}{2}\left\|\nabla\cdot(\mathbf{u}-\bar{\mathbf{u}})\right\|^2_{L^2(\Omega)}-\frac{1}{2}\|\mathbf{v}-\bar{\mathbf{v}}\|_{L^2(\Omega)}^2\nonumber\\
		&-\frac{\delta}{2}\left\|\nabla\cdot(\mathbf{v}-\bar{\mathbf{v}})\right\|_{L^2(\Omega)}^2 - \frac{\tau}{2\beta}\|\nabla(c-\bar{c})\|^2_{L^2(\Omega)}~ds,
	\end{align}
	with the remainder
	\begin{align}\label{eq:residual}
		\mathcal{R}:=&~\alpha\mathcal{R}_1(p-\bar{p})+\mathcal{R}_2(\mathbf{u}-\bar{\mathbf{u}})+\frac{1}{\beta}\mathcal{R}_3(c-\omega-(\bar{c}-\bar{\omega}))+\delta\mathcal{R}_4(\mathbf{v}-\bar{\mathbf{v}})\nonumber\\
		&+\mathcal{R}_5(\omega_t-\bar{\omega}_t)+\mathcal{R}_3 (W''(\bar{c})+1)(c-\bar{c}),
	\end{align}
	and where the constant $K$ is such that,
	\begin{align}\label{eq:K scaling}
		K\lesssim&~(\lambda+\nu)^{-1}+\gamma^{-1}+\|\nabla\cdot(\bar{c}\bar{\mathbf{u}}+\bar{\mathbf{v}})\|_\infty+\|\bar{\mathbf{u}}\|^2_\infty+\|\nabla\cdot\bar{\mathbf{u}}\|_\infty+\|\nabla\bar{\mathbf{u}}+\nabla\bar{\mathbf{u}}^\intercal\|_\infty\nonumber\\
        &+\|\nabla\bar{c}\|_\infty+\|\bar{\mathbf{u}}\cdot\nabla\bar{c}\|_\infty+\beta\gamma\|\nabla\Delta\bar{\omega}\|_\infty.
	\end{align}
\end{Lemma}

\begin{rem}\label{remark 1.5}
	The a-priori assumption $\bar{c}\in L^\infty(0,T;\Omega)$ is used to simply bound $\|W''(\bar{c})\|_\infty$ and $\|W'''(\bar{c})\|_\infty$, uniformly in $\varepsilon$, and thus conveniently omit these prefactors in the description of $K$.
\end{rem}
\begin{rem}
	Only bulk viscosity has been used to absorb problematic terms. As such, Lemma \ref{lem:growth rate reduced relen} holds even in the case of vanishing shear viscosity, $\nu=0$, as long as $\lambda>0$ which only appears in the approximating system anyhow. This is reflected in the scaling of the constant $K$ in \eqref{eq:K scaling}. 
\end{rem}

We will now derive an a-priori error estimate of the difference between a solution to \eqref{eq:NSCHr} and a solution to \eqref{eq:NSCH} via Lemma \ref{lem:growth rate reduced relen}. To this end, we have to construct a suitable solution $\bar{\mathbf{U}}$ to \eqref{eq:residualNSCHr} from the exact solution to the NSCH system \eqref{eq:NSCH}. 

\begin{rem}\label{rem:reduced}
Note that if $\bar{\mathbf{U}}$ is chosen straightforwardly, the resulting convergence rate is reduced. Indeed, setting $\bar{p}:=\widetilde{p}$ and $\bar{\mathbf{u}}=\widetilde{\mathbf{u}}$ with $(\widetilde{p},\widetilde{\mathbf{u}},\widetilde{c})$ the exact solution to the NSCH system \eqref{eq:NSCH} would yield $\mathcal{R}_1=- \tilde p_t$.
This would lead to a convergence rate of order $1/2$ in $\alpha$ due to the term $\alpha \mathcal{R}_1 (p - \bar p)$ in \eqref{eq:growth rate full relative energy} and the fact that the relative energy only controls $\alpha (p - \tilde p)^2$, i.e., the control vanishes in the $\alpha \rightarrow 0$ limit.
\end{rem}

Instead, we do the following.

\begin{Proposition}\label{prop:high-order}
	Let $\varepsilon\in(0,1)^4$ and $\tau\lesssim \beta\lesssim \delta$. Let $(\widetilde{p},\widetilde{\mathbf{u}},\widetilde{c})$ denote a classical solution of the NSCH system \eqref{eq:NSCH} and let $\theta:(0,T) \times \Omega \rightarrow \mathbb{R}$ satisfy 
	\begin{align}\label{eq:theta}
		\Delta \theta(t,\cdot)=\widetilde{p}_t(t,\cdot),
	\end{align}
	and  $\int_\Omega \theta(t,\cdot) ~d\vecx=0$ for almost all $t \in (0,T)$.
    Set 
	\begin{align}\label{eq:choice of barU}
		\bar{p}:=\widetilde{p},\quad \bar{\mathbf{u}}:=\widetilde{\mathbf{u}}-\alpha\nabla\theta,\quad \bar{c}:=\widetilde{c},\quad \bar{\mathbf{v}}:=-\nabla\mu(\widetilde{c})+\alpha \widetilde{c}\,\nabla\theta+\tau\nabla\widetilde{c},
	\end{align} 
	and choose $\bar{\omega}$ as the solution to
    \begin{equation}\label{eq:elliptic reconstruction}
        -\gamma\Delta\bar{\omega}+\frac{1}{\beta}\bar{\omega}=\frac{1}{\beta}\tilde{c}.
    \end{equation}
    Let $\mathbf{U}$ be a weak solution to \eqref{eq:NSCHr} and let  $\bar{\mathbf{U}}$ be defined by \eqref{eq:choice of barU}-\eqref{eq:elliptic reconstruction}. Then, there exist constants $K,k,C>0$ such that for all times $t\leq T$, there holds
	\begin{align}\label{eq:growth rate high order}
		\int\limits_\Omega \eta^\varepsilon(t,\mathbf{x})~d\mathbf{x}\leq&~ \int\limits_\Omega \eta^\varepsilon(0,\mathbf{x})~d\mathbf{x} + \int_0^t K\int\limits_\Omega \eta^\varepsilon(s,\mathbf{x})~d\vecx+C\alpha^2\|\widetilde{p}_t\|^2_{H^1}\left(\|\widetilde{\mathbf{u}}\|^2_{H^1}+\alpha^2\|\widetilde{p}_t\|^2_{L^2}\right)\nonumber\\
		&+ C\alpha^2\|\widetilde{p}_t\|^2_{L^2}\left(\|\widetilde{\mathbf{u}}\|^2_\infty+\alpha^2|\widetilde{p}_t|^2_{H^1}+k(\lambda+\nu)+k\|\widetilde{c}\|^2_\infty\right)+C\alpha^2\nu\|\widetilde{p}_t\|^2_{H^1}\nonumber\\
		&+\beta^2\gamma^2\left(\left(\frac{k\|\widetilde{c}\|_\infty^2}{\lambda+\nu}+\|\nabla\widetilde{c}\|_\infty^2\right)|\widetilde{c}|^2_{H^4}+|\widetilde{c}|^2_{H^5}\right)+C\frac{k\alpha^2}{(\lambda+\nu)}\|\widetilde{p}_{tt}\|_{L^2}^2\nonumber\\
		&+ k\delta^2\left(|\mu(\widetilde{c})_t|^2_{H^1}+\tau^2|\widetilde{c}_t|^2_{H^1}+\tau^2|\widetilde{c}|^2_{H^3}+|\mu(\widetilde{c})|^2_{H^3}+C\alpha^2\|\nabla\widetilde{c}\|^2_\infty\|\widetilde{p}_t\|^2_{L^2}\right)\nonumber\\
		&+ Ck\alpha^2\delta^2\left(|\widetilde{p}_t|^2_{H^1}\left(|\widetilde{c}|^2_{H^2}+\|\widetilde{c}_t\|^2_{L^2}\right)+\|\widetilde{c}\|^2_\infty\left(|\widetilde{p}_t|^2_{H^1}+\|\widetilde{p}_{tt}\|^2_{L^2}\right)\right)+Ck\tau^2|\widetilde{c}|_{H^1}^2~ds.
	\end{align}
	The constant $K>0$ is as in \eqref{eq:K scaling}, $k>0$ is large enough to facilitate certain Young's inequalities and $C>0$ is a constant coming from elliptic estimates and Sobolev inequalities.
\end{Proposition}
\begin{proof}
By definition, $(\bar{p},\bar{\mathbf{u}},\bar{c},\bar{\mathbf{v}},\bar{\omega})$ is a solution to \eqref{eq:residualNSCHr} with residuals
\begin{align*}
	\mathcal{R}_1:=&~0,\\
	\mathcal{R}_2:=&~\alpha\nabla\theta_t+\alpha(\widetilde{\mathbf{u}}\cdot\nabla)\nabla\theta+\alpha(\nabla\theta\cdot\nabla)\bar{\mathbf{u}}+\frac{\alpha}{2}\bar{\mathbf{u}}\Delta\theta-\alpha\nu\Delta\nabla\theta-\alpha(\lambda+\nu)\nabla\Delta\theta+\gamma \bar{c}\nabla\Delta(\bar{c}-\bar{\omega}),\\
	\mathcal{R}_3:=&~0,\\
	\mathcal{R}_4:=&-\bar{\mathbf{v}}_t-\frac{\gamma}{\delta}\nabla\Delta(\bar{c}-\bar{\omega})-\frac{\alpha}{\delta}\bar{c}\nabla\theta-\frac{\tau}{\delta}\nabla\bar{c}-\nabla(\nabla\cdot\bar{\mathbf{v}}),\\
	\mathcal{R}_5:=&~0.
\end{align*}
We start from \eqref{eq:growth rate full relative energy} and note that the contribution of the second residual can be bounded as follows:
\begin{align}
	\int_\Omega \mathcal{R}_2(\mathbf{u}-\bar{\mathbf{u}})~d\vecx=& \int_\Omega-\alpha\theta_t\nabla\cdot(\mathbf{u}-\bar{\mathbf{u}}) -\alpha(\widetilde{\mathbf{u}}\cdot\nabla)\nabla\theta\cdot(\mathbf{u}-\bar{\mathbf{u}})-\alpha(\nabla\theta\cdot\nabla)\bar{\mathbf{u}}\cdot(\mathbf{u}-\bar{\mathbf{u}})\nonumber\\
	&-\frac{\alpha}{2}\bar{\mathbf{u}}\cdot(\mathbf{u}-\bar{\mathbf{u}})\Delta\theta+\alpha\nu\Delta\nabla\theta\cdot(\mathbf{u}-\bar{\mathbf{u}})+\alpha(\lambda+\nu)(\nabla\cdot(\mathbf{u}-\bar{\mathbf{u}}))\Delta\theta\nonumber\\
	&-\gamma\left(\bar{c}\Delta(\bar{c}-\bar{\omega})\right)\nabla\cdot(\mathbf{u}-\bar{\mathbf{u}})-\gamma\left(\nabla\bar{c}\cdot(\mathbf{u}-\bar{\mathbf{u}})\right)\Delta(\bar{c}-\bar{\omega})~d\vecx\nonumber\\
	\leq&~\frac{k\alpha^2}{2(\lambda+\nu)}\|\theta_t\|^2_{L^2}+\frac{\lambda+\nu}{2k}\|\nabla\cdot(\mathbf{u}-\bar{\mathbf{u}})\|^2_{L^2}+\frac{\alpha^2}{2}\|\widetilde{\mathbf{u}}\|^2_\infty|\theta|^2_{H^2}+\frac{1}{2}\|\mathbf{u}-\bar{\mathbf{u}}\|^2_{L^2}\nonumber\\
	&+\frac{\alpha^2}{2}\|\nabla\theta\|^2_\infty\|\bar{\mathbf{u}}\|^2_{H^1}+\frac{1}{2}\|\mathbf{u}-\bar{\mathbf{u}}\|^2_{L^2}+\frac{\alpha^2}{4}\|\bar{\mathbf{u}}\|^2_\infty|\theta|^2_{H^2}+\frac{1}{4}\|\mathbf{u}-\bar{\mathbf{u}}\|^2_{L^2}\nonumber\\
	&+\frac{\alpha^2\nu}{2}|\theta|^2_{H^3}+\frac{\nu}{2}\|\mathbf{u}-\bar{\mathbf{u}}\|^2_{L^2}+\frac{k\alpha^2(\lambda+\nu)}{2}|\theta|^2_{H^2}+\frac{\lambda+\nu}{2k}\|\nabla\cdot(\mathbf{u}-\bar{\mathbf{u}})\|^2_{L^2}\nonumber\\
	&+\frac{k\gamma^2}{2(\lambda+\nu)}\|\bar{c}\|^2_\infty|\bar{c}-\bar{w}|_{H^2}^2+\frac{\lambda+\nu}{2k}\|\nabla\cdot(\mathbf{u}-\bar{\mathbf{u}})\|^2_{L^2}+\frac{\gamma^2}{2}\|\nabla\bar{c}\|_\infty^2|\bar{c}-\bar{\omega}|_{H^2}^2\nonumber\\
	&+\frac{1}{2}\|\mathbf{u}-\bar{\mathbf{u}}\|_{L^2}^2,
\end{align}
where choosing $k>0$ large enough allows for absorption of the divergence terms via bulk viscosity. 

The other residual is bounded as follows, 
\begin{align}
	\int_\Omega\delta\mathcal{R}_4(\mathbf{v}-\bar{\mathbf{v}})~d\vecx\leq&~\frac{k\delta^2}{2}\|\bar{\mathbf{v}}_t\|_{L^2}^2+\frac{k\gamma^2}{2}|\bar{c}-\bar{\omega}|^2_{H^3}+\frac{k\alpha^2}{2}|\bar{c}|^2_\infty|\theta|_{H^1}^2\nonumber\\
	&+\frac{k\tau^2}{2}|\bar{c}|_{H^1}^2
	+\frac{k\delta^2}{2}|\bar{\mathbf{v}}|_{H^2}^2+\frac{5}{2k}\|\mathbf{v}-\bar{\mathbf{v}}\|^2_{L^2},
\end{align} 
where choosing $k>0$ large enough allows for absorption of the $\|\mathbf{v}-\bar{\mathbf{v}}\|^2_{L^2}$ terms via damping.

It remains to apply the elliptic estimates
\begin{align}\label{eq:jan elliptic estimate}
		|\bar{c}-\bar{w}|^2_{H^2(\Omega)}\leq \beta^2|\bar{c}|^2_{H^4(\Omega)},\qquad |\bar{c}-\bar{w}|^2_{H^3(\Omega)}\leq \beta^2|\bar{c}|^2_{H^5(\Omega)},
\end{align}
 that can be found in \cite{Giesselmann14}, and to note that on the one hand
\begin{align}
	\|\theta_t\|_{L^2}\lesssim\|\widetilde{p}_{tt}\|_{L^2},
\end{align}
while for $s>1/2$, one has the bound
\begin{align}
	\|\nabla\theta\|_\infty\lesssim\|\theta\|_{H^{2+s}}\lesssim\|\widetilde{p}_t\|_{H^s},
\end{align}
by Sobolev embedding and elliptic estimates.
\end{proof}

Finally, we obtain the following a-priori error estimate.

\begin{thm}\label{thm:high order convergence result}
	Let $\varepsilon\in(0,1)^4$, $\tau\lesssim \beta\lesssim \delta$, and $\mathbf{U}$ denote a weak solution to \eqref{eq:NSCHr}.
	 Let $(\widetilde{p},\widetilde{\mathbf{u}},\widetilde{c})$ denote a classical solution  to the NSCH system \eqref{eq:NSCH} complemented with periodic boundary conditions, such that
     \begin{align}
     \widetilde{p}\in&~ C^1\left([0,T];H^2(\Omega)\right)\cap C^2\left([0,T];L^2(\Omega)\right)\label{eq:norms for constant C 2},\\
        \widetilde{\mathbf{u}}, |\nabla\widetilde{\mathbf{u}}+\nabla\widetilde{\mathbf{u}}^\intercal|_F\in&~  L^\infty\left((0,T)\times\Omega\right),\label{eq:norms for constant C 1}\\
           \widetilde{c}\in&~ L^\infty\left(0,T;W^{4,\infty}(\Omega)\right)\cap L^\infty\left(0,T;H^5(\Omega)\right)\cap C^1\left([0,T];H^3(\Omega)\right),
       \end{align} 
       and choose $\widetilde{\omega}$ as the solution to \eqref{eq:elliptic reconstruction} with periodic boundary conditions.
	 
     Then, under the initial conditions
     \begin{align}
        &p(0,\cdot)=\widetilde{p}(0,\cdot)\in H^1(\Omega),\quad &\mathbf{u}(0,\cdot)=\widetilde{\mathbf{u}}(0,\cdot)\in W^{1,\infty}(\Omega;\R^d),\\
        &c(0,\cdot)=\widetilde{c}(0,\cdot)\in W^{4,\infty}(\Omega) \cap H^5(\Omega),\quad &\mathbf{v}(0,\cdot)=-\nabla\mu(c(0,\cdot)),
     \end{align}
	 there holds for all times $t\in(0,T)$,
\begin{multline}\label{eq:NSCH high order convergence bound}
\Biggl(
   \sqrt{\alpha}\|p-\widetilde{p}\|
   +\|c-\widetilde{c}\|
   +\|(c-\widetilde{c})^2\| 
   +\sqrt{\beta}\left\|\frac{c-\omega}{\beta}
     +\gamma\Delta\widetilde{c}\right\|
  \\
   +\|\mathbf{u}-\widetilde{\mathbf{u}}\| 
   +\sqrt{\delta}\left\|\mathbf{v}
     +\nabla\mu(\widetilde{c})\right\|
   +\sqrt{\gamma}\left\|\nabla\omega
     -\nabla\widetilde{c}\right\|
\Biggr)(t)\\
\leq{}
C\exp(tK)\left(
   \frac{1}{K}\left[
      \left(\frac{\alpha+\beta\gamma}{\sqrt{\lambda+\nu}}
      +\delta\right)
      +\alpha+\tau
   \right]
   +\alpha+\beta+\tau
\right)<\infty .
\end{multline}
	   with $\|\cdot\|:=\|\cdot\|_{L^2(\Omega)}$, and where $K>0$ is as in \eqref{eq:K scaling} and the constant $C>0$ depends solely on the norms \eqref{eq:norms for constant C 1}-\eqref{eq:norms for constant C 2} of the solution $[\widetilde{\mathbf{u}},\widetilde{c},\widetilde{p}]$ to the NSCH system \eqref{eq:NSCH}.
\end{thm}
\begin{proof}
    The proof is a mere application of Gronwall's Lemma to \eqref{eq:growth rate high order} in the context of Proposition \ref{prop:high-order} since, indeed, the left-hand side of \eqref{eq:NSCH high order convergence bound} is bounded by a multiple of $\sqrt{\eta^\varepsilon(\mathbf{U},\bar{\mathbf{U}})}+\alpha+\beta+\tau$ with $\bar{\mathbf{U}}$ as defined in Proposition \ref{prop:high-order} due to the regularity assumptions \eqref{eq:norms for constant C 1}-\eqref{eq:norms for constant C 2}. Furthermore, note that $K$ as defined in \eqref{eq:K scaling} is independent of $\varepsilon$ since, by definition, the relaxation parameters are bounded from above.
\end{proof}

\begin{rem}
    From \eqref{eq:NSCH high order convergence bound} it is evident that the relaxation system provides approximations of the pressure, the Laplacian of the phase field variable and the gradient of the chemical potential of the target system as
    \begin{align}
        \|p-\widetilde{p}\|_{L^\infty L^2}\leq&~C_0 \sqrt{\alpha},\quad \text{if}~\delta\lesssim\alpha,\\
        \left\|\frac{c-\omega}{\beta}+\gamma\Delta\widetilde{c}\right\|_{L^\infty L^2} \leq&~C_0\sqrt{\beta},\quad\text{if}~\alpha,\delta\lesssim\beta,\\
        \left\|\mathbf{v}+\nabla\mu(\widetilde{c})\right\|_{L^\infty L^2} \leq&~C_0\sqrt{\delta},\quad\text{if}~\alpha\lesssim\delta,
    \end{align}
    where the constant $C_0$ depends on $T,\gamma,\lambda,\nu$ and norms of the target solution as specified in Theorem \ref{thm:high order convergence result}.
\end{rem}

\begin{rem}[Convergence in physical domains]\label{rem:convergence in physical domains}
    By similar considerations, a convergence result in the spirit of Theorem \ref{thm:high order convergence result} can be proven in physical domains with convergence order $O(\alpha^{1/2}+\beta+\delta+\tau)$. The sub-optimal rate in $\alpha$ is due to a difficulty that arises via the no-flux boundary conditions for $\bar{\mathbf{u}}$. Using the same strategy as in the proof of \eqref{eq:growth rate high order} does not work in this case since additional residual terms appear that are concentrated on the boundary and cannot be absorbed.
\end{rem}

\begin{rem}[Euler-Cahn-Hilliard]
    Let us stress that Theorem \ref{thm:high order convergence result} remains valid in the inviscid regime $\nu=0$, as long as $\lambda>0$. Even more so, in the context of an Euler-Cahn-Hilliard system (i.e., referring to \eqref{eq:NSCH} with $\nu=0$), the $O(\alpha+\beta+\delta+\tau)$ convergence result can be generalized
    to physical domains if the impermeability condition ${\mathbf{u}\cdot\mathbf{n}=0}$ for the velocity fields is prescribed on the boundary. This result is validated by numerical tests in Section \ref{sec:numerics}.
\end{rem}

\section{Numerical justification of the relaxation approximation}\label{sec:numerics}
In this section, we first introduce tailored finite-difference methods for the NSCH system and the relaxation approximations. 
Then, numerical experiments are presented that confirm the convergence result from Theorem \ref{thm:high order convergence result} obtained in Section \ref{sec:analysis}.
For this, simulations in one and two spatial dimensions are performed with a variety of initial data including a single droplet which relaxes into its equilibrium, the Ostwald ripening of two bubbles, the merging of two bubbles and a droplet-wall interaction.
Besides the convergence of the relaxation formulation of the NSCH system \eqref{eq:NSCHr} to the original NSCH system \eqref{eq:NSCH} for $|\varepsilon| \rightarrow 0$, the effect of the grid resolution and time step size on this convergence is investigated.
However, only the inviscid limit systems of \eqref{eq:NSCH} and \eqref{eq:NSCHr} with $\visc = 0$ are considered numerically.
The structure of this section is as follows: First, in Sections \ref{sec:numNSCH}, \ref{sec:numNSCHr}, the numerical discretization of the original NSCH system \eqref{eq:NSCH} and its relaxation formulation \eqref{eq:NSCHr} are outlined, respectively.
This is followed by numerical experiments in one dimension which are utilized to investigate the convergence of \eqref{eq:NSCHr} to \eqref{eq:NSCH} for $|\varepsilon| \rightarrow 0$.
In Section \ref{sec:num1D}, beyond the convergence for $|\varepsilon| \rightarrow 0$ at a fixed grid resolution and time step size, its sensitivity with respect to the grid resolution and time step size is examined.
Finally, in Section \ref{sec:num2D}, two-dimensional convergence results are presented.
\subsection{A numerical approach for the NSCH system and the relaxation  approximation\label{sec:num}}
\subsubsection{A finite-difference MAC method  for the NSCH system} \label{sec:numNSCH}
For the numerical discretization of NSCH system \eqref{eq:NSCH}, a second-order accurate conservative finite-difference scheme in space is utilized following ideas from \cite{Shu2006,Dhaouadi24}.
Advancement in time is realized by the use of a first-order accurate implicit method and a constant time step size $\Dt$.
For this, the domain $\domain = (\xmin,\xmax) \times (\ymin,\ymax)$ is tessellated into a Cartesian equispaced marker-and-cell (MAC) grid with $\nCells = \Nx \times \Ny$ elements.
Each element $\domainCell$ has a volume given as $\vol{\domainCell} = \Dx \Dy$ with $\xindex = 1, \dots , \Nx$ and $\yindex = 1, \dots , \Ny$.
The length of the cell edges in the $x$- and $y$-direction are defined as $\Dx = \frac{\xmax-\xmin}{\Nx}$ and $\Dy = \frac{\ymax-\ymin}{\Ny}$, respectively.
As common for a MAC grid, the \textit{thermodynamic} variables, i.e. the phase field variable $\phas$ and the pressure $\pres$, are defined at the cell centers $\Vecx_{\xindex,\yindex} = (x_\xindex,y_\yindex)$ of each element $\domainCell$.
\begin{figure}[htpb]
  \centering
  \begin{subfigure}[t]{0.45\textwidth}
        \centering
        \includegraphics[width=\textwidth]{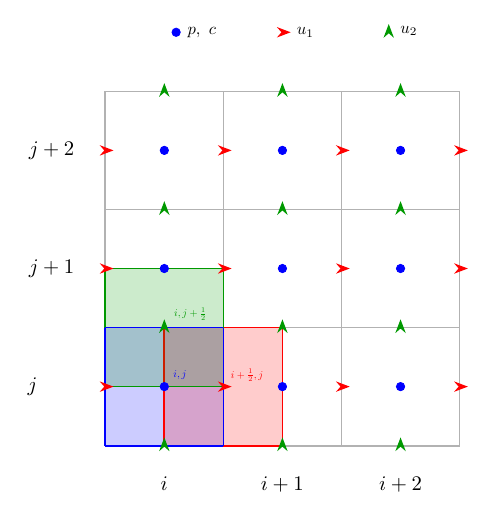}
        \caption{Marker-and-cell grid utilized for the discretization of the NSCH system \eqref{eq:NSCH}.}
        \label{fig:MACgridNSCH}
  \end{subfigure}
  \hfill
  \begin{subfigure}[t]{0.45\textwidth}
        \centering
        \includegraphics[width=\textwidth]{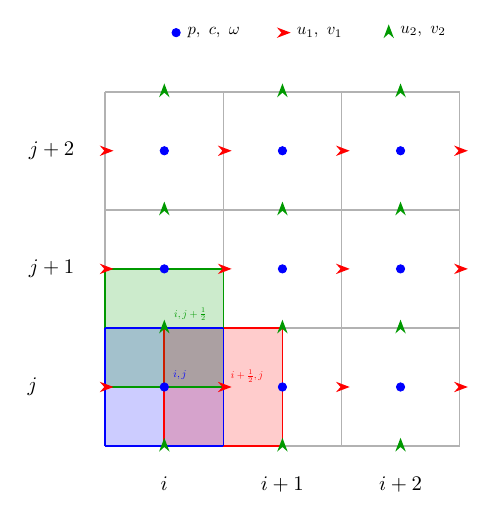}
        \caption{Marker-and-cell grid utilized for the discretization of the relaxation approximation of the NSCH system \eqref{eq:NSCHr}.}
        \label{fig:MACgridNSCHr}
  \end{subfigure}
  \caption{Comparison of the spatial distribution of the solution variables for the NSCH system \eqref{eq:NSCH} and its relaxation approximation\eqref{eq:NSCHr}.}
  \label{fig:MACgrid}
\end{figure}
In contrast, the \textit{hydrodynamic} variables which means the velocities in the $x$ and $y$-directions, $\velv = (\velx,\vely)$, are defined on a staggered grid with the cell centers of these corresponding elements $\domainCellx$ and $\domainCelly$ defined on the edges of the elements $\domainCell$, see  Fig.~\ref{fig:MACgridNSCH} for a graphical illustration.
\newline
The semi-discretization of the governing equations \eqref{eq:NSCH} in time reads as
\begin{subequations}
  \label{eq:NSCH_approx}
  \begin{align}
    \label{eq:NSCHphasstar_approx}
    c^{\star} + \Dt \Big[ \div \left( c^{\star} \velv^{n} \right) - \div \left( \DoubWellPot^{\prime\prime} (c^{n}) \grad c^{\star} - \gamma \grad \lap c^{\star} \right) \Big] =& ~c^{n}, \\
    \label{eq:NSCHvelstar_approx}
    \velv^{\star} + \Dt ~ \div \left( \velv^{\star} \otimes \velv^{n} \right) =& ~\velv^{n} - \Dt ~ c^{\star} \Big[ \DoubWellPot^{\prime\prime} (c^{n}) \grad c^{\star} - \gamma \grad \lap c^{\star} \Big], \\
    \label{eq:NSCHpres_approx}
    \Dt~\lap \pres^{n+1} =& ~\div \velv^{\star}, \\
    \label{eq:NSCHvel_approx}
    \velv^{n+1} =& ~\velv^{\star} - \Dt~\div \pres^{n+1}, \\
    \label{eq:NSCHphas_approx}
    c^{n+1} + \Dt \Big[ \div \left( c \velv \right)^{n+1} - \div \left( \DoubWellPot^{\prime\prime} (c^{\star}) \grad c^{n+1} - \gamma \grad \lap c^{n+1} \right) \Big] =& ~c^{n}.
  \end{align}
\end{subequations}
The scheme in \eqref{eq:NSCH_approx} is a straightforward extension of a classical second-order accurate semi-implicit conservative finite-difference scheme for the inviscid case, i.e., \eqref{eq:NSCH} with $\nu =0$.
Herein, the projection method proposed in \cite{Chorin1967} which is based on a Helmholtz-Hodge decomposition of the velocity field is used in \eqref{eq:NSCHpres_approx} and \eqref{eq:NSCHvel_approx} to guarantee the divergence-free condition of the velocity field.
Due to the semi-implicit discretization, the equations in \eqref{eq:NSCH_approx} can be solved independently in the order as they appear in each time step.
Moreover, for each solution variable only a linear system has to be solved which can be easily realized by the use of the GMRES algorithm \cite{Saad1986}.
A detailed discussion of the utilized stencils for the approximation of the spatial derivatives in \eqref{eq:NSCH_approx} is given in the Appendix  \ref{appx:NSCH_FD}.
\subsubsection{A finite-difference MAC method for the  relaxation approximation} \label{sec:numNSCHr}
Similar to the discretization of the original NSCH system, also for the relaxation formulation of the NSCH system, a second-order accurate semi-implicit conservative finite-difference scheme is used.
However, as outlined in Section \ref{subsec:approx}, system \eqref{eq:NSCHr} has several modifications in comparison to the relaxation formulation originally proposed in \cite{KKR25}.
In order to demonstrate that the convergence result of Theorem \ref{thm:high order convergence result} can be numerically reproduced also without these modifications, the relaxation formulation as proposed in \cite{KKR25} is considered in the following.
Consequently, the investigated system is given as
\begin{align}
	\pres_t^\numepsilon +  \frac{1}{\alpha} \div \velv^\numepsilon =& ~ 0, \nonumber \\
  \velv_t^\numepsilon +  (\velv^\numepsilon \cdot \grad) \velv^\numepsilon + \frac{1}{2} (\div \velv^\numepsilon) \velv^\numepsilon + \grad \pres^\numepsilon =& ~ -\phas^\numepsilon \grad \left( \DoubWellPot^{\prime} (\phas^\numepsilon) + \frac{1}{\beta} ( \phas^\numepsilon - \omega^\numepsilon ) \right), \nonumber \\
  \label{eq:NSCHr_num}
	\phas_t^\numepsilon +  \div (\phas^\numepsilon \velv^\numepsilon) + \div \vecmu^\numepsilon =& ~ 0, \\
	\vecmu_t^\numepsilon + \frac{1}{\delta} \grad \left( \DoubWellPot^{\prime} (\phas^\numepsilon) + \frac{1}{\beta} \phas^\numepsilon \right) =& ~ -\frac{\vecmu^\numepsilon}{\delta} + \frac{1}{\delta \beta} \grad \omega^\numepsilon, \nonumber \\
	-\gamma \lap \ordp^\numepsilon + \frac{1}{\beta} \ordp^\numepsilon =& ~ \frac{1}{\beta} \phas^\numepsilon. \nonumber
\end{align}
For the sake of clarity, the index $\numepsilon = (\alpha, \beta, \delta) \in (0,1)^3$ will be omitted in the following if there is no risk of confusion with  the unknows of the  limit system \eqref{eq:NSCH}.
The semi-discretization of the governing equations \eqref{eq:NSCHr_num} is done in a similar manner as in \eqref{eq:NSCH_approx}:
First, the phase field variable $\phas^\numepsilon$ is evolved in time.
The corresponding semi-discrete evolution equations read as
\begin{subequations}
  \label{eq:NSCHr_approx}
  \begin{align}
    \label{eq:NSCHr_cstarOne}
    \phas^{\star} + \Dt~\Big[ \div (\phas^{\star} \velv^{n}) + \div \vecmu^{\star} \Big] =& ~\phas^{n} , \\
    \label{eq:NSCHr_mustar}
    \vecmu^{\star} + \frac{\Dt}{\delta + \Dt}~\Big[ \DoubWellPot^{\prime \prime} (\phas^{n}) \grad \phas^{\star} + \frac{1}{\beta} \grad \left( \phas^{\star} - \ordp^{\star} \right) \Big] =& ~\frac{\delta}{\delta + \Dt} \vecmu^{n}, \\
    \label{eq:NSCHr_omegastar}
    {\color{white}{\Big[}} ( \mathcal{I} -\gamma \beta \lap ) \ordp^{\star} = ~\OperatorORDP \ordp^{\star} =& ~ \phas^{\star} {\color{white}{\Big]}},
  \end{align}
\end{subequations}
where $\OperatorORDP$ defines an operator given as $\OperatorORDP =  \mathcal I - \gamma \beta \lap$, $\mathcal I$ being the identity operator.
From \eqref{eq:NSCHr_approx} it is obvious that for the evolution of the phase field variable, a system of three equations has to be solved.
In the numerical implementation, we solve for $\phas^{\star}$ by substituting \eqref{eq:NSCHr_mustar} and \eqref{eq:NSCHr_omegastar} into \eqref{eq:NSCHr_cstarOne}.
The resulting update formula for $\phas^{\star}$ is then given by
\begin{align}
  \label{eq:NSCHr_cstarTwo}
  \phas^{\star} + \Dt~\left[ \div (\phas^{\star} \velv^{n}) - \frac{\Dt}{\delta + \Dt}~\div \Big[ \DoubWellPot^{\prime \prime} (\phas^{n}) \grad \phas^{\star} + \frac{1}{\beta} \grad \left( \phas^{\star} - \OperatorORDP^{-1} \phas^{\star} \right) \Big] \right] =& ~\phas^{n} - \frac{\delta \Dt}{\delta + \Dt} \div \vecmu^{n}.
\end{align}
While on the continuous level, the following equalities hold
\begin{align}
  \label{eq:omegaelliptic}
  \frac{1}{\beta} \left( \phas^{\star} - \OperatorORDP^{-1} \phas^{\star} \right) = \frac{1}{\beta} \left( \phas^{\star} - \ordp^{\star} \right) = \frac{1}{\beta} \left( \OperatorORDP \ordp^{\star} - \ordp^{\star} \right) = - \gamma \lap \ordp^{\star} = - \gamma \lap \OperatorORDP^{-1} \phas^{\star},
\end{align}
our numerical experiments indicate that the discretized version of the very first term in \eqref{eq:omegaelliptic} is prone to machine accuracy for $\beta \leq 10^{-7}$.
This limits the overall convergence of the relaxation system \eqref{eq:NSCHr_num} to the system \eqref{eq:NSCH} for $|\numepsilon| \rightarrow 0$.
The reason for this is that the discrete version of $\OperatorORDP^{-1}$ is only accurate up to machine precision.
Consequently, the discretized version of \eqref{eq:omegaelliptic} reads as
\begin{align}
  \label{eq:errorMP}
  \frac{1}{\beta} \left( \phas^{\star} - \widetilde{\OperatorORDP}^{-1} \phas^{\star} \right) = - \gamma \left( \widetilde{\lap} + \frac{\cR_M}{\beta} \right) \ordp^{\star} + \cR_M,
\end{align}
where the superscript $\widetilde{\boldsymbol{\cdot}}$ indicates the discrete approximation of an operator and $\varepsilon_M \coloneqq \nicefrac{\cR_M}{\beta}$ is an error term with $\cR_M \sim \cO ( 10^{-16} )$ the error due to machine precision on a 64-bit architecture.
Since $\cR_M$ is approximately a constant, with decreasing $\beta$, the error term $\varepsilon_M$ will eventually dominate the overall error and prevents convergence.
To overcome this shortcoming and enable a sound convergence of the relaxation formulation \eqref{eq:NSCHr_num}, instead of the very first term in \eqref{eq:omegaelliptic}, the very last one is used.
Consequently, the final update formula for the phase field variable reads as
\begin{align}
  \label{eq:NSCHr_cstarFinal}
  \phas^{\star} + \Dt~\left[ \div (\phas^{\star} \velv^{n}) - \frac{\Dt}{\delta + \Dt}~\div \Big[ \DoubWellPot^{\prime \prime} (\phas^{n}) \grad \phas^{\star} - \grad \lap \OperatorORDP^{-1} \phas^{\star} \Big] \right] =& ~\phas^{n} - \frac{\delta \Dt}{\delta + \Dt} \div \vecmu^{n}.
\end{align}
Therefore, this discretization resembles a fourth-order spatial operator similar to the discretized version of the original formulation of the NSCH system \eqref{eq:NSCHr_num}.
However, given that the purpose of the present work is to investigate the convergence of  solutions for system \eqref{eq:NSCHr_num} to those of \eqref{eq:NSCH} in the limit of $|\numepsilon| \rightarrow 0$, which would not be possible with \eqref{eq:NSCHr_cstarTwo} for arbitrarily small values of $\beta$, this increased stencil size is accepted as a drawback of the chosen discretization.
After the evolution of the phase field variable in time, the second step in the algorithm involves the update of the order parameter by solving \eqref{eq:NSCHr_omegastar} for $\ordp^{\star}$.
In the third step, the velocities are evolved in time.
For this, similar to the discretization of the original NSCH system, first a predictor step, given as
\begin{align}
  \label{eq:NSCHrvelstar_approx}
  \velv^{\star} + \Dt ~ \div \left( \velv^{\star} \otimes \velv^{n} \right) =& ~\velv^{n} - \Dt ~ c^{\star} \Big[ \DoubWellPot^{\prime\prime} (c^{n}) \grad c^{\star} - \gamma \grad \lap \ordp^{\star} \Big],
\end{align}
is performed.
Note that in \eqref{eq:NSCHrvelstar_approx}, the additional term $\frac{1}{2}(\div \velv) \velv$ has been omitted and will be accounted for in a later step.
Second, following \cite{Kwatra2009}, the pressure is evolved in time by combining
\begin{subequations}
  \begin{align}
    \alpha \pres^{n+1} =& ~\alpha \pres^{n} -  \Dt~\div \velv^{\star \star} , \\
    \label{eq:NSCHrvelstarstar}
    \velv^{\star \star} =& ~\velv^{\star} - \Dt~\grad \pres^{n+1}.
  \end{align}
\end{subequations}
This yields the following diffusion equation with a source term for the pressure
\begin{align}
  \label{eq:NSCHr_pres_approx}
  \alpha \pres^{n+1} - \Dt^2 \lap \pres^{n+1} =& ~\alpha \pres^{n} -  \Dt~\div \velv^{\star}.
\end{align}
Finally, the velocities are updated to the next time level by first applying \eqref{eq:NSCHrvelstarstar} and then accounting for the missing term $\frac{1}{2}(\div \velv) \velv$ by the use of
\begin{align}
  \label{eq:NSCHr_vel_approx}
  \velv^{n+1} =& ~\velv^{\star\star} - \Dt~\frac{1}{2} (\div \velv^{\star\star}) \velv^{\star\star}.
\end{align}
The final step involves the computation of the phase field variable at the next time step by
\begin{align}
  \phas^{n+1} + \Dt~\left[ \div (\phas \velv)^{n+1} - \frac{\Dt}{\delta + \Dt}~\div \Big[ \DoubWellPot^{\prime \prime} (\phas^{\star}) \grad \phas^{n+1} - \grad \lap \OperatorORDP^{-1} \phas^{n+1} \Big] \right] = ~\phas^{n} - \frac{\delta \Dt}{\delta + \Dt} \div \vecmu^{n},
\end{align}
which is then followed by an update of $\ordp$ and $\vecmu$ obtained by the subsequent evaluation of
\begin{subequations}
  \label{eq:NSCHr_omega_mu}
  \begin{align}
    \label{eq:ordp_initial}
    {\color{white}{\Big[}} \ordp^{n+1} =& ~ \OperatorORDP^{-1} \phas^{n+1} {\color{white}{\Big]}}, \\
    \label{eq:mu_initial}
    \vecmu^{n+1} =& ~\frac{\delta}{\delta + \Dt} \vecmu^{n} - \frac{\Dt}{\delta + \Dt}~\Big[ \DoubWellPot^{\prime \prime} (\phas^{\star}) \grad \phas^{n+1} + \frac{1}{\beta} \grad \lap \ordp^{n+1} \Big].
  \end{align}
\end{subequations}
Summarizing, the evolution of the numerical approximation of \eqref{eq:NSCHr_num} involves the subsequent evaluation of \eqref{eq:NSCHr_cstarFinal}, \eqref{eq:NSCHrvelstar_approx}, \eqref{eq:NSCHr_pres_approx}, \eqref{eq:NSCHr_vel_approx} and \eqref{eq:NSCHr_omega_mu}.
A detailed discussion of the used stencils for the approximation of the spatial derivatives in \eqref{eq:NSCHr_num} is given in Appendix \ref{appx:NSCHr_FD}.
\subsection{Numerical experiments in 1D} \label{sec:num1D}
\subsubsection{1D Ostwald ripening}
\begin{figure}[htpb]
  \centering
  \includegraphics[width=0.9\linewidth]{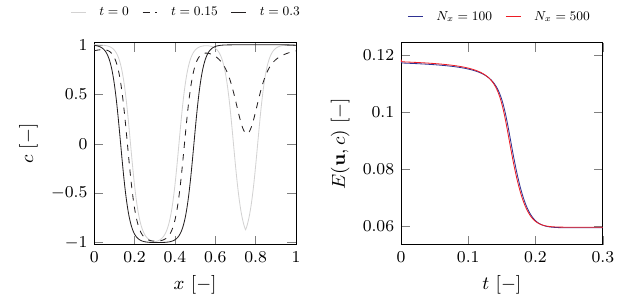}
  \caption{Left: Solution of the phase field variable at different time instances $t = 0$, $t = 0.15$, $t = 0.3$ with $\Nx = 500$. Right: Temporal evolution of the total Helmholtz energy \eqref{eq:energy NSCH} for different grid resolutions.}
  \label{fig:Ostwald1D}
\end{figure}

\begin{figure}[htpb]
  \centering
  \includegraphics[width=0.9\linewidth]{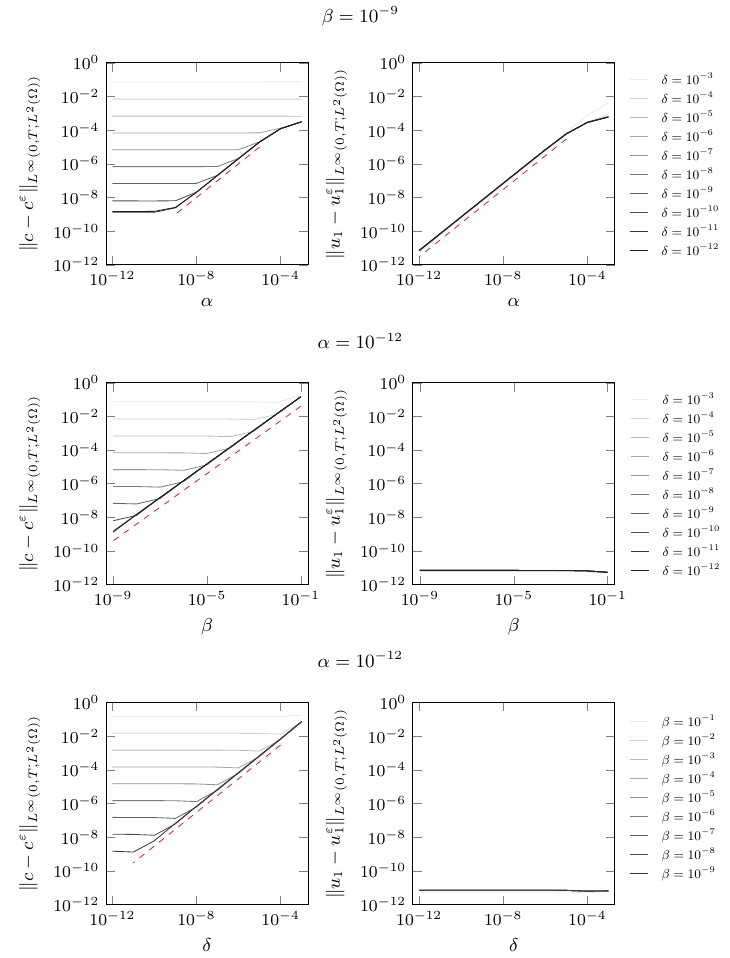}
  \caption{Convergence of the relaxation system \eqref{eq:NSCHr_num} to the NSCH system \eqref{eq:NSCH} for the Ostwald ripening test case with $\Nx = 100$. In the top row, the convergence in terms of $\alpha$ for a fixed value of $\beta$ and varying $\delta$ is depicted. The center row shows the convergence in terms of $\beta$ for a fixed $\alpha$ and varying $\delta$. In the bottom row, the convergence error in terms of $\delta$ for a fixed value of $\alpha$ and varying $\delta$ is plotted. The red dashed lines indicate convergence with order $\cO(\alpha)$, $\cO(\beta)$ and $\cO(\delta)$, depending on the respective running variable on the abscissa.}
  \label{fig:Error_Ostwald1D_N100}
\end{figure}

For the numerical investigation of the convergence behavior of \eqref{eq:NSCHr_num} to \eqref{eq:NSCH} in the limit $|\numepsilon| \rightarrow 0$ in one spatial dimension, the Ostwald ripening of two bubbles with initial data
\begin{align}
  (\phas, \velx )\transpose (0,x) = ( A(x) , 0 )\transpose,
\end{align}
where
\begin{align}
  A (x) = -1 + \sum_{i = 1}^2 \tanh \left( \frac{| x - x_i | - r_i}{\sqrt{2 \gamma}} \right)
\end{align}
and an end time of $t = 0.3$ was simulated.
The setup for the Ostwald ripening is adopted from \cite{Dhaouadi24,KKR25}, where the initial droplet's positions and radii are given as $x_1 = 0.3$, $x_2 = 0.75$ and $r_1 = 0.12$, $r_2 = 0.06$.
For the initialization of the evolution variables of the relaxation system \eqref{eq:NSCHr}, we used $\phas^\numepsilon(0,x) = \phas(0,x)$, $\velv^\numepsilon(0,x) = \velv(0,x)$ and $\pres^\numepsilon(0,x) = 0$.
The relaxation variable $\ordp^\numepsilon$ is initialized by solving the corresponding elliptic constraint based on $\phas^\numepsilon(0,x)$, see e.g. \eqref{eq:ordp_initial}, and the flux variable $\vecmu^\numepsilon$ was initialized as $\vecmu^\numepsilon = - \grad \chempot$ by discretely approximating the gradient of the chemical potential with the spatial discretization introduced in Appendix \ref{appx:NSCH_FD}.
A periodic computational domain 
$\Omega = (0,1)$ and a capillary parameter $\gamma = 10^{-3}$ were chosen.
The simulations were performed for different grid resolutions $\Nx \in \lbrace 100,~500 \rbrace$ and permutations of the relaxation parameters $\alpha \in [10^{-3},10^{-12}]$, $\beta \in [10^{-1},10^{-9}]$ and $\delta \in [10^{-3},10^{-12}]$ resulting in a total amount of 1802 simulations.
The time step size $\Dt = 10^{-3}$ was constant.
\newline
The solution of the phase field variable at different time instances and the temporal evolution of the total Helmholtz energy over time are depicted in Fig.~\ref{fig:Ostwald1D}.
Over time, the expected vanishing of the smaller bubble on  the right and the corresponding growth of the bubble on the left can be observed.
The thermodynamic admissibility of the results are indicated by a monotonic decreasing total Helmholtz energy over time.
A further refinement of the spatial resolution did not affect the results.
The error of the discretized equations \eqref{eq:NSCHr_semidiscrete} in comparison to the discretized equations \eqref{eq:NSCH_approx} for different choices of the  approximation  parameters in $\numepsilon$ and the corresponding convergence behavior is depicted in Fig.~\ref{fig:Error_Ostwald1D_N100}.
In the left and right columns, the error in the phase field variable and the velocity are plotted, respectively.
In the top row, the convergence with respect to $\alpha$ for a fixed value of $\beta = 10^{-9}$ and decreasing values of $\delta \in [10^{-3},10^{-12}]$ is presented.
The center row highlights the convergence with respect to $\beta$ for a fixed value of $\alpha = 10^{-12}$ and decreasing values of $\delta \in [10^{-3},10^{-12}]$.
In the bottom row, the convergence with respect to $\delta$ for a fixed value of $\alpha = 10^{-12}$ and decreasing values of $\beta \in [10^{-1},10^{-9}]$ is demonstrated.
For all solution variables, the expected convergence of order $\cO (\alpha + \beta + \delta )$ is achieved.
Exemplarily for the top left plot:
For sufficiently small values of $\beta$ and $\delta$, an error convergence of order $\cO ( \alpha )$ is observed, indicated by the red dashed line.
However, if $\delta$ is not sufficiently small, the error stagnates even for decreasing values of $\alpha$.
The stagnant error convergence observed in the top left and bottom left plots, where the error does not decrease with $\alpha$ even with decreasing values of $\delta$, can be attributed to a predominant error due to the relaxation parameter $\beta$.
A further reduction of $\beta$ would resolve this problem.
Due to the one-dimensional setup, the error in the velocity is mainly governed by the choice of $\alpha$, i.e., for sufficiently small values of $\alpha$, the error in $\velx$ is unaffected by the choice of $\beta$ and $\delta$.
Additional plots, investigating the convergence behavior also for the grid resolution $\Nx=500$, are given in Appendix \ref{appx:1D} in Fig. \ref{fig:Error_Ostwald1D_N500}.
However, different grid resolutions as well as time step sizes only affect the absolute error, not the convergence behavior.
%
%
%
%
\subsection{Numerical experiments in 2D} \label{sec:num2D}
\begin{figure}[htpb]
  \centering
  \begin{subfigure}[t]{0.32\textwidth}
        \centering
        \includegraphics[width=\textwidth]{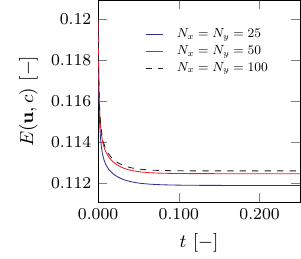}
  \end{subfigure}
  \begin{subfigure}[t]{0.32\textwidth}
        \centering
        \includegraphics[width=\textwidth]{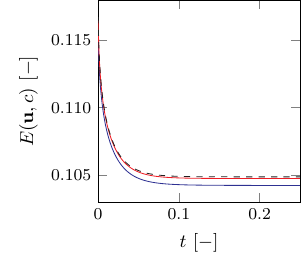}
  \end{subfigure}
  \begin{subfigure}[t]{0.32\textwidth}
        \centering
        \includegraphics[width=\textwidth]{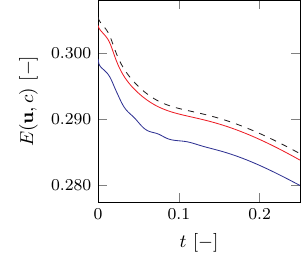}
  \end{subfigure}
  \caption{Temporal evolution of the total Helmholtz energy \eqref{eq:energy NSCH} of the NSCH system \eqref{eq:NSCH} for the single bubble, the merging droplets and the colliding droplets test cases, from left to right, each for three different grid resolutions.}
  \label{fig:Energy2D}
\end{figure}

For the numerical investigation of the convergence behavior of \eqref{eq:NSCHr_num} to \eqref{eq:NSCH} in the limit $|\numepsilon| \rightarrow 0$ in two spatial dimensions, the relaxation of a bubble to its equilibrium with the initial data
\begin{align}
  (\phas, \velx, \vely )\transpose (0,\mathbf{x}) = ( B (\mathbf{x}) , 0 , 0 )\transpose,
\end{align}
where
\begin{align}
  B (\mathbf{x}) = \begin{cases} - \cos (2 \pi r ) &: r = \sqrt{{(x-0.5)}^2 + {(y-0.5)}^2} \leq 0.5 \\
                        1 &: \mathrm{else}
          \end{cases}
\end{align}
and an end time of $t = 0.25$ and the merging of two droplets with initial data
\begin{align}
  (\phas, \velx, \vely )\transpose (0,\mathbf{x}) = ( C(\mathbf{x}) , 0 , 0 )\transpose,
\end{align}
where
\begin{align}
  \label{eq:MergingDroplets}
  C (\mathbf{x}) = -1 + \sum_{i=1}^2 \left( - \tanh \left( \frac{r - r_i}{\sqrt{2 \gamma}} \right) + \tanh \left( \frac{r + r_i}{\sqrt{2 \gamma}} \right) \right)
\end{align}
and an end time of $t = 0.25$ were simulated.
In \eqref{eq:MergingDroplets}, the symbol $r = \sqrt{{(x-x_i)}^2 + {(y-y_i)}^2}$ indicates the Euclidean distance to the droplet's center with the droplets' centers and radii given as $\mathbf{x}_1 = (0.4,0.5)\transpose,~\mathbf{x}_2 = (0.7,0.5)\transpose$ and $r_1 = 0.2,~r_2=0.1$, respectively.
In a final test case, a binary droplet collision with initial data
\begin{align}
  (\phas, \velx, \vely )\transpose (0,\mathbf{x}) = ( C(\mathbf{x}) , D_1(\mathbf{x}) , D_2(\mathbf{x}) )\transpose,
\end{align}
was investigated, where the droplets' centers and radii were defined as $\mathbf{x}_1 = (0.5,0.7)\transpose,~\mathbf{x}_2 = (0.5,0.3)\transpose$ and $r_1 = r_2=0.15$, respectively.
The divergence free velocity field was initialized as
\begin{align}
  D_1 (\mathbf{x}) =& \sin( 2 \pi x ) \cos( 2 \pi y ), \\
  D_2 (\mathbf{x}) =& \cos( 2 \pi x ) \sin( 2 \pi y )
\end{align}
and the simulation was performed until $t = 0.25$.
The initialization of the evolution variables of the relaxation system \eqref{eq:NSCHr} in 2D follows the same procedure as in the 1D case.
For all three numerical experiments, a periodic computational domain  $\Omega = (0,1)^2$ and a capillary parameter $\gamma = 6 \cdot 10^{-3}$ was chosen, if not stated otherwise.
The simulations were performed for two grid resolutions $\Nx = \Ny \in \lbrace 25,~50 \rbrace$ and permutations of the relaxation parameters $\alpha \in [10^{-3},10^{-12}]$, $\beta \in [10^{-1},10^{-9}]$ and $\delta \in [10^{-3},10^{-12}]$ resulting in a total amount of 5406 simulations.
The time step size $\Dt = 10^{-3}$ was constant.
\newline
The temporal evolution of the total Helmholtz energy according to the definition in \eqref{eq:energy NSCH} for all three test cases is depicted in Fig.~\ref{fig:Energy2D}.
The thermodynamic admissibility of the results is indicated by the monotonically decreasing total Helmholtz energy for the investigated grid resolutions $\Nx = \Ny = 25$ and $\Nx = \Ny = 50$.
In addition, the evolution of the total Helmholtz energy is shown for a grid resolution of $\Nx = \Ny = 100$ to demonstrate the grid convergence.
However, due to the corresponding computational costs and the amount of simulations, this grid resolution is not used for the further analysis of the convergence in terms of $\numepsilon$.
%
%
%
%
\subsubsection{Evolution of a 2D bubble}
\begin{figure}[htpb]
  \centering
  \begin{subfigure}[t]{0.32\textwidth}
        \centering
        \includegraphics[height=5cm]{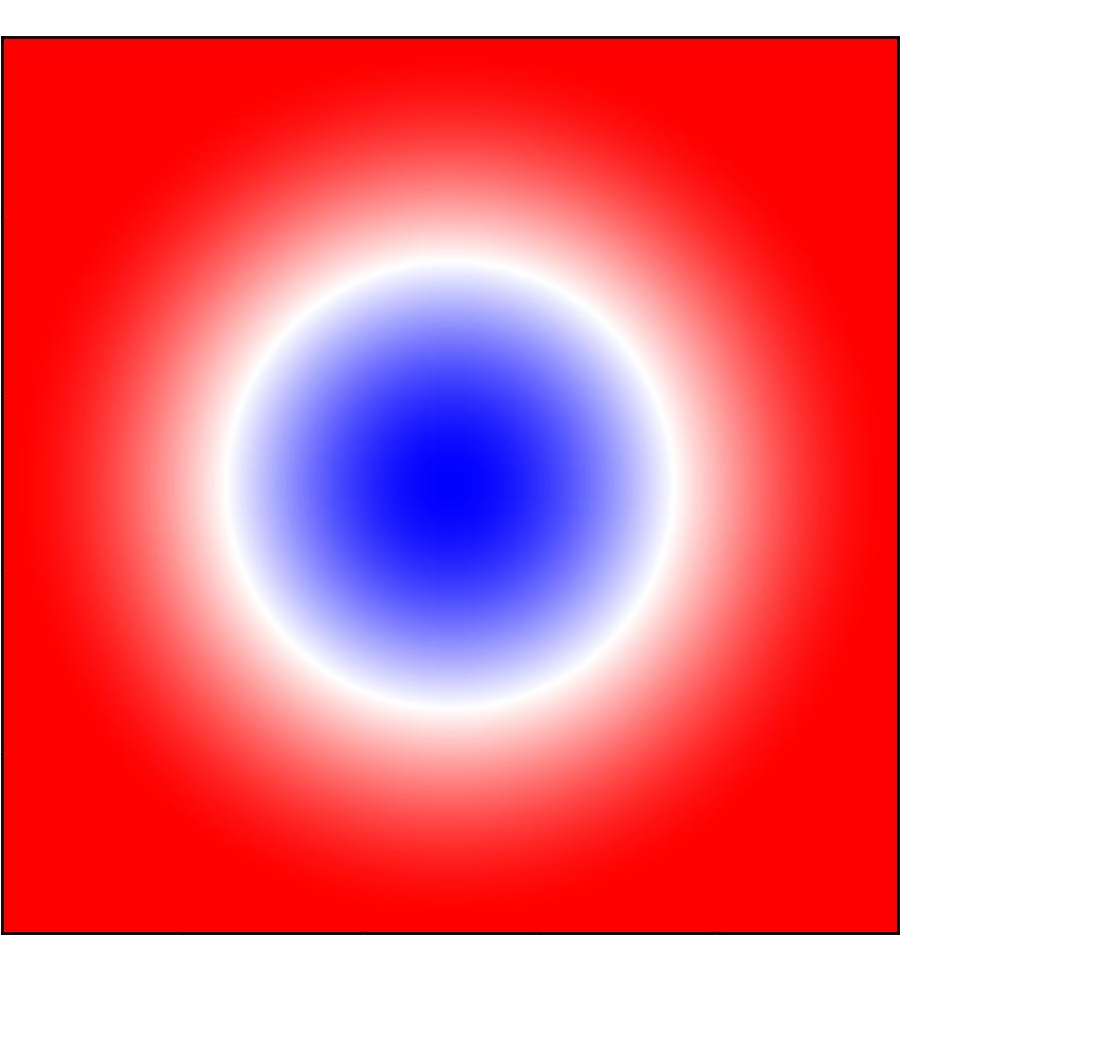}
  \end{subfigure}
  \begin{subfigure}[t]{0.32\textwidth}
        \centering
        \includegraphics[height=5cm]{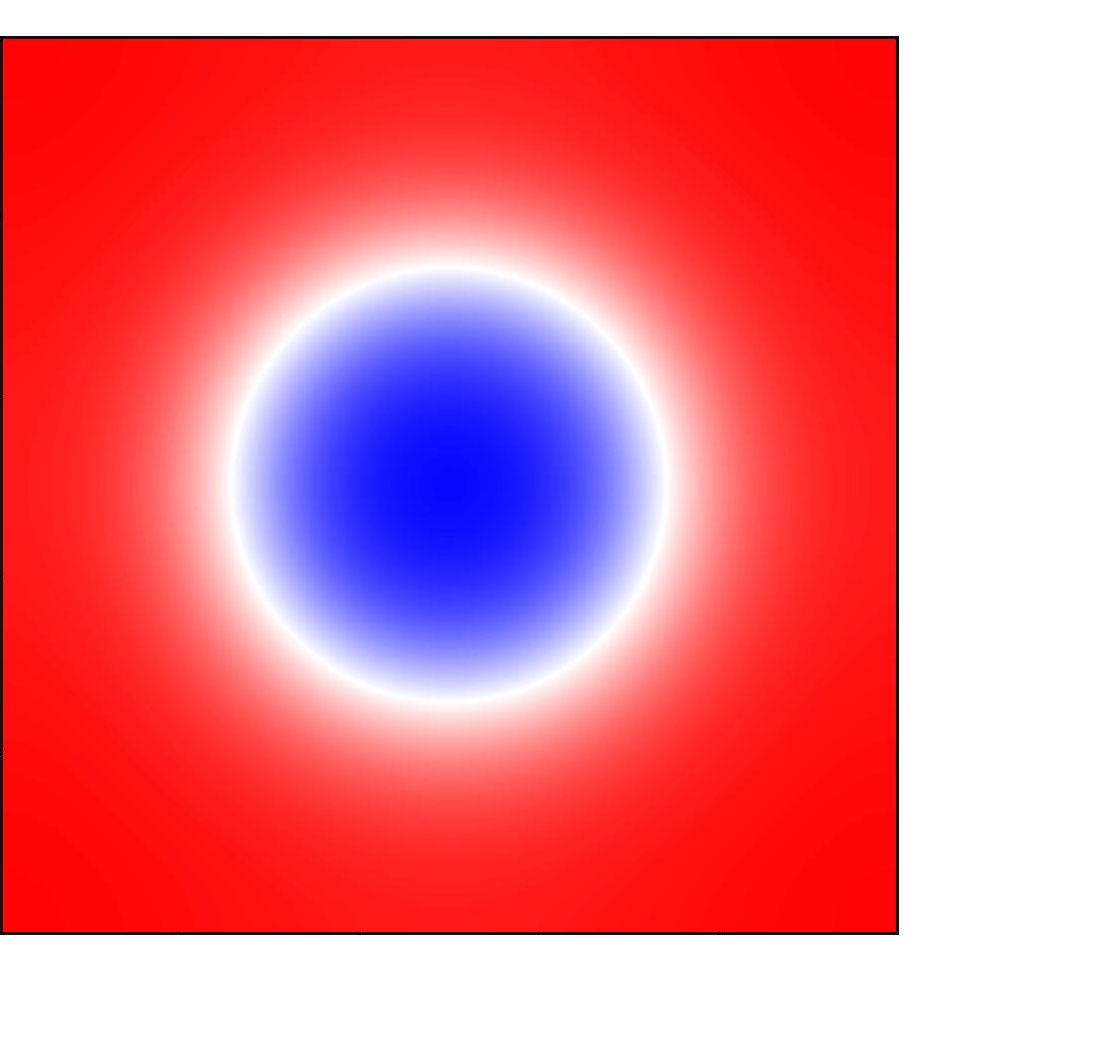}
  \end{subfigure}
  \begin{subfigure}[t]{0.32\textwidth}
        \centering
        \includegraphics[height=5cm]{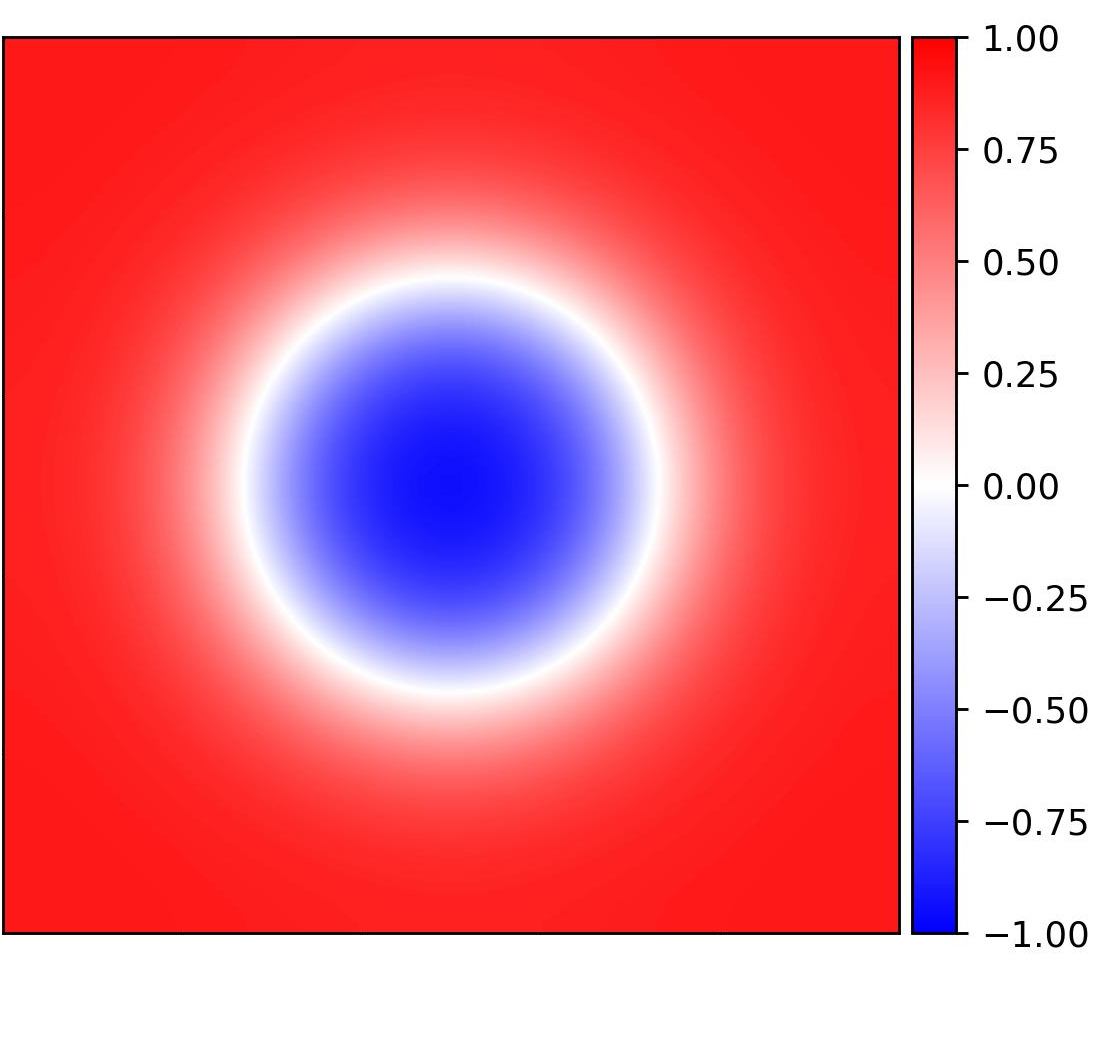}
  \end{subfigure}
  \caption{Solution of the phase field variable for the 2D bubble test case at $t = 0$ (left), $t = 0.01$ (center) and $t = 0.25$ (right) computed with $\Nx = \Ny = 50$.}
  \label{fig:Solution_Cosine2D}
\end{figure}

\begin{figure}[htpb]
  \centering
  \includegraphics[width=\linewidth]{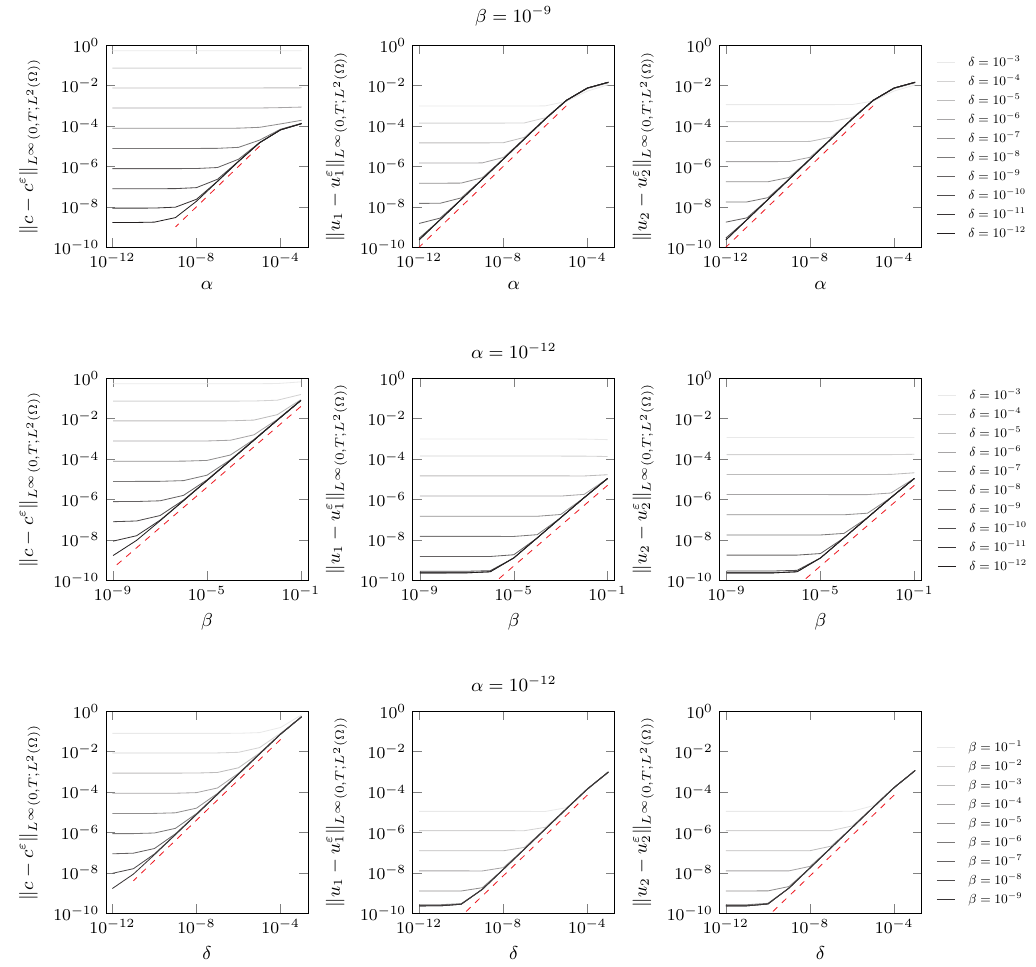}
  \caption{Convergence of the relaxation system \eqref{eq:NSCHr_num} to the NSCH system \eqref{eq:NSCH} for the single bubble test case $\Nx = \Ny = 25$. In the top row, the convergence in terms of $\alpha$ for a fixed value of $\beta$ and varying $\delta$ is depicted. The center row shows the convergence in terms of $\beta$ for a fixed $\alpha$ and varying $\delta$. In the bottom row, the convergence error in terms of $\delta$ for a fixed value of $\alpha$ and varying $\delta$ is plotted. The red dashed lines indicate convergence with order $\cO(\alpha)$, $\cO(\beta)$ and $\cO(\delta)$, depending on the respective running variable on the abscissa.}
  \label{fig:Error_Cosine2D_N25}
\end{figure}

As a first test case, the relaxation of an initially perturbed bubble to its equilibrium was simulated.
The solution of the phase field variable at $t = 0$, $t = 0.01$ and $t = 0.25$ computed with $\Nx = \Ny = 50$ is depicted in Fig.~\ref{fig:Solution_Cosine2D}.
For the visualization, a linear interpolation between the solution points of the phase field variable was used.
The relaxation to a spherical bubble over time is evident.
The error of the discretized equations \eqref{eq:NSCHr_semidiscrete} in comparison to the discretized equations \eqref{eq:NSCH_approx} for different choices of the relaxation parameters in $\numepsilon$ and the corresponding convergence behavior is depicted in Fig.~\ref{fig:Error_Cosine2D_N25}.
In the left column, the error in the phase field variable is plotted, while in the center and right columns, the velocities in $x$ and $y$ are shown, respectively.
Similar to the one-dimensional results:
In the top row, the convergence with respect to $\alpha$ for a fixed value of $\beta = 10^{-9}$ and decreasing values of $\delta \in [10^{-3},10^{-12}]$ is presented.
The center row highlights the convergence with respect to $\beta$ for a fixed value of $\alpha = 10^{-12}$ and decreasing values of $\delta \in [10^{-3},10^{-12}]$.
In the bottom row, the convergence with respect to $\delta$ for a fixed value of $\alpha = 10^{-12}$ and decreasing values of $\beta \in [10^{-1},10^{-9}]$ is demonstrated.
For all solution variables, the expected convergence of order $\cO (\alpha + \beta + \delta )$ is achieved,
Again, the red dashed lines indicate convergence in the order of the running variable on the abscissa for the remaining two relaxation parameters sufficiently small.
The stagnant error convergence observed in the top left as well as all plots in the center and bottom row can be attributed to a predominant error 
caused by the relaxation parameters $\alpha$ and $\beta$.
A further reduction of $\alpha$ and $\beta$ would resolve this problem.
Additional plots, investigating the convergence behavior also for the grid resolution $\Nx=\Ny=50$, are given in Appendix \ref{appx:2D} in Fig. \ref{fig:Error_Cosine2D_N50}.
Similar to the one-dimensional experiments, the order of convergence remains unaffected by the grid resolution.
%
%
%
%
\subsubsection{Merging droplets in 2D}
\begin{figure}[htpb]
  \centering
  \includegraphics[width=\linewidth]{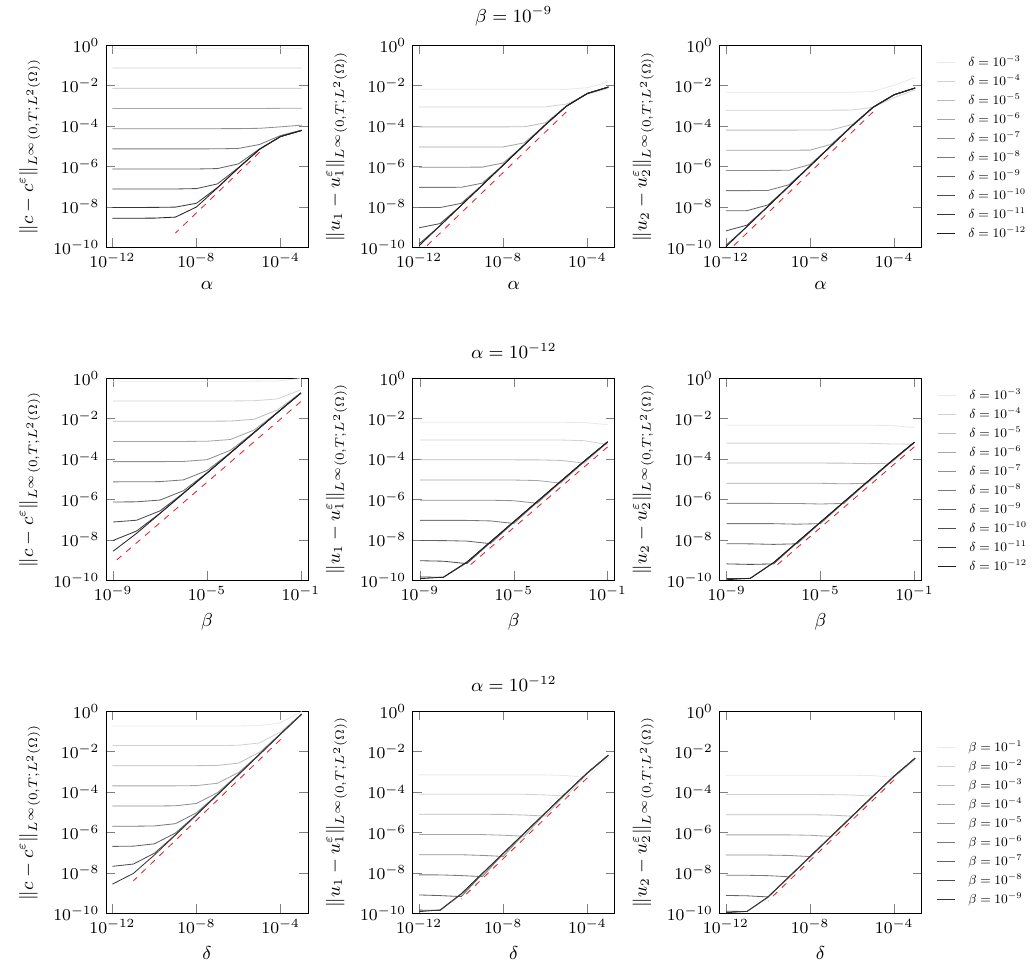}
  \caption{Convergence of the relaxation system \eqref{eq:NSCHr_num} to the NSCH system \eqref{eq:NSCH} for the merging droplets test case with $\Nx = \Ny = 25$. In the top row, the convergence in terms of $\alpha$ for a fixed value of $\beta$ and varying $\delta$ is depicted. The center row shows the convergence in terms of $\beta$ for a fixed $\alpha$ and varying $\delta$. In the bottom row, the convergence error in terms of $\delta$ for a fixed value of $\alpha$ and varying $\delta$ is plotted. The red dashed lines indicate convergence with order $\cO(\alpha)$, $\cO(\beta)$ and $\cO(\delta)$, depending on the respective running variable on the abscissa.}
  \label{fig:Error_MergingDroplets2D_N25}
\end{figure}

As a second numerical experiment, a slightly more dynamic setup, the merging of two droplets initially in contact, was simulated.
The solution of the phase field variable at $t = 0$, $t = 0.01$ and $t = 0.25$ computed with $\Nx = \Ny = 50$ is depicted in  Fig.~\ref{fig:Solution_Cosine2D}.
For the visualization, a linear interpolation between the solution points of the phase field variable was used.
\begin{figure}[htpb]
  \centering
  \begin{subfigure}[t]{0.32\textwidth}
        \centering
        \includegraphics[height=5cm]{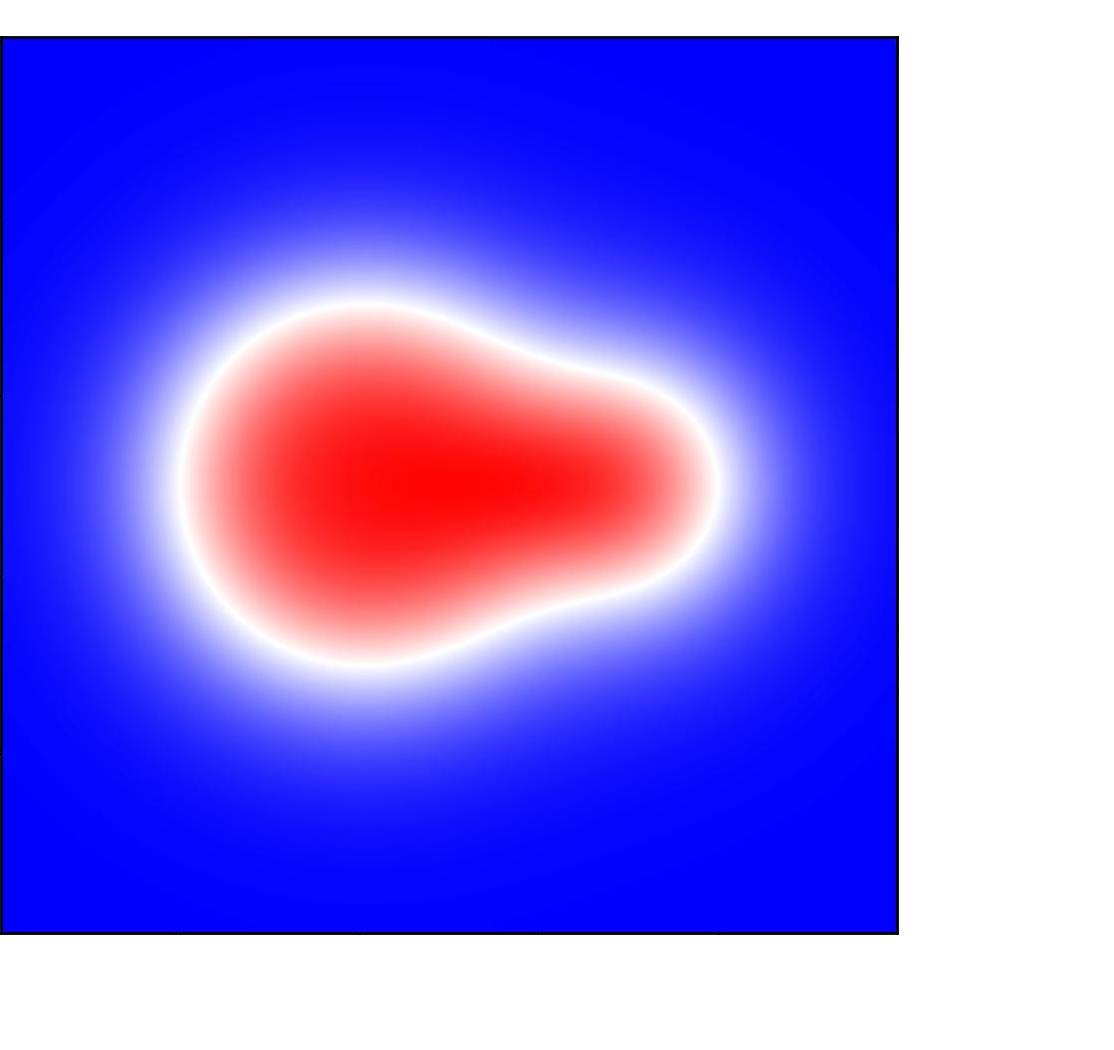}
  \end{subfigure}
  \begin{subfigure}[t]{0.32\textwidth}
        \centering
        \includegraphics[height=5cm]{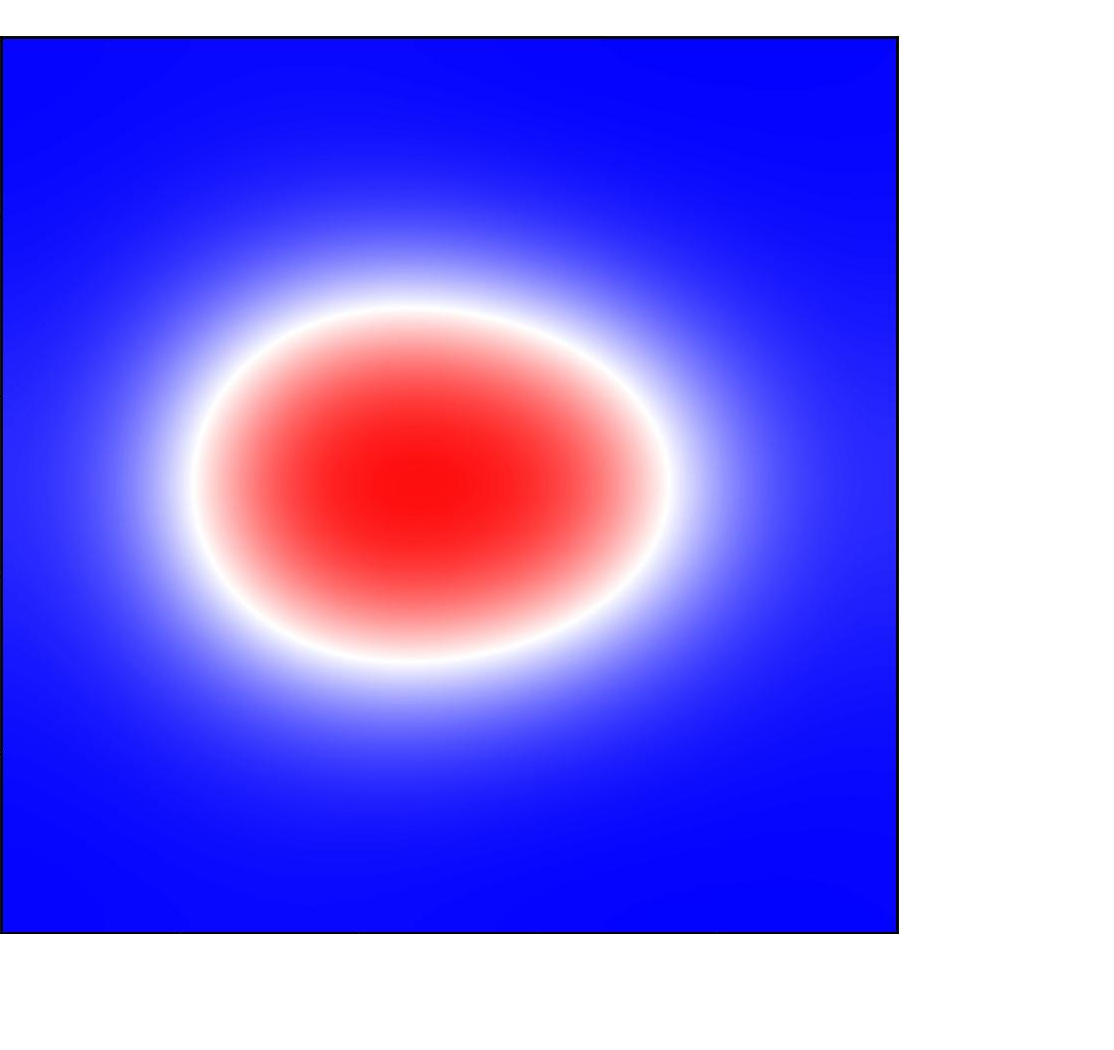}
  \end{subfigure}
  \begin{subfigure}[t]{0.32\textwidth}
        \centering
        \includegraphics[height=5cm]{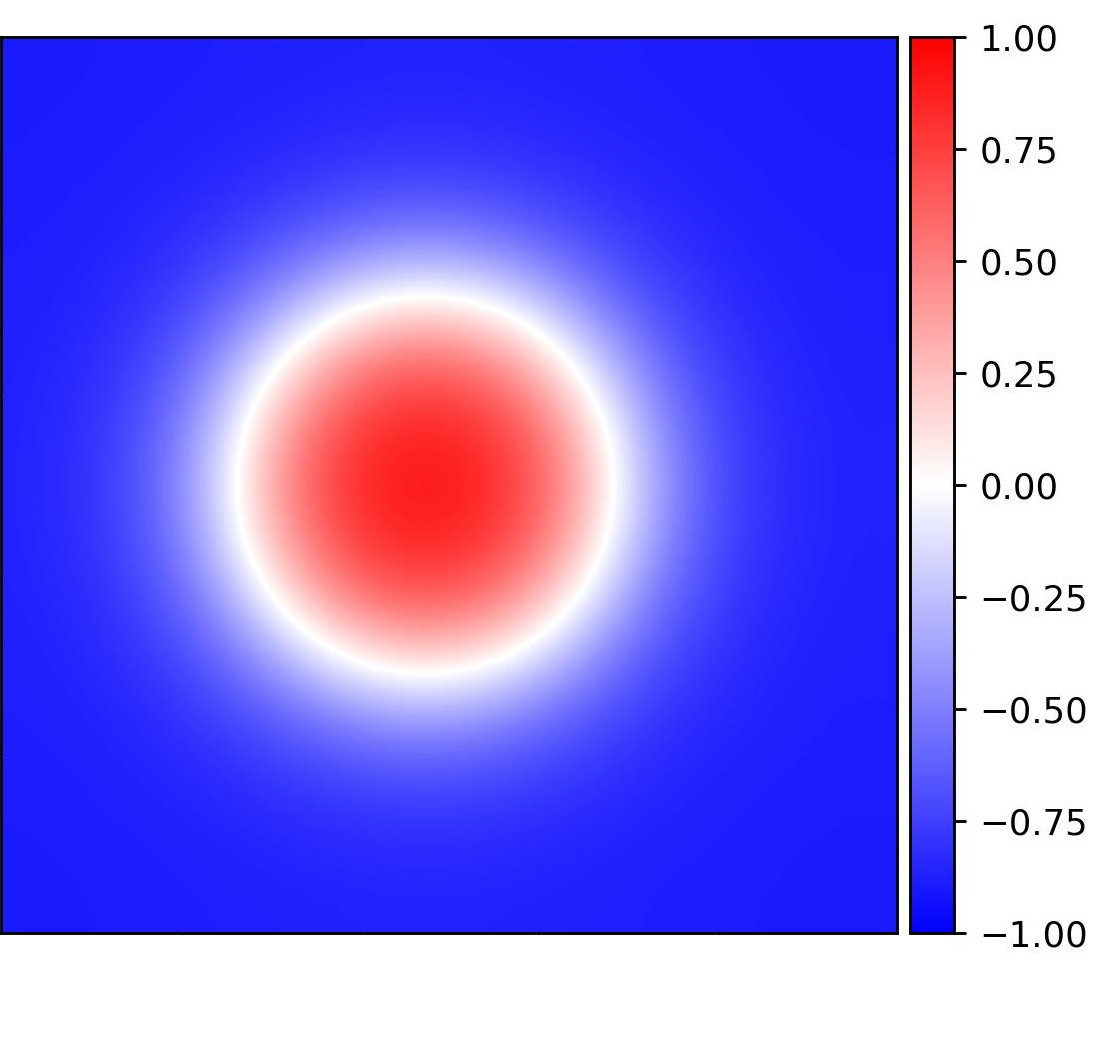}
  \end{subfigure}
  \caption{Solution of the phase field variable for the merging droplets test case at $t = 0$ (left), $t = 0.01$ (center) and $t = 0.25$ (right) computed with $\Nx = \Ny = 50$.}
  \label{fig:Solution_MergingDroplets2D}
\end{figure}
Over time, the droplets merge and form a single spherical droplet at the end.
The error of the discretized equations \eqref{eq:NSCHr_semidiscrete} in comparison to the discretized equations \eqref{eq:NSCH_approx} for different choices of the relaxation parameters in $\numepsilon$ and the corresponding convergence behavior is depicted in Fig.~\ref{fig:Error_MergingDroplets2D_N25}.
A similar arrangement of the plots to the previous section is used.
Similar to the previous test cases, the anticipated convergence of order $\cO (\alpha + \beta + \delta )$ is once more evident.
In comparison to the single bubble test case, the overall error level is slightly elevated due to the higher velocities in this numerical experiment.
Also in this setup, the parameters $\alpha$ and $\beta$ dominate the error for small values of $\delta$.
Additional plots investigating the convergence behavior also for the grid resolution $\Nx=\Ny=50$ are given in Appendix \ref{appx:2D} in Fig. \ref{fig:Error_MergingDroplets2D_N50}.
Similar to the previous observations, the order of convergence remains unaffected.
%
%
%
%
\subsubsection{2D droplet collision / droplet-wall interaction}
\begin{figure}[htpb]
  \centering
  \begin{subfigure}[t]{0.32\textwidth}
        \centering
        \includegraphics[height=5cm]{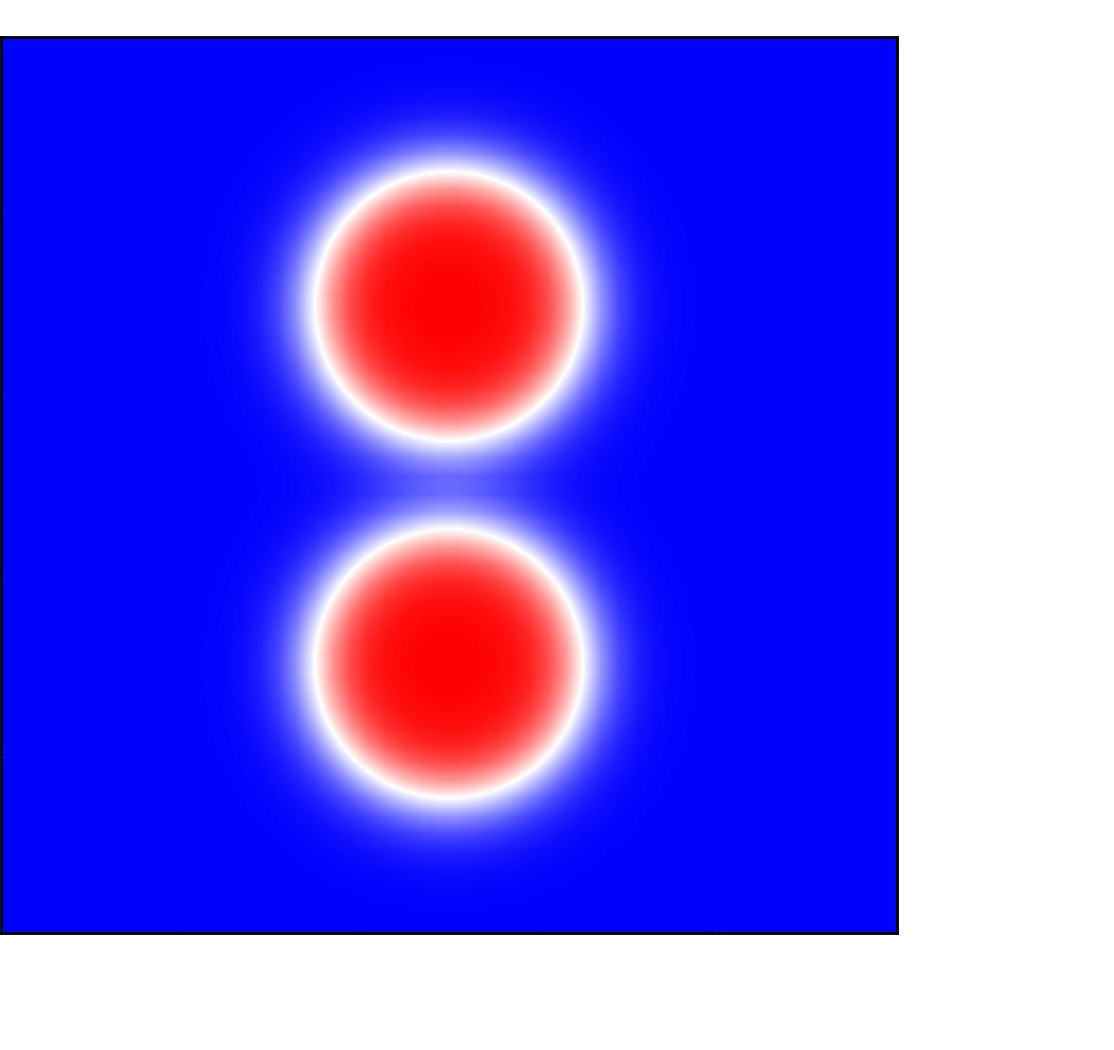}
  \end{subfigure}
  \begin{subfigure}[t]{0.32\textwidth}
        \centering
        \includegraphics[height=5cm]{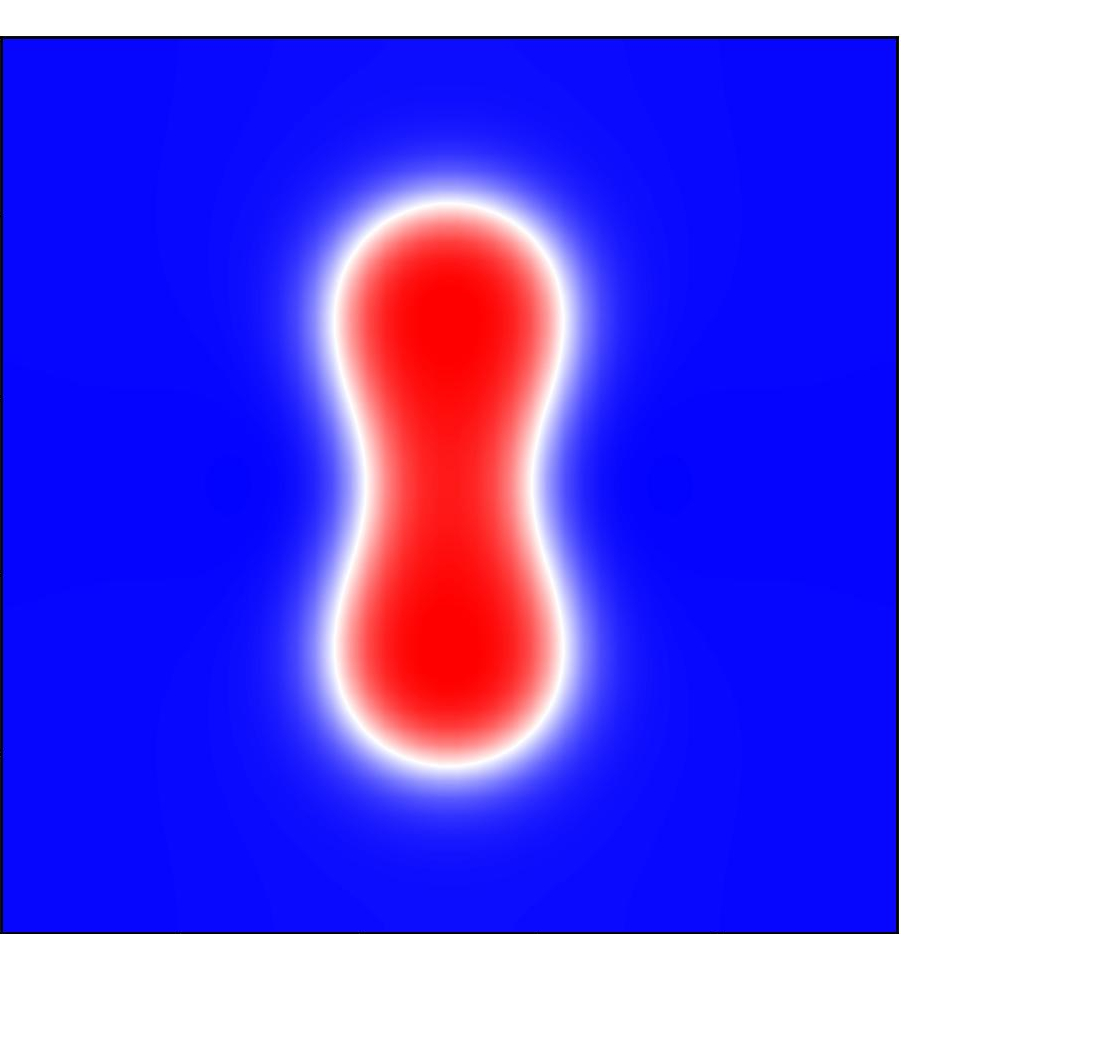}
  \end{subfigure}
  \begin{subfigure}[t]{0.32\textwidth}
        \centering
        \includegraphics[height=5cm]{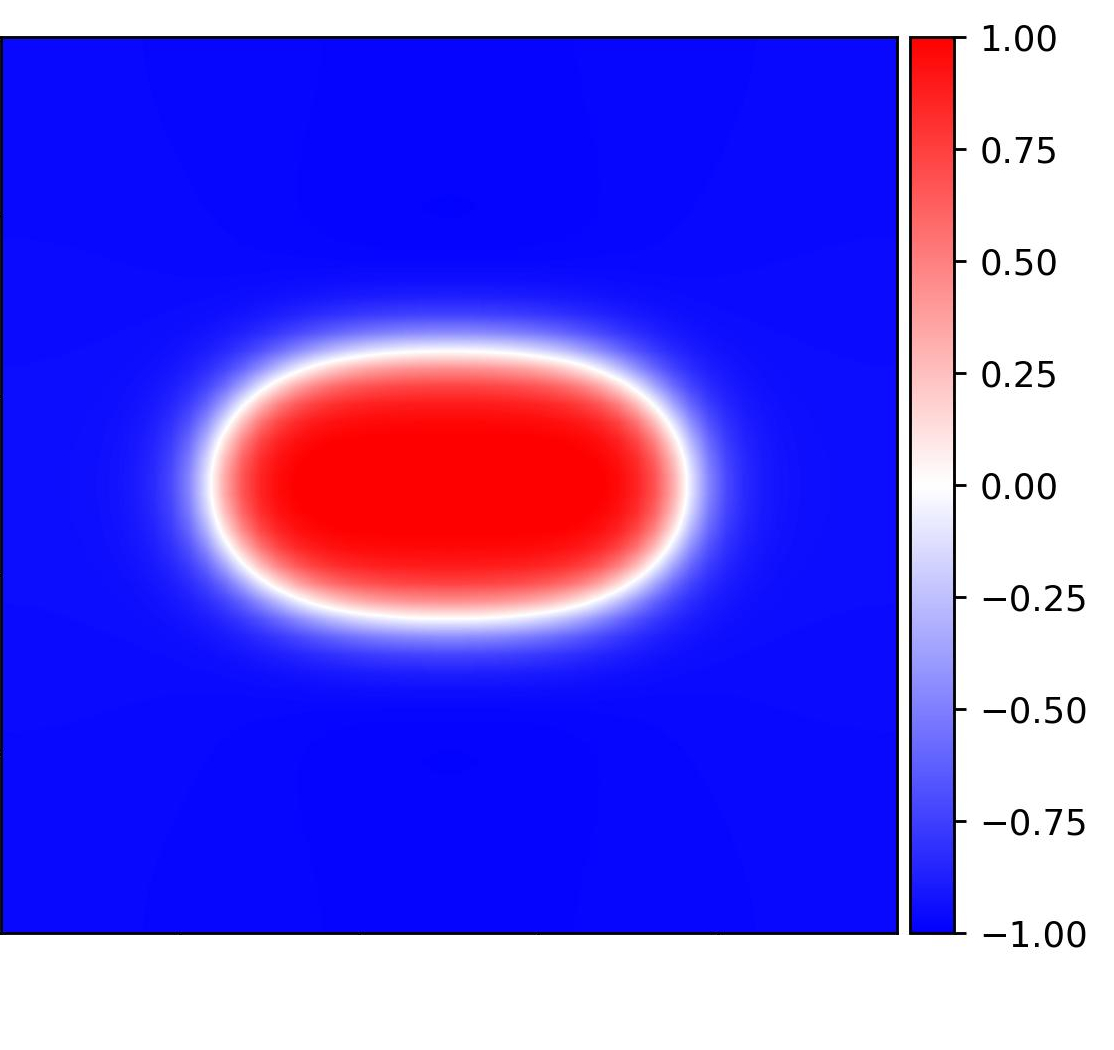}
  \end{subfigure}
  \caption{Solution of the phase field variable for the colliding droplets test case at $t = 0$ (left), $t = 0.03$ (center) and $t = 0.25$ (right) computed with $\Nx = \Ny = 50$.}
  \label{fig:Solution_CollidingDroplets2D}
\end{figure}

\begin{figure}[htpb]
  \centering
  \includegraphics[width=\linewidth]{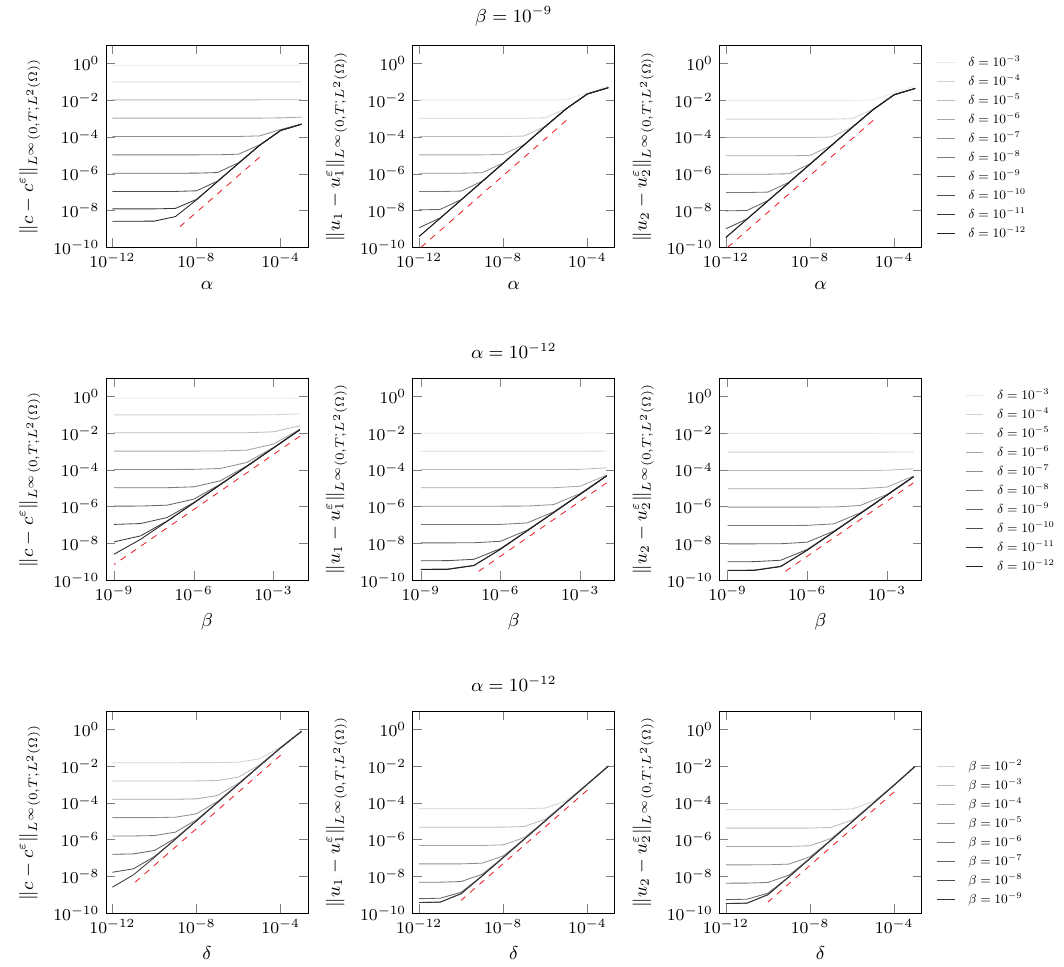}
  \caption{Convergence of the relaxation system \eqref{eq:NSCHr_num} to the NSCH system \eqref{eq:NSCH} for the colliding droplets test case with $\Nx = \Ny = 25$. In the top row, the convergence in terms of $\alpha$ for a fixed value of $\beta$ and varying $\delta$ is depicted. The center row shows the convergence in terms of $\beta$ for a fixed $\alpha$ and varying $\delta$. In the bottom row, the convergence error in terms of $\delta$ for a fixed value of $\alpha$ and varying $\delta$ is plotted. The red dashed lines indicate convergence with order $\cO(\alpha)$, $\cO(\beta)$ and $\cO(\delta)$, depending on the respective running variable on the abscissa.}
  \label{fig:Error_CollidingDroplets2D_N25}
\end{figure}

In the final test case, a highly dynamic droplet collision is investigated.
This setup is interesting for two reasons.
First, it exhibits considerably higher velocities in comparison to the preceding test cases.
Second, due to the symmetry of this numerical experiment, the test case is analogous to a droplet-wall interaction.
At the symmetry plane, the constraints of a droplet-wall interaction with a slip wall and a $90^{\degree}$ contact angle
\begin{align}
  \velv \cdot \mathbf{n} = 0, \quad \grad \phas \cdot \mathbf{n} = 0, \quad \grad \chempot \cdot \mathbf{n} = 0
\end{align}
are satisfied.
In contrast to the previous numerical experiments, a capillary parameter $\gamma = 10^{-3}$ was utilized to avoid an initial contact of the two droplets.
The solution of the phase field variable at $t = 0$, $t = 0.03$ and $t = 0.25$ computed with $\Nx = \Ny = 50$ is depicted in Fig.~\ref{fig:Solution_CollidingDroplets2D}.
Similar to the previous test cases, a linear interpolation between the solution points of the phase field variable was used for the visualization.
Over time, the droplets collide, merge and form the shape of an ellipsoid at $t = 0.25$.
If the simulation would be performed further, the remaining droplet would oscillate and relax to a spherical shape.
\newline
The error of the discretized equations \eqref{eq:NSCHr_semidiscrete} in comparison to the discretized equations \eqref{eq:NSCH_approx} for different choices of the relaxation parameters in $\numepsilon$ and the corresponding convergence behavior is displayed in Fig.~\ref{fig:Error_CollidingDroplets2D_N25}.
A similar arrangement of the plots to the previous section is used.
Also for this slightly modified setup, the same convergence behavior as for the previous test cases is observed, where the error decreases with $\cO ( \alpha + \beta + \delta )$.
Also in this setup, the parameters $\alpha$ and $\beta$ dominate the error for small values of $\delta$.
Additional plots investigating the convergence behavior also for the grid resolution $\Nx=\Ny=50$ are given in Appendix \ref{appx:2D} in Fig. \ref{fig:Error_CollidingDroplets2D_N50}.
Similar to the previous observations, the order of convergence remains unaffected.

\section{Conclusions}\label{chap:p2:conclusion}
We    investigated a relaxation approximation for the  classical NSCH system from \cite{hohenberg1977theory}   which provides a diffuse-interface approach for incompressible two-phase flow.  Based on the relative-entropy approach and under appropriate smoothness assumptions,
we  have shown that the  $L^2$-difference between 
the approximate solution and the solution of the limit NSCH system vanishes linearly with respect 
to the relaxation parameters. Note that the relaxation system also provides approximations of $\nabla c$ converging (again in the  $L^2$-norm) with  the same rate and of the pressure, the Laplacian of $c$ and the Cahn-Hilliard flux converging with lower rates. 
For a series of interesting two-phase settings we confirmed the 
theoretical convergence results by numerical simulations. The numerical results have been obtained 
by a new MAC-type finite-difference method  that is applicable in a seamless way for the limit 
NSCH system and the relaxation approximation. \\ 
As we have shown, our  numerical method from the realm of solvers for hyperbolic balance laws 
allows  high-velocity  settings and is applicable even to  inviscid flow.  In this sense it 
can act as a numerical method to solve the limit NSCH system directly in such regimes. 
Then, for practical computations, one needs  criteria to minimize the complete 
error in terms of the model error and the discretization error. This might be achieved by 
the development of appropriate a-posteriori error estimators \cite{GiesselmannPryer_2017}. Finally, for the  NSCH model from \cite{hohenberg1977theory}   various  improvements have been suggested  in the last decade. It remains to  investigate whether these models allow also for a hyperbolic relaxation approximation like our system \eqref{eq:NSCHr}  which is accessible  to relative-entropy analysis and  appropriate extensions of our numerical approach. Currently we work on specific extensions for problems with high-density contrasts and chemically reactive  multi-component mixtures, 
see e.g.~\cite{ABG,RvWDR26,RovWo}.

\paragraph{Acknowledgements}{Financial support by the German Research Foundation (DFG), within  the projects No. 525866748 and No. 526024901 of the Priority Programme - SPP 2410 Hyperbolic Balance Laws in Fluid Mechanics: Complexity, Scales, Randomness (CoScaRa) is acknowledged.}

\section*{Appendix}
\appendix
In the following, the second-order accurate conservative finite difference scheme utilized for the approximation of the spatial derivatives in \eqref{eq:NSCH_approx} and \eqref{eq:NSCHr_num} is detailed.
Although, this discretization is rather standard, for the sake of clarity and completeness, all employed stencils are defined.
Herein, indices refer to the definitions in Fig.~\ref{fig:MACgrid}.
\section{Spatial discretization of the NSCH system}\label{appx:NSCH_FD}
\subsection{Discretization of the evolution equations for the phase field variable \eqref{eq:NSCHphasstar_approx} and \eqref{eq:NSCHphas_approx}} \label{sec:NSCH_phas}
%
%
The $x$ contribution of the convective divergence operator in \eqref{eq:NSCHphasstar_approx} is approximated as
\begin{align}
  \left( \pfrac{\phas \velx}{x} \right)_{i,j} \approx \frac{ (\phas \velx)_{i+\oot,j} - (\phas \velx)_{i-\oot,j} }{\Dx} \quad \mathrm{with} \quad (\phas \velx)_{i+\oot,j} \approx \oot \Big( (\phas \velx)_{i+1,j} + (\phas \velx)_{i,j} \Big)
\end{align}
and
\begin{align}
  \left( \velx \right)_{i,j} \approx \frac{ \left(\velx\right)_{i+\oot,j} \Dx\Dy + \left(\velx\right)_{i-\oot,j} \Dx\Dy }{ 2 \Dx\Dy } = \oot \Big( \left(\velx\right)_{i+\oot,j} + \left(\velx\right)_{i-\oot,j} \Big).
\end{align}
Equivalently, the contribution in the $y$-direction reads as
\begin{align}
  \left( \pfrac{\phas \vely}{y} \right)_{i,j} \approx \frac{ (\phas \vely)_{i,j+\oot} - (\phas \vely)_{i,j-\oot} }{\Dy} \quad \mathrm{with} \quad (\phas \velx)_{i,j+\oot} \approx \oot \Big( (\phas \vely)_{i,j+1} + (\phas \vely)_{i,j} \Big)
\end{align}
and
\begin{align}
  \left( \vely \right)_{i,j} \approx \frac{ \left(\vely\right)_{i,j+\oot} \Dx\Dy + \left(\vely\right)_{i,j-\oot} \Dx\Dy }{ 2 \Dx\Dy } = \oot \Big( \left(\vely\right)_{i,j+\oot} + \left(\vely\right)_{i,j-\oot} \Big).
\end{align}
%
%
In the following, the spatial discretization of the terms related to the Cahn-Hilliard part in \eqref{eq:NSCHphasstar_approx} are detailed, starting with the contribution of the double-well potential
\begin{align}
  \left( \div \Big( \DoubWellPot^{\prime\prime} (\phas) \grad \phas \Big) \right)_{i,j}  = \left( \pfrac{\DoubWellPot^{\prime\prime} (\phas) \pfrac{\phas}{x}}{x} \right)_{i,j} + \left( \pfrac{\DoubWellPot^{\prime\prime} (\phas) \pfrac{\phas}{y}}{y} \right)_{i,j}.
\end{align}
Here, the first term ($x$-direction) is approximated as
\begin{align}
  \left( \pfrac{\DoubWellPot^{\prime\prime} (\phas) \pfrac{\phas}{x}}{x} \right)_{i,j} \approx \frac{ \left( \DoubWellPot^{\prime\prime} (\phas) \pfrac{\phas}{x} \right)_{i+\oot,j} - \left( \DoubWellPot^{\prime\prime} (\phas) \pfrac{\phas}{x} \right)_{i-\oot,j} }{\Dx}
\end{align}
with
\begin{align}
  \label{eq:DoubWellPot_ipootj}
  \Big( \DoubWellPot^{\prime\prime} (\phas) \Big)_{i+\oot,j} \approx \frac{ - \DoubWellPot^{\prime\prime} (\phas_{i-1,j}) + 7 \DoubWellPot^{\prime\prime} (\phas_{i,j}) + 7 \DoubWellPot^{\prime\prime} (\phas_{i+1,j}) - \DoubWellPot^{\prime\prime} (\phas_{i+2,j}) }{12}
\end{align}
and
\begin{align}
  \label{eq:gradxc_ipootj}
  \left( \pfrac{\phas}{x} \right)_{i+\oot,j} \approx \frac{ \phas_{i-1,j} - 7 \phas_{i,j} + 7 \phas_{i+1,j} - \phas_{i+2,j} }{12 \Dx}.
\end{align}
Similarly, the second term ($y$-direction) reads as
\begin{align}
  \left( \pfrac{\DoubWellPot^{\prime\prime} (\phas) \pfrac{\phas}{y}}{y} \right)_{i,j} \approx \frac{ \left( \DoubWellPot^{\prime\prime} (\phas) \pfrac{\phas}{x} \right)_{i,j+\oot} - \left( \DoubWellPot^{\prime\prime} (\phas) \pfrac{\phas}{y} \right)_{i,j-\oot} }{\Dy}
\end{align}
with
\begin{align}
  \label{eq:DoubWellPot_ijpoot}
  \Big( \DoubWellPot^{\prime\prime} (\phas) \Big)_{i,j+\oot} \approx \frac{ - \DoubWellPot^{\prime\prime} (\phas_{i,j-1}) + 7 \DoubWellPot^{\prime\prime} (\phas_{i,j}) + 7 \DoubWellPot^{\prime\prime} (\phas_{i,j+1}) - \DoubWellPot^{\prime\prime} (\phas_{i,j+2}) }{12}
\end{align}
and
\begin{align}
  \label{eq:gradyc_ijpoot}
  \left( \pfrac{\phas}{y} \right)_{i,j+\oot} \approx \frac{ \phas_{i,j-1} - 7 \phas_{i,j} + 7 \phas_{i,j+1} - \phas_{i,j+2} }{12 \Dy}.
\end{align}
The remaining fourth-order contribution
\begin{align}
  \Big( \div (\gamma \grad \lap \phas ) \Big)_{i,j} = \pfrac{ }{x} \left( \pppfrac{\phas}{x} + \pppxyyfrac{\phas}{x}{y} \right)_{i,j} + \pfrac{ }{y} \left( \pppxxyfrac{\phas}{x}{y} + \pppfrac{\phas}{y} \right)_{i,j}
\end{align}
is approximated as
\begin{align}
  \label{eq:biharmonic_approx}
  \Big( \div (\gamma \grad \lap \phas ) \Big)_{i,j} \approx &{\color{white}{+}} \frac{ \left( \pppfrac{\phas}{x} + \pppxyyfrac{\phas}{x}{y} \right)_{i+\oot,j} - \left( \pppfrac{\phas}{x} + \pppxyyfrac{\phas}{x}{y} \right)_{i-\oot,j} }{\Dx} \nonumber \\
                                                            &+ \frac{ \left( \pppxxyfrac{\phas}{x}{y} + \pppfrac{\phas}{y} \right)_{i,j+\oot} - \left( \pppxxyfrac{\phas}{x}{y} + \pppfrac{\phas}{y} \right)_{i,j-\oot} }{\Dy}
\end{align}
with
\begin{align}
  \label{eq:gradxlapc}
  \left( \pppfrac{\phas}{x} + \pppxyyfrac{\phas}{x}{y} \right)_{i+\oot,j} \approx &{\color{white}{+}} \frac{ -\phas_{i-1,j} + 3 \phas_{i,j} - 3 \phas_{i+1,j} + \phas_{i+2,j} }{\Dx^3} \nonumber \\
                                                                                  &+ \frac{ -\phas_{i,j-1} + \phas_{i+1,j-1} + 2 \phas_{i,j} - 2 \phas_{i+1,j} - \phas_{i,j+1} + \phas_{i+1,j+1} }{\Dx^2\Dy}
\end{align}
and
\begin{align}
  \label{eq:gradylapc}
  \left( \pppxxyfrac{\phas}{x}{y} + \pppfrac{\phas}{y} \right)_{i,j+\oot} \approx &{\color{white}{+}} \frac{ -\phas_{i,j-1} + 3 \phas_{i,j} - 3 \phas_{i,j+1} + \phas_{i,j+2} }{\Dy^3} \nonumber \\
                                                                                  &+ \frac{ -\phas_{i-1,j} + \phas_{i-1,j+1} + 2 \phas_{i,j} - 2 \phas_{i,j+1} - \phas_{i+1,j} + \phas_{i+1,j+1} }{\Dx\Dy^2}
\end{align}
Straightforward computations show that combining \eqref{eq:biharmonic_approx} to \eqref{eq:gradylapc} yields a second-order accurate central approximation of the biharmonic operator.
This stencil choice for the Cahn-Hilliard part in \eqref{eq:NSCHphasstar_approx} has previously been used in \cite{Dhaouadi24} and forms an overall second-order accurate approximation.
\subsection{Discretization of the momentum equations \eqref{eq:NSCHvelstar_approx}} \label{sec:NSCH_velstar}
%
%
The approximation of the convective term of the momentum equation in the $x$-direction
\begin{align}
  \Big( \div ( \velx \velv ) \Big)_{i+\oot,j} = \pfrac{}{x} \Big( \velx \velx \Big)_{i+\oot,j} + \pfrac{}{y} \Big( \velx \vely \Big)_{i+\oot,j}
\end{align}
reads as
\begin{align}
  \Big( \div ( \velx \velv ) \Big)_{i+\oot,j} \approx \frac{ \left( \velx \velx \right)_{i+1,j} - \left( \velx \velx \right)_{i,j} }{\Dx} + \frac{ \left( \velx \vely \right)_{i+\oot,j+\oot} - \left( \velx \vely \right)_{i+\oot,j-\oot} }{\Dy}
\end{align}
with
\begin{align}
  \left( \velx \velx \right)_{i,j} \approx \frac{1}{2} \left( (\velx \velx)_{i+\oot,j} + (\velx \velx)_{i-\oot,j} \right)
\end{align}
and
\begin{align}
  \left( \velx \vely \right)_{i+\oot,j+\oot} \approx \frac{1}{2} \left( (\velx \vely)_{i+\oot,j+1} + (\velx \vely)_{i+\oot,j} \right),
\end{align}
where
\begin{align}
  \left( \vely \right)_{i+\oot,j} \approx \frac{ (\vely)_{i,j-\oot} + (\vely)_{i+1,j-\oot} + (\vely)_{i,j+\oot} + (\vely)_{i+1,j+\oot} }{4}.
\end{align}
%
%
Similarly to the $x$-direction, the approximation of the convective term of the momentum equation in the $y$-direction
\begin{align}
  \Big( \div ( \vely \velv ) \Big)_{i,j+\oot} = \pfrac{}{x} \Big( \vely \velx \Big)_{i,j+\oot} + \pfrac{}{y} \Big( \vely \vely \Big)_{i,j+\oot}
\end{align}
reads as
\begin{align}
  \Big( \div ( \vely \velv ) \Big)_{i,j+\oot} \approx \frac{ \left( \vely \velx \right)_{i+\oot,j+\oot} - \left( \vely \velx \right)_{i-\oot,j+\oot} }{\Dx} + \frac{ \left( \vely \vely \right)_{i,j+1} - \left( \vely \vely \right)_{i,j} }{\Dy}
\end{align}
with
\begin{align}
  \left( \vely \vely \right)_{i,j} \approx \frac{1}{2} \left( (\vely \vely)_{i,j+\oot} + (\velx \velx)_{i,j-\oot} \right)
\end{align}
and
\begin{align}
  \left( \vely \velx \right)_{i+\oot,j+\oot} \approx \frac{1}{2} \left( (\velx \vely)_{i+1,j+\oot} + (\velx \vely)_{i,j+\oot} \right),
\end{align}
where
\begin{align}
  \left( \velx \right)_{i,j+\oot} \approx \frac{ (\velx)_{i-\oot,j} + (\velx)_{i+\oot,j} + (\velx)_{i-\oot,j+1} + (\vely)_{i+\oot,j+1} }{4}.
\end{align}
%
%
For the approximation of the Cahn-Hilliard contributions
\begin{align}
  \left( \phas \left[ \DoubWellPot^{\prime\prime} (\phas) \pfrac{\phas}{x} - \gamma \left( \pppfrac{c}{x} + \pppxyyfrac{c}{x}{y} \right) \right] \right)_{i+\oot,j} \quad \mathrm{and} \quad \left( \phas \left[ \DoubWellPot^{\prime\prime} (\phas) \pfrac{\phas}{y} - \gamma \left( \pppxxyfrac{c}{x}{y} + \pppfrac{c}{y} \right) \right] \right)_{i,j+\oot}
\end{align}
in the momentum equations in the $x$- and $y$-directions, respectively, except from $(\phas)_{i+\oot,j}$ and $(\phas)_{i,j+\oot}$, all definitions were already introduced in \eqref{eq:DoubWellPot_ipootj}, \eqref{eq:gradxc_ipootj}, \eqref{eq:DoubWellPot_ijpoot}, \eqref{eq:gradyc_ijpoot}, \eqref{eq:gradxlapc} and \eqref{eq:gradylapc}.
The missing two contributions for $(\phas)_{i+\oot,j}$ and $(\phas)_{i,j+\oot}$ are approximated as
\begin{align}
  (\phas)_{i+\oot,j} \approx \frac{ - \phas_{i-1,j} + 7 \phas_{i,j} + 7 \phas_{i+1,j} - \phas_{i+2,j} }{12}, \quad
  (\phas)_{i,j+\oot} \approx \frac{ - \phas_{i,j-1} + 7 \phas_{i,j} + 7 \phas_{i,j+1} - \phas_{i,j+2} }{12}.
\end{align}
\subsection{Discretization of the pressure equation \eqref{eq:NSCHpres_approx}} \label{sec:NSCH_pres}
%
%
For the approximation of the Laplace and divergence operators in the pressure equation \eqref{eq:NSCHpres_approx}, standard second-order central stencils are utilized, given as
\begin{align}
  (\lap \pres)_{i,j} \approx \frac{ \pres_{i-1,j} - 2 \pres_{i,j} + \pres_{i+1,j} }{\Dx^2} + \frac{ \pres_{i,j-1} - 2 \pres_{i,j} + \pres_{i,j+1} }{\Dy^2}
\end{align}
and
\begin{align}
  (\div \velv)_{i,j} \approx \frac{ (\velx)_{i+\oot,j} - (\velx)_{i-\oot,j} }{\Dx} + \frac{ (\vely)_{i,j+\oot} - (\vely)_{i,j-\oot} }{\Dy}.
\end{align}
\subsection{Discretization of the momentum equations \eqref{eq:NSCHvel_approx}} \label{sec:NSCH_vel}
%
%
In the correction step of the momentum equations \eqref{eq:NSCHvel_approx}, for the gradients of the pressure in the $x$- and $y$-direction, a standard central second-order stencil is applied, reading as
\begin{align}
  \left( \pfrac{\pres}{x} \right)_{i+\oot,j} \approx \frac{ \pres_{i+1,j} - \pres_{i,j} }{\Dx} \quad \mathrm{and} \quad \left( \pfrac{\pres}{y} \right)_{i+\oot,j} \approx \frac{ \pres_{i,j+1} - \pres_{i,j} }{\Dy},
\end{align}
respectively.

\section{Spatial discretization of the relaxation approximation  of the NSCH system}\label{appx:NSCHr_FD}
In this appendix, the spatial discretization of the relaxation formulation of the NSCH equations \eqref{eq:NSCHr_num} is detailed.
For this, first, the utilized update formulas and their sequence is summarized in the following
\begin{subequations}
  \tiny
  \label{eq:NSCHr_semidiscrete}
  \begin{align}
    \label{eq:NSCHr_phasstar}
    \phas^{\star} + \Dt~\left[ \div (\phas^{\star} \velv^{n}) - \frac{\Dt}{\delta + \Dt}~\div \Big[ \DoubWellPot^{\prime \prime} (\phas^{n}) \grad \phas^{\star} - \grad \lap \OperatorORDP^{-1} \phas^{\star} \Big] \right] =& ~\phas^{n} - \frac{\delta \Dt}{\delta + \Dt} \div \vecmu^{n}, \\
    \label{eq:NSCHr_ordpstar}
    {\color{white}{\Big[}} \ordp^{\star} =& ~ \OperatorORDP^{-1} \phas^{\star} {\color{white}{\Big]}}, \\
    \label{eq:NSCHr_velstar}
  \velv^{\star} + \Dt ~ \div \left( \velv^{\star} \otimes \velv^{n} \right) =& ~\velv^{n} - \Dt ~ c^{\star} \Big[ \DoubWellPot^{\prime\prime} (c^{n}) \grad c^{\star} - \gamma \grad \lap \ordp^{\star} \Big], \\
    \label{eq:NSCHr_pres}
    \alpha \pres^{n+1} - \Dt^2 \lap \pres^{n+1} =& ~\alpha \pres^{n} -  \Dt~\div \velv^{\star}, \\
    \label{eq:NSCHr_velstarstar}
    \velv^{\star \star} =& ~\velv^{\star} - \Dt~\grad \pres^{n+1}, \\
    \label{eq:NSCHr_vel}
    \velv^{n+1} =& ~\velv^{\star\star} - \Dt~\frac{1}{2} (\div \velv^{\star\star}) \velv^{\star\star}, \\
    \label{eq:NSCHr_phas}
    \phas^{n+1} + \Dt~\left[ \div (\phas \velv)^{n+1} - \frac{\Dt}{\delta + \Dt}~\div \Big[ \DoubWellPot^{\prime \prime} (\phas^{\star}) \grad \phas^{n+1} - \grad \lap \OperatorORDP^{-1} \phas^{n+1} \Big] \right] =& ~\phas^{n} - \frac{\delta \Dt}{\delta + \Dt} \div \vecmu^{n}, \\
    \label{eq:NSCHr_ordp}
    {\color{white}{\Big[}} \ordp^{n+1} =& ~ \OperatorORDP^{-1} \phas^{n+1} {\color{white}{\Big]}}, \\
    \label{eq:NSCHr_mu}
    \vecmu^{n+1} =& ~\frac{\delta}{\delta + \Dt} \vecmu^{n} - \frac{\Dt}{\delta + \Dt}~\Big[ \DoubWellPot^{\prime \prime} (\phas^{\star}) \grad \phas^{n+1} + \frac{1}{\beta} \grad \lap \ordp^{n+1} \Big]
  \end{align}
\end{subequations}
with
\begin{align}
  \label{eq:ordp_elliptic}
  \OperatorORDP \ordp =  ( \mathcal{I} - \gamma \beta \lap ) \ordp = \phas.
\end{align}
Starting with the discretization of the Laplace operator in \eqref{eq:ordp_elliptic}, a standard second-order accurate central stencil, reading as
\begin{align}
    \label{eq:omega_FD}
  \left( \lap \ordp \right)_{i,j} = \frac{ \ordp_{i-1,j} - 2 \ordp_{i,j} + \ordp_{i+1,j} }{\Dx^2} + \frac{ \ordp_{i,j-1} - 2 \ordp_{i,j} + \ordp_{i,j+1} }{\Dy^2},
\end{align}
has been used.
The corresponding inverse of the operator $\OperatorORDP$ utilized in \eqref{eq:NSCHr_phasstar}, \eqref{eq:NSCHr_ordpstar}, \eqref{eq:NSCHr_phas}, \eqref{eq:NSCHr_ordp} and \eqref{eq:NSCHr_mu} is computed by first performing an LU-factorization and second, applying the Jacobi-Gauss algorithm.
\subsection{Discretization of the evolution equations for the phase field variable \eqref{eq:NSCHr_phasstar} and \eqref{eq:NSCHr_phas}} \label{sec:NSCHr_phas}
For the discretization of \eqref{eq:NSCHr_phasstar} and \eqref{eq:NSCHr_phas}, the employed stencils have already been introduced in \cref{sec:NSCH_phas}.
The only missing contribution is the divergence operator $\div \vecmu$ which, analogously to the velocity, is approximated by a central second-order accurate finite difference, given as
\begin{align}
  ( \div \vecmu)_{i,j} = \frac{ (\vecmux)_{i+\oot,j} - (\vecmux)_{i-\oot,j} }{\Dx} + \frac{ (\vecmuy)_{i,j+\oot} - (\vecmuy)_{i,j-\oot} }{\Dy}.
\end{align}
\subsection{Discretization of the pressure and momentum equations \eqref{eq:NSCHr_velstar} to \eqref{eq:NSCHr_vel}}
Similar to the previous section, the majority of the used stencils has been introduced in \cref{sec:NSCH_velstar,sec:NSCH_pres,sec:NSCH_vel}.
Herein, for the approximation of the third-order contribution $\grad \lap \ordp$, an analogous stencil as in \eqref{eq:gradxlapc} and \eqref{eq:gradylapc} is utilized only in terms of $\ordp$ instead of $\phas$.
The only missing contribution is the non-conservative product in \eqref{eq:NSCHr_vel}, where the stencils
\begin{align}
  {\color{white}{\Bigg[}} \left( (\div \velv) \velx \right)_{i+\oot,j} =& ~\frac{ (\div \velv)_{i+1,j} + (\div \velv)_{i,j} }{2} ~ (\velx)_{i+\oot,j}, {\color{white}{\Bigg]}} \\
  {\color{white}{\Bigg[}} \left( (\div \velv) \vely \right)_{i,j+\oot} =& ~\frac{ (\div \velv)_{i,j+1} + (\div \velv)_{i,j} }{2} ~ (\vely)_{i,j+\oot} {\color{white}{\Bigg]}}
\end{align}
are employed.
\subsection{Discretization of the evolution equation for the flux variable \eqref{eq:NSCHr_mu}}
All stencils utilized for the discretization of \eqref{eq:NSCHr_mu} have previously been introduced in section \ref{sec:NSCH_phas} and \eqref{eq:omega_FD}.

\section{1D Experiments}\label{appx:1D}
\begin{figure}[htpb]
  \centering
  \includegraphics[width=0.8\linewidth]{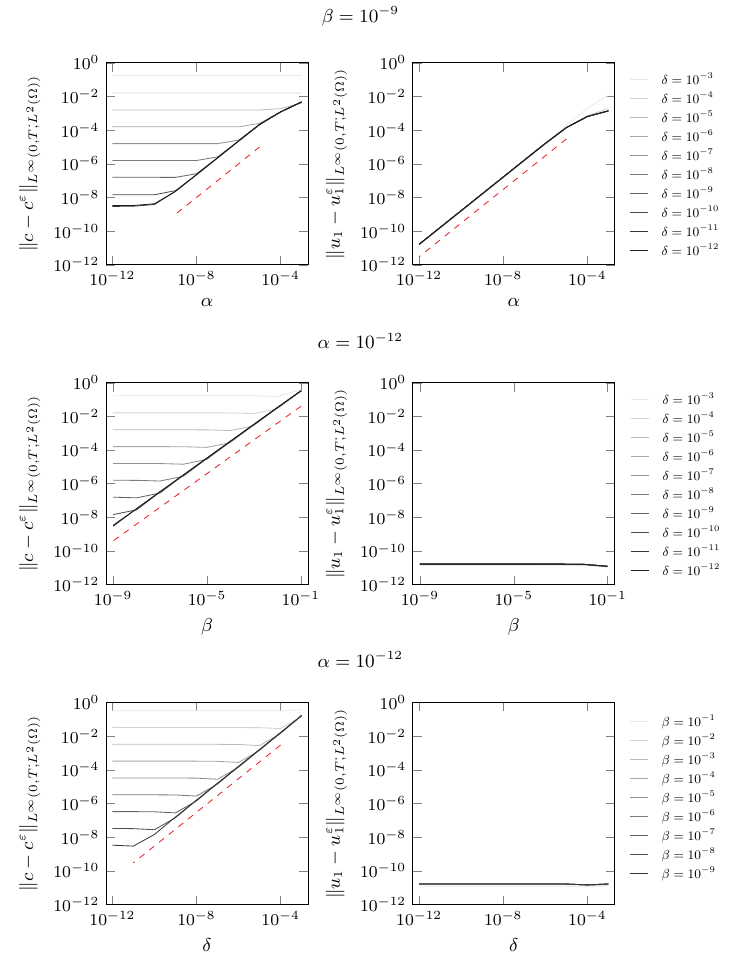}
  \caption{Convergence of the relaxation system \eqref{eq:NSCHr_num} to the NSCH system \eqref{eq:NSCH} for the Ostwald ripening test case with $\Nx = 500$. In the top row, the convergence in terms of $\alpha$ for a fixed value of $\beta$ and varying $\delta$ is depicted. The center row shows the convergence in terms of $\beta$ for a fixed $\alpha$ and varying $\delta$. In the bottom row, the convergence error in terms of $\delta$ for a fixed value of $\alpha$ and varying $\delta$ is plotted. The red dashed lines indicate convergence with order $\cO(\alpha)$, $\cO(\beta)$ and $\cO(\delta)$, depending on the respective running variable on the abscissa.}
  \label{fig:Error_Ostwald1D_N500}
\end{figure}

\newpage
\section{2D Experiments}\label{appx:2D}
\begin{figure}[htpb]
  \centering
  \includegraphics[width=\linewidth]{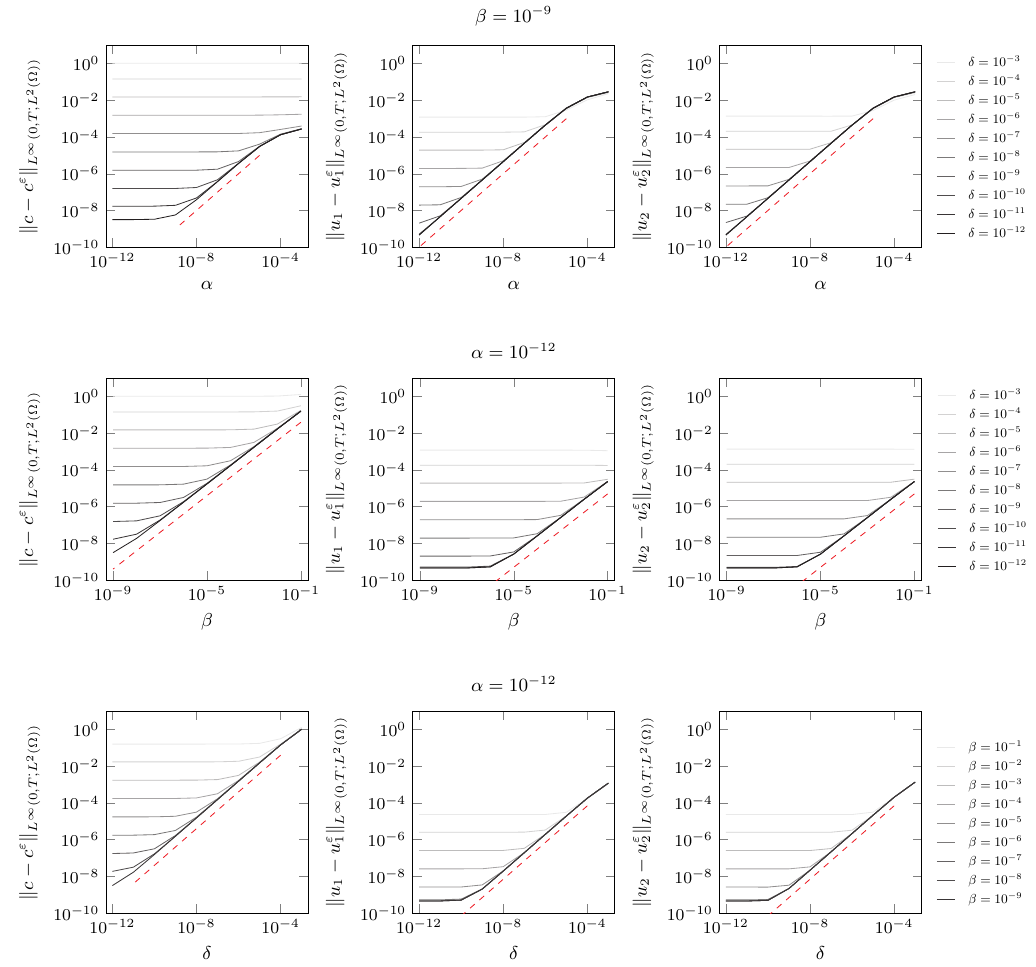}
  \caption{Convergence of the relaxation system \eqref{eq:NSCHr_num} to the NSCH system \eqref{eq:NSCH} for the single bubble test case $\Nx = \Ny = 50$. In the top row, the convergence in terms of $\alpha$ for a fixed value of $\beta$ and varying $\delta$ is depicted. The center row shows the convergence in terms of $\beta$ for a fixed $\alpha$ and varying $\delta$. In the bottom row, the convergence error in terms of $\delta$ for a fixed value of $\alpha$ and varying $\delta$ is plotted. The red dashed lines indicate convergence with order $\cO(\alpha)$, $\cO(\beta)$ and $\cO(\delta)$, depending on the respective running variable on the abscissa.}
  \label{fig:Error_Cosine2D_N50}
\end{figure}
\begin{figure}[htpb]
  \centering
  \includegraphics[width=\linewidth]{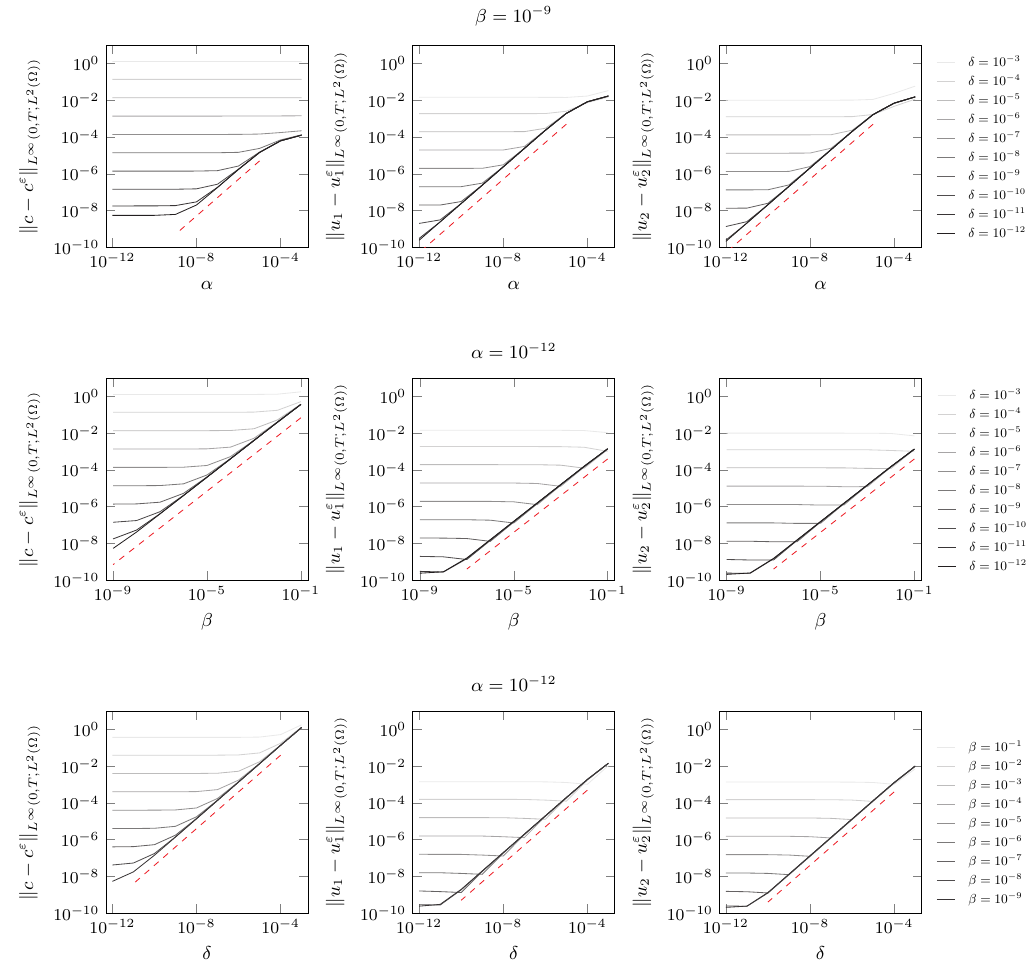}
  \caption{Convergence of the relaxation system \eqref{eq:NSCHr_num} to the NSCH system \eqref{eq:NSCH} for the merging droplets test case with $\Nx = \Ny = 50$. In the top row, the convergence in terms of $\alpha$ for a fixed value of $\beta$ and varying $\delta$ is depicted. The center row shows the convergence in terms of $\beta$ for a fixed $\alpha$ and varying $\delta$. In the bottom row, the convergence error in terms of $\delta$ for a fixed value of $\alpha$ and varying $\delta$ is plotted. The red dashed lines indicate convergence with order $\cO(\alpha)$, $\cO(\beta)$ and $\cO(\delta)$, depending on the respective running variable on the abscissa.}
  \label{fig:Error_MergingDroplets2D_N50}
\end{figure}
\begin{figure}[htpb]
  \centering
  \includegraphics[width=\linewidth]{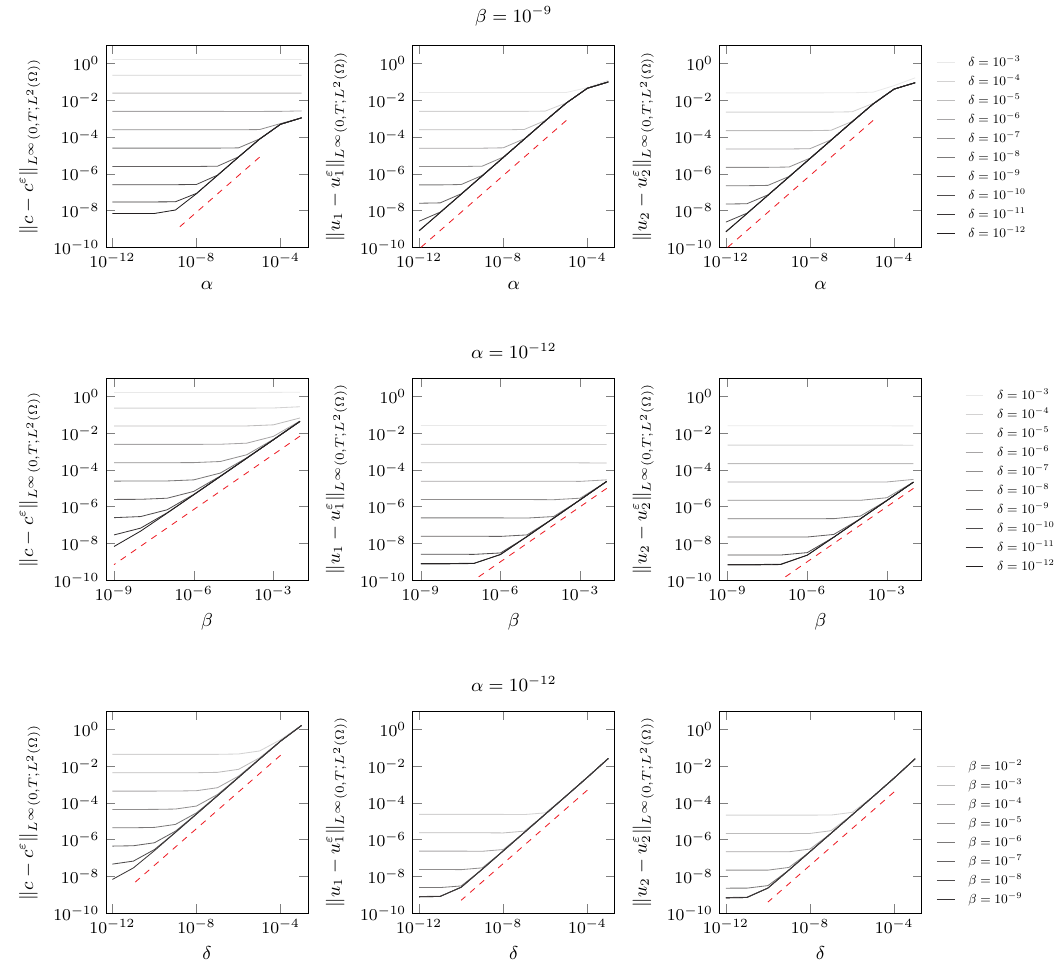}
  \caption{Convergence of the relaxation system \eqref{eq:NSCHr_num} to the NSCH system \eqref{eq:NSCH} for the merging droplets test case with $\Nx = \Ny = 50$. In the top row, the convergence in terms of $\alpha$ for a fixed value of $\beta$ and varying $\delta$ is depicted. The center row shows the convergence in terms of $\beta$ for a fixed $\alpha$ and varying $\delta$. In the bottom row, the convergence error in terms of $\delta$ for a fixed value of $\alpha$ and varying $\delta$ is plotted. The red dashed lines indicate convergence with order $\cO(\alpha)$, $\cO(\beta)$ and $\cO(\delta)$, depending on the respective running variable on the abscissa.}
  \label{fig:Error_CollidingDroplets2D_N50}
\end{figure}


\newpage
\printbibliography

\end{document}